%1.04.06
\input amstex
\input amsppt.sty
\input epsf
\NoBlackBoxes
\magnification=\magstep1
%\advance\vsize0.5cm\voffset=-1cm\advance\hsize1.5cm\hoffset-0.5cm
%arxive
\advance\vsize-1cm\voffset=-0.5cm\advance\hsize1cm\hoffset0cm

\def\diag{\mathop{\fam0 diag}}

\def\Emb{\mathop{\fam0 Emb\phantom{}}}
\def\rk{\mathop{\fam0 rk}}
\def\st{\mathop{\fam0 st}}
\def\pr{\mathop{\fam0 pr}}
\define\rel{\allowbreak\mkern 12mu {\fam0 rel}\,\,}
\def\delet{\mathaccent"7017 }
\def\rel#1{\allowbreak\mkern8mu{\fam0 rel}\,\,#1}
\def\id{\mathop{\fam0 id}}
\def\adiag{\mathop{\fam0 adiag}}

\def\Int{\mathop{\fam0 Int}}
\def\Con{\mathop{\fam0 Con}}
\def\Cl{\mathop{\fam0 Cl}}
\def\link{\mathop{\fam0 link}}
\def\lk{\mathop{\fam0 lk}}
\def\sign{\mathop{\fam0 sign}}
\def\im{\mathop{\fam0 im}}
\def\coim{\mathop{\fam0 coim}}
\def\ker{\mathop{\fam0 ker}}

\def\rmop#1{\expandafter\def\csname#1\endcsname{\operatorname{#1}}}
\def\R{\Bbb R} \def\Z{\Bbb Z} \def\C{\Bbb C} \def\H{\Bbb H}

\def\t{\widetilde}

\def\intr{1.}    \def\Intr{\S1}
\def\wi{2.}    \def\Wi{\S2}
\def\lc{3.}    \def\Lc{\S3}
\def\vk{4.}   \def\Vk{\S4}
\def\wu{5.}    \def\Wu{\S5}
\def\kt{6.}    \def\Kt{\S6}
\def\ex{7.}    \def\Ex{\S7}
\def\web{8.}    \def\Web{\S8}

\topmatter
\title Embedding and knotting of manifolds \\
in Euclidean spaces \endtitle
\author Arkadiy  Skopenkov \endauthor
\address 
Department of Differential Geometry, Faculty of Mechanics and Mathematics, 
Moscow State University, Moscow, Russia 119992, and Independent University of 
Moscow, B. Vlasyevskiy, 11, 119002, Moscow, Russia. 
e-mail: skopenko\@mccme.ru
\endaddress
\thanks
This paper is to be published in London Mathematical Society Lecture Notes,
{\it Surveys in geometry and number theory: Reports on contemporary Russian
mathematics.}
The author gratefully acknowledges the support of a grant from the London
Mathematical Society via the programme \lq Invitation of young
Russian mathematicians',
by INTAS Grant No. YSF-2002-393, by the Russian Foundation for
Basic Research, Grants No 05-01-00993 and 04-01-00682, President of Russian
Federation Grants MD-3938.2005.1 and NSH-1988.2003.1, and by the Pierre Deligne
fund based on his 2004 Balzan prize in mathematics.
\endthanks

\abstract
A clear understanding of topology of higher-dimensional objects is important in
many branches of both pure and applied mathematics.
%There are several books (or chapters in books) on two- and three-dimensional
%intuitive topology.
%Following the pattern of the Kirby and the Freedman-Quinn books on 4-manifolds,
In this survey we attempt to present some results of higher-dimensional
topology in a way which makes clear the visual and intuitive part of the
constructions and the arguments.
In particular, we show how abstract algebraic constructions
%(cohomology, characteristic classes and equivariant maps)
appear naturally in the study of geometric problems.
Before giving a general construction, we illustrate the main ideas in simple
but important particular cases, in which the essence is not veiled by
technicalities.

More specifically, we present several classical and modern results on
the embedding and knotting of manifolds in Euclidean space.
We state many concrete results (in particular, recent explicit
classification of knotted tori).
Their statements (but not proofs!) are simple and accessible to
non-specialists.
We outline a general approach to embeddings via the classical van
Kampen-Shapiro-Wu-Haefliger-Weber 'deleted product' obstruction.
This approach reduces the isotopy classification of embeddings to the homotopy
classification of equivariant maps, and so implies the above concrete results.
We describe the revival of interest in this beautiful branch of topology, by
presenting new results in this area (of Freedman, Krushkal, Teichner, Segal,
Spie\D z and the author): a generalization the Haefliger-Weber embedding
theorem below the metastable dimension range and examples showing that
other analogues of this theorem are false outside the metastable dimension
range.
\endabstract
\endtopmatter

%We also state the Browder-Haefliger-Levine-Novikov embedding theorem,
%which was obtained by application of surgery technique.
%construction of the Hudson-Habegger invariant and its applications.
%By proving some of these results and constructing some examples we illustrate
%classical and modern tools of geometric topology (engulfing, the Whitney
%trick, van Kampen and Casson finger moves and their generalizations).
%We also show how the embedding theory is related to other branches of
%mathematics, in particular, homotopy topology and functional analysis.

%We prove beyond the metastable dimension the PL cases of the
%classical theorems due to Haefliger, Harris, Hirsch and Weber) on
%the deleted product criteria for embeddings and immersions.
%The isotopy and regular homotopy versions of the above theorems
%are also improved.
%We show by examples that they cannot be improved further.
%These results have many interesting corollaries, e.g.
%1) Any closed homologically 2-connected smooth 7-manifold smoothly
%embeds in $\R^{11}$.
%2) If $p\le q$ and $m\ge\frac{3q}2+p+2$ then the set of PL
%embeddings $S^p\times S^q\to\R^m$ up to PL  isotopy is in 1--1
%correspondence with $\pi_q(V_{m-q,p+1})\oplus\pi_p(V_{m-p,q+1})$.

\document
%\head I. Obstruction theory for beginners \endhead
%Introduction
%1. Characteristic classes for beginners (orientability, vector
%fields and immersions).
%3. Stability of self-intersections of paths and cycles in the plane
%\hcm\ Homotopy classification of maps
%\head I. Obstructions to embedding and isotopy \endhead
%\pro\ Embeddings into the plane: prohibited examples

\intr\ Introduction

\wi\ The simplest-to-state results on embeddings

\lc\ Links and knotted tori

\vk\ The van Kampen obstruction

\wu\ The Haefliger-Wu invariant

\kt\ On the deleted product of the torus

%\com\ Complements
%\nb\ Neighborhoods and classification of immersions
%\lb\ The Browder-Haefliger-Levine-Novikov embedding theorem.
%\head II. Proofs of embedding theorems  \endhead

\ex\ Borromean rings and the Haefliger-Wu invariant

%\eng\ Engulfing

\web\ The disjunction method

%\gen\ Generalization of the Weber Theorem
%\cla\ Modification of an almost embedding to an embedding
%\pla\ Modification of an embedding to an immersion
%\dis\ The Disjunction Theorem
%\cyl\ Reduction of knotting theorems to embedding theorems
%\imr\ Construction of an immersion
%\head III. Related problems of geometric topology \endhead
%\apr\ Approximability by embeddings
%\bas\ Basic embeddings
%\lm\ Link maps
%\smallskip

References

%\newpage
%\head I. Obstructions to embeddability and isotopy \endhead
\head \intr\ Introduction \endhead \subhead Embedding and Knotting
Problems \endsubhead Many theorems in mathematics state that an
arbitrary object of a given abstractly defined class is a
subobject of a certain \lq standard' object of this class. Such
are the Cayley theorem on the embedding of finite groups into the
symmetric groups, the theorem on the representation of compact Lie
groups as subgroups of $GL(V)$ for a certain linear space $V$, the
Urysohn theorem on the embedding of normal spaces with countable
basis into the Hilbert space, the general position theorem on the
embedding of finite polyhedra into $\R^m$, the
Menger--N\"obeling--Pontryagin theorem on the embedding of
finite-dimensional compact spaces into $\R^m$, the Whitney theorem
on the embedding of smooth manifolds into $\R^m$, the Nash theorem
on the embedding of Riemannian manifolds into $\R^m$, the Gromov
theorem on the embedding of symplectic manifolds into $\R^{2n}$,
etc. The solution of the 13th Hilbert problem by Kolmogorov and
Arnold can also be formulated in terms of embeddings. Although
interesting in themselves, these embeddability theorems also prove
to be powerful tools for solving other problems. Subtler problem
about the embeddability of a given space into $\R^m$ for a {\it
given $m$}, as well as about counting such embeddings, are among
the most important classical problems of topology.

According to Zeeman, the classical problems of topology are the following.

1) {\it The Homeomorphism Problem:} When are two given spaces homeomorphic?

2) {\it The Embedding Problem:} When does a given space embed into $\R^m$?

3) {\it The Knotting Problem:} When are two given embeddings isotopic?

Definitions of an embedding and an isotopy are recalled in the next subsection.

Embedding and Knotting Problems have played an outstanding role in the
development of topology.
Various methods for the investigation of these problems were created by such
classical figures as G.~Alexander, H.~Hopf, E.~van Kampen, K.~Kuratowski,
S.~MacLane, L.S.~Pontryagin, R.~Thom, H.~Whitney, W. Browder, A. Haefliger,
M. Hirsch, J. F. P. Hudson, M. Irwin, J. Levine, S. Novikov, R. Penrose,
J. H. C. Whitehead, C. Zeeman and others.
For surveys see [Wu65, Introduction, Gi71, Ad93, RS96, RS99].
Nowadays interest in this subject is reviving.

There are only a few cases in which a concrete answer to the Embedding and
Knotting Problems can be given.
E. g. for the best known specific case of the Knotting Problem, i.e. for the
theory of codimension 2 embeddings (in particular, for the classical theory of knots
in $\R^3$), a complete concrete classification is neither known nor expected.
%In this survey we consider mostly the case of embeddings in codimension
%greater than 2 (i.e. $m-\dim N>2$).
A concrete complete description of a (non-empty) set of embeddings
of a given manifold up to isotopy is {\it only} known either just
below the stable range, or for highly-connected manifolds, or for
links and knotted tori. (An concrete complete answer to the
Embedding Problem was obtained additionally for projective spaces
[Gi71], for products of low-dimensional manifolds or of graphs
[ARS01, Sk03] and for certain twisted products [Re71, RS02].)
Therefore Knotting Problem is one of the hardest problems in
topology. The Embedding Problem is also hard for similar reasons.
However, the statements (but not the proofs!) are simple and
accessible to non-specialists. One of the purposes of this survey
is to list such statements. They are presented in \Wi\ and \Lc.
Statements analogous to those presented (e.g. for non-closed
manifolds) are often omitted.
%, statements which  more complicated (and thus often less useful) statements

Another purpose of this survey is to outline a general approach
useful for obtaining such concrete complete results.
There are interesting approaches giving nice {\it theoretical} results.
The author is grateful to M. Weiss for indicating that the approach of [GW99, We]
gives also concrete results on homotopy type of the space of embeddings
$S^1\to\R^n$.
The application of surgery [Le65, Br68, Ha62, Ha66, Ha66', Ha86, CRS04]
gives good concrete results for simplest manifolds.
The advantage of the surgery approach (comparatively to the deleted product approach,
see below) is that it works in the presence of smooth knots $S^n\to\R^m$.
The disadvantage is that it uses the homotopy type of the complement and description
of possible normal bundles, and so faces computational difficulties even for relatively
simple manifolds like tori, see \Lc\ (for a successful attempt to overcome this
problem see [KS05, Sk06]).
According to Wall, surgery only reduces geometric problems on embeddings
to algebraic problems which are even harder to solve [Wa70].

The method of the {\it Haefliger-Wu invariant} (or the {\it deleted product}
method) gives many concrete results.
We introduce this invariant in \Wu.
In particular, most of the results from \Wi, \Lc\ and \Vk\ can be deduced by
the deleted product method (although originally some of them were proved
directly, sometimes in a weaker form).
The deleted product method is a demonstration of a general
mathematical idea of `complements of diagonals' and the `Gauss mapping' which
appeared in works of Borsuk and Lefschetz around 1930.
The Haefliger-Wu invariant generalizes the linking coefficient, the Whitney
obstruction and the van Kampen obstruction.
The deleted product method in the  theory of embeddings was developed by van
Kampen (1932), Shapiro (1957), Wu (1957-59), Haefliger (1962), Weber (1967),
Harris (1969) and others.
The Van-Kampen-Shaprio-Wu approach works for embeddings of {\it polyhedra},
but is closely related to embeddings of manifolds and so is presented in \Vk.
The classical Haefliger--Weber Theorem \wu4 asserts the bijectivity of the
Haefliger-Wu invariant for embeddings of $n$-dimensional polyhedra and
manifolds into $\R^m$ under the \lq metastable' dimension restriction
$$2m\ge3n+4.$$
Other embedding invariants (obtained using $p$-fold deleted products, see
end of \Wu, or using complement together with normal bundle [Br68, CRS04])
are hard to compute.
So the investigation of embeddings for $2m<3n+4$ naturally leads to the problem
of finding conditions under which the Haefliger-Wu invariant is complete
without the metastable dimension assumption.
There were many examples showing that for embeddings of manifolds the
metastable dimension restriction is sharp in many senses (Boechat, Freedman,
Haefliger, Hsiang, Krushkal, Levine, Mardesic, Segal, Skopenkov, Spiez,
Szarba, Teichner, Zeeman, see \Wu).
So it is surprising (Theorem \wu5) that for $d$-connected $n$-dimensional
manifolds and in the piecewise linear category the metastable restriction in
the Haefliger-Weber Theorem can be weakened to
$$2m\ge3n+3-d.$$
We present many beautiful {\it examples} motivated by the Embedding and
Isotopy Problems.
In particular, in \Wi\ we present a construction of the Hudson torus (which is
simpler than the original one), in \Lc\ we construct examples illustrating the
distinction between piecewise linear  and smooth embeddings.
In \Kt\ we prove some results on the deleted product of 'torus' $S^p\times S^q$
and on the Haefliger-Wu invariant of knotted tori.
In \Kt\ and \Ex\ we construct most of the examples of the incompleteness of the van
Kampen obstruction and the Haefliger-Wu invariant, announced in \Wu\ and \Vk.
The construction of these examples is based on knotted tori (\Kt) or on
(higher-dimensional) Casson finger moves (\Ex).
For some other examples we only give references.

The Haefliger-Weber theorem and its analogue below the metastable case
were obtained by a combination of and the improvement of methods and results
from various parts of topology: the theory of immersions, homotopy theory,
engulfing, the Whitney trick, van Kampen finger moves, the
Freedman-Krushkal-Teichner trick and their generalizations.
The most important method is the {\it disjunction method} (end of \Vk\ and
\Web).
These methods are also applied in other areas.
In \Web\ we prove the surjectivity of the Haefliger-Wu invariant in the piecewise
linear case.
For the reader's convenience, we take a historical approach to the exposition:
the disjunction method is applied in its complete generality only after
illustration in simpler particular cases.
We also prove the analogue of the Haefliger-Weber theorem below the
metastable range for the simplest case.
We do not prove many other results of \S2-\S5 but give references and sometimes
sketch the proofs.

Sections \Kt, \Ex\ and \Web\ depend on \Wu.
Otherwise the sections are independent of each other except for minor details
that can well be omitted during the first reading.

\subhead Definitions and notations \endsubhead
A smooth {\it embedding} is a smooth injective map $f:N\to\R^m$ such that
$df$ is a monomorphism at each point.

By a polyhedron we shall understand a {\it compact} polyhedron.
A map $f:N\to\R^m$ of a polyhedron $N$ is {\it piecewise-smooth} if it is
smooth on each simplex of some smooth triangulation of $N$.
We denote the piecewise-smooth category by PL.
This is the usual notation for the piecewise-linear category but the
classification of piecewise-smooth embeddings (or immersions) coincides with
the classification of piecewise linear embeddings (or immersions) [Ha67].
A PL {\it embedding} is a PL injective map $f:N\to\R^m$.

We write CAT for DIFF or PL.
We often omit CAT if a statement holds in both PL and DIFF categories.

Two embeddings $f,g:N\to\R^m$ are said to be {\it (ambient) isotopic}
(Figure \intr1), if there exists a homeomorphism onto (an {\it ambient
isotopy}) $F:\R^m\times I\to\R^m\times I$ such that

$F(y,0)=(y,0)\quad\text{for each}\quad y\in\R^m,$

$F(f(x),1)=g(x)\quad\text{for each}\quad x\in N,$ \quad and

$F(\R^m\times\{t\})=\R^m\times\{t\}\quad\text{for each}\quad t\in I.$

\bigskip
\centerline{\epsffile{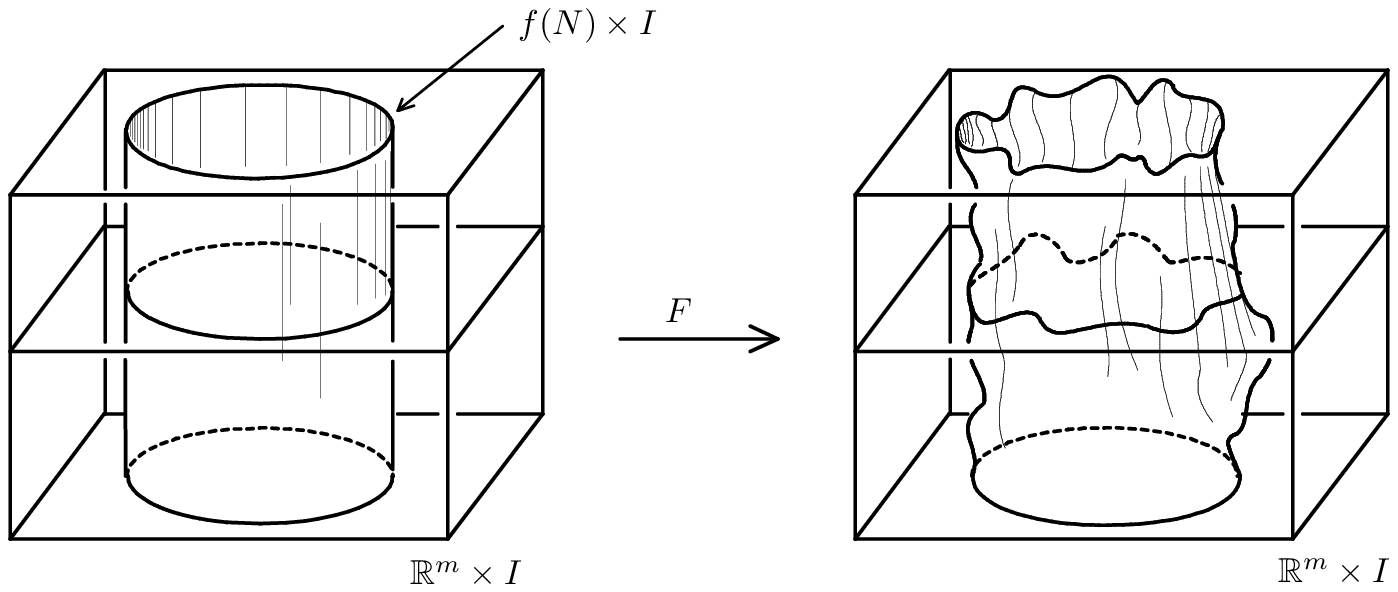}}
\centerline{\it Figure \intr1} %approximately here}
\bigskip

An ambient isotopy is also a homotopy $\R^m\times I\to\R^m$ or a family of
maps $F_t:\R^m\to\R^m$, generated by the map $F$ in the obvious manner.
Evidently, ambient isotopy is an equivalence relation on the set of
embeddings of $N$ into $\R^m$.
An embedding is called {\it trivial}, if it is isotopic to the standard
embedding (the latter is evidently defined from the context).

We use the notations of [RS72].
An equality between sets denotes a 1--1 correspondence.
Denote by $V_{m,n}$ the Stiefel manifold of $n$-frames in $\R^m$.
Let $\Z_{(k)}$ be $\Z$ for $k$ even and $\Z_2$ for $k$ odd.
Note that $\Z_{(k)}=\pi_k(V_{n,n-k})$ for $1<k<n$.
If the coefficients are omitted from the notation of (co)homology groups,
then $\Z$-coefficients are assumed.
For a manifold or a polyhedron $N$ we denote its dimension by $n=\dim N$.
Denote by $\Emb^m_{CAT}(N)$ the set of CAT embeddings $N\to\R^m$ up to ambient
CAT isotopy.
If $|\Emb^m_{CAT}(N)|=1$, we shall say that $N$ {\it CAT unknots} in $\R^m$.
An embedding is {\it trivial} if it is isotopic to the standard embedding
(whose choice is clear).
For a map $f:N\to\R^m$ we denote by
%$$\t\Delta(f)=\{(x,y)\in N\times N\ |\ x\ne y,\ fx=fy\}\quad\text{and}\quad
$\Sigma(f)=\Cl\{x\in N\ :\ |f^{-1}fx|\ge2\}$ its {\it self-intersection set}.

A closed manifold $N$ is called {\it homologically $k$-connected}, if $N$ is
connected and $H_i(N)=0$ for each $i=1,\dots,k$.
This condition is equivalent to $\tilde H_i(N)=0$ for
each $i=0,\dots,k$, where $\tilde H_i$ are reduced homology groups.
A pair $(N,\partial N)$ is called {\it homologically $k$-connected}, if
$H_i(N,\partial N)=0$ for each $i=0,\dots,k$.
Note that if $H_0(N,\partial N)=0$, then the manifold $N$ has no closed connected
components.
We use
the following conventions: 0-connectedness is equivalent to homological
0-connectedness and to connectedness; $k$-connectedness for $k<0$ is an empty condition.
We put $\pi^S_l=0$ for $l<0$.

\subhead Other equivalence relations and problems \endsubhead
Ambient isotopy is a stronger equivalence relation than non-ambient isotopy,
isoposition, concordance, bordism, etc.
Two embeddings $f,g:N\to\R^m$ are called {\it (non-ambient) isotopic}, if
there exists an embedding $F:N\times I\to\R^m\times I$ such that

$F(x,0)=(f(x),0)$,

$F(x,1)=(g(x),0)$ for each $x\in N$ and

$F(N\times\{t\})\subset\R^m\times\{t\}$ for each $t\in I$.

In the DIFF category or for $m-n\ge3$ in the PL (or TOP) category
{\it isotopy implies ambient isotopy} [HZ64, Hu66, Ak69, Ed75, \S7].
For $m-n\le2$ this is not so: e.g., any knot $S^1\to S^3$ is non-ambiently PL
isotopic to the trivial one, but not necessarily ambiently PL isotopic to it.

Two embeddings $f,g:N\to\R^m$ are said to be (orientation preserving)
{\it isopositioned}, if there is an (orientation preserving) homeomorphism
$h:\R^m\to\R ^m$ such that $h\circ f=g$.
For embeddings into $\R^m$
{\it orientation preserving isoposition is equivalent to isotopy}
(the Alexander-Guggenheim Theorem) [RS72, 3.22].
%smooth case?

Two embeddings $f,g:N\to\R ^m$ are said to be {\it ambiently concordant}
if there is a homeomorphism (onto) $F:\R^m\times I\to\R ^m\times I$ (which is
called a {\it concordance}) such that

$F(y,0)=(y,0)$ for each $y\in\R^m$ and

$F(f(x),1)=g(x)$ for each $x\in N$.

The definition of {\it non-ambient concordance} is analogously obtained from
that of non-ambient isotopy by dropping the last condition of
level-preservation.
In the DIFF category or for $m-n\ge3$ in the PL (or TOP) category
{\it concordance implies ambient concordance and isotopy} [Li65, Hu70, HL71]
(this is not so in codimension 2).
This result allows a reduction of the problem of isotopy to the relativized
problem of embeddability.

Let us give a (by no means complete) list of references for closely related
problems in geometric topology (the references inside the papers listed here
could also be useful for a reader).
In the problems of embeddability and isotopy the space $\R^m$ can be replaced
by an {\it arbitrary} space~$Y$.
The case when $Y$ is a manifold has been studied most extensively; for the
case when $Y$ is a product of trees see [St89, Theorem 4.6 and Remark, GR92,
GMR94, Zh94, Ku00].
For embeddings {\it up to cobordism} see [Br71, Li75].
For embeddings {\it up to homotopy} see [Co69, St, Wa70, \S11, Hu70', CW78, Ha84].
For the classification of {\it link maps} see [Mi54, MR86, Ko88, Ko90, Ma90, HK98, Sk00].
For embeddings of polyhedra in {\it some} manifolds see [Wa66, Ne68, LS69, RBS99].
For the problem of {\it embeddability of compacta} and the close problem of
{\it approximability by embeddings} see [Ch69, Mc67, Si69, SS83, KW85,
Da86, Ak96, Ak96', RS96, \S9, Mi97, RS98, Ak00, RS01', ARS02, Me02, RS02, Sk03',
Me04] (the author is grateful to P. Akhmetiev for indicating that
the paper [Ak96] contains a mistake for $n=3,7$ and that the paper [Ak96']
contains a 'preliminary version' of the proof, the complete version being
submitted elsewhere).
For the problem of {\it intersection of compacta} see [DRS93, ST91].
For {\it basic embeddings} see [St89, Sk95, RS99, \S5, Ku00].
For {\it immersions} see [Gi71, Ad93, cf. Sk02].

\subhead Acknowledgements \endsubhead
This survey was based on lectures the author had given at various times at
the Independent University of Moscow, Moscow State University, the Steklov
Mathematical Institute (Moscow and St.\ Petersburg branches), the Technical
University of Berlin, the Ruhr University of Bochum, the Lorand E\"otvos
University of Budapest,  the University of Geneva, the University of Heidelberg,
the University of Ljubljana, the University of Siegen, the University of Uppsala,
the University of Warsaw and the University of Zagreb.
The preliminary version was
prepared in January 2002 after a series of lectures at the Universities of
Aberdeen, Cambridge, Edinburgh and Manchester, sponsored by London
Mathematical Society.  I would like to acknowledge all these institutions
for their hospitality and, personally, P. M. Akhmetiev, V. M.  Buchstaber, A. V.
Chernavskiy, P.  Eccles, K. E. Feldman, A. T. Fomenko, M.  Kreck, U.
Koschorke, R.  Lickorish, A.  Haefliger, R. Levy, S. Mardesic, A.  S.
Mischenko, N. Yu.  Netsvetaev, V. M.  Nezhinskiy, M. M. Postnikov, E.
Rees, D. Repovs, E. V.  Schepin, Yu. P.  Solovyov, A. Sz\"ucs, V. A.
Vassiliev, O. Ya. Viro, C.  Weber, M. Weiss, G.  Ziegler and H.  Zieschang
for their invitations and useful discussions.
It is a pleasure to express special gratitude to P. Eccles for his
many remarks on the preliminary version of this paper.

%\newpage
\head \wi\ The simplest-to-state results on embeddings \endhead

\subhead Embeddings just below the stable range \endsubhead

\proclaim{General Position Theorem \wi1.a} Every $n$-polyhedron or $n$-manifold
embeds into $\R^{2n+1}$ and unknots in $\R^m$ for $m\ge2n+2$.
\endproclaim

\bigskip
\centerline{\epsffile{3-1.eps}}
\centerline{\it Figure \wi1} %approximately here}
\bigskip

The restriction $m\ge2n+2$ in Theorem \wi1.a is sharp as the Hopf linking
$S^n\sqcup S^n\to \R ^{2n+1}$ shows (Figure~\wi1).
The number $2n+1$ in Theorem \wi1.a is the minimal possible for polyhedra.

\smallskip
{\bf Example \wi1.b.} {\it For each $n$ there exists an $n$-polyhedron,
non-embeddable in $\R^{2n}$.}

\smallskip
In Example \wi1.b one can take the $n$-th power of a non-planar graph [Sk03],
the $n$-skeleton of a $(2n+2)$-simplex [Ka32, Fl34] or the $(n+1)$-th join power
of the three-point set (see the proof for this case in \Wu).

Note that

$N\times I$ embeds into $\R ^{2n+1}$ for each $n$-polyhedron $N$ [RSS95];

$N\times I$ unknots in $\R^{2n+2}$ for each $n$-polyhedron $N$ (let us sketch a proof
which can though be omitted for the first reading: the result follows because by general
position every two embeddings $N\times I\to\R^{2n+2}$ are regular homotopic and their
restrictions to the spine $N\times\{\frac12\}$ are isotopic).

\smallskip
{\bf Theorem \wi2.a.} {\it Every $n$-manifold embeds into $\R^{2n}$}
[Ka32, Wh44].

\smallskip
Theorem \wi2.a (as well as Theorem \wi2.b below) is proved using general
position and the Whitney trick; the proof in the smooth and PL case is
sketched in \Vk\ and in [RS72, RS99, \S8], respectively.

The dimension $2n$ in Theorem \wi2.a is the best possible when $n=2^k$ because
$\R P^{2^k}$ does not embed into $\R ^{2^{k+1}-1}$
(this is proved using the mod2 Whitney obstruction defined below [MS74, RS00,
RS02']).
But is not the best possible for other $n$ by Theorem \wi3 below.
A celebrated and difficult conjecture is that every closed $n$-manifold
embeds into $\R ^{2n+1-\alpha(n)}$, where $\alpha(n)$ is the number of units
in the dyadic expansion of $n$.
The analogous conjecture for immersions was proved in [La82, Co85].
Note that if $n=2^{k_1}+\dots+2^{k_{\alpha(n)}}$ and $k_1<\dots<k_{\alpha(n)}$,
then the $n$-manifold
$N=\R P^{2^{k_1}}\times\dots\times\R P^{2^{k_{\alpha(n)}}}$
does not embed into $\R^{2n-\alpha(n)}$
(this is proved again using the mod2 Whitney obstruction defined below).

\smallskip
{\bf Theorem \wi3.}
{\it (a) Every $n$-manifold (except that if $n=2^k$ and the manifold is closed,
we need to assume that it is orientable) embeds into $\R^{2n-1}$.

(b) Suppose that $N$ is a closed $n$-manifold, $n$ is even,
$n\ne2^k(2^h +1)$ for integers $k,h\ge2$ and $H_1(N)=0$.
Then $N$ embeds into $\R^{2n-2}$, provided $n\ge6$ in the PL category or
$n\ge8$ in the smooth category.

(c) Suppose that $N$ is a closed $n$-manifold, $\alpha(n)\ge5$ and
either $n=0,1(4)$ and $N$ is orientable, or $n=2,3(4)$ and $N$ is
non-orientable.
Then $N$ embeds into $\R^{2n-2}$.}

\smallskip
Classical Theorems \wi3 are much more complicated to prove than Theorem \wi2.a:
one needs a generalization of the Whitney trick and calculation of characteristic
classes, both non-trivial.
Theorems \wi3.ab follow from the case $d=0$ of Theorem \wi8.a below (which is true
for orientable manifolds) and [Ma60, Theorem 1.c and
Corollary 2, Ma62, Theorem 1], except that Theorem \wi3.a for $n=3,4$ has to
be proved separately [Hi61, Hi65, Ro65, Wa65, BH70, Do87, Fu94, cf.\ No61,
Fu02].
(From the proof [Ma60, p.~100] it follows that the restriction [Ma60, Theorem
1.c] should be stated as \lq the number of $h_i$'s which are equal to $h_q+1$
is even', cf. [Hi61].)
Theorem \wi3.c follows from the Haefliger-Weber Theorem \wu4 below and [Ba75,
Theorem 45].

The condition \lq $\alpha(n)\ge5$' in Theorem \wi3.c can be weakened to
\lq $n\ge7$ and $\bar w_{n-i}(N)=0$ for each $i\le4$' (the classes
$\bar w_i(N)$ are defined in the subsection 'the Whitney obstruction' below).

\smallskip
{\bf Theorem \wi2.b.}
{\it Every connected $n$-manifold unknots in $\R^{2n+1}$ for $n>1$} [Wu58].

\smallskip
Here for each $n$ the dimension $2n+1$ is the best possible and the
connectedness assumption is indeed necessary, as the Hopf linking above and the
Hudson Torus Example \wi6.c below show.
The Hopf linking is distinguished from the standard linking using the simplest
($\Z$-valued or $\Z_2$-valued) {\it linking coefficient} (whose definition is
obtained from the definition of the Whitney invariant by setting $m=2n+1$ in
the subsection 'the Whitney invariant' below).

\smallskip
{\bf Theorem \wi4.}
{\it Suppose that $n\ge2$ and a compact $n$-manifold $N$ has $s$
closed orientable connected components and $t$ closed non-orientable connected
components (and, possibly, some non-closed components).
Then the set of pairwise linking coefficients defines a 1--1 correspondence
$\Emb\phantom{}^{2n+1}(N)\to\Z^{\frac{s(s-1)}2}\oplus\Z_2^{st+\frac{t(t-1)}2}$.}

\smallskip
Note that {\it every $n$-polyhedron $N$ such that $H^n(N)=0$

(a) PL (if $n=2$, only TOP) embeds into $\R^{2n}$;

(b) PL unknots in $\R^{2n+1}$.}

Assertion (a) for $n\ge3$ is deduced from Theorems \vk1 or \wu4 below [Wu58,
We68, see also Ho71].
For $n=2$ it was proved independently [Ki84] and for $n=1$ it is trivial.
Assertion (b) for $n\ge2$ is deduced from Theorems \vk4 or \wu4 below [Pr66]
and for $n=1$ it is trivial.
In these assertions for each $n$ the dimensions are the best possible and the
$H^n(N)=0$ assumption is indeed necessary.

%\newpage
\subhead Embedding and unknotting of highly-connected manifolds \endsubhead

\smallskip
{\bf Theorem \wi5.} {\it The sphere $S^n$, or even any homology $n$-sphere,

(a) PL unknots in $\R^m$ for $m-n\ge3$ [Ze60, St63, Gl63, Sc77];

(b) DIFF unknots in $\R ^m$ for $m\ge\frac{3n}2+2$ [Ha61, Ha62'', Ad93, \S7];

(c) PL (if $n=3$, only TOP) embeds into $\R^{n+1}$ [follows from Ke69] and

(d) DIFF embeds into $\R ^{[3n/2]+2}$ [Ha61, Ha62'', Ad93, \S7].}

\smallskip
Theorem \wi5.a is also true in the TOP locally flat category (see the definition in
\Lc) [Ru73, Sc77].
Here the local flatness assumption is indeed necessary.

Knots in codimension 2 and the Trefoil Knot Example \lc5 below show that the
dimension restrictions are sharp (even for standard spheres) in Theorems
\wi5.a and \wi5.b, respectively.
By [Le65, HLS65, cf. MT95, pp.\ 407--408] the dimension restriction in
Theorem \wi5.d is indeed necessary (and conjecturally almost sharp) even
for {\it homotopy spheres}.
However, from [Bo71], it follows that

{\it any $4k$-dimensional homotopy sphere embeds into $\R^{6k+1}$.}

Theorems \wi2.b and \wi5 may be generalized as follows.

\smallskip
{\bf Theorem \wi6.} {\it For $n\ge2d+2$, every closed homologically
$d$-connected $n$-manifold

(a) embeds into $\R^{2n-d}$ ($n\ne2d+2$ in the DIFF case) and

(b) unknots in $\R^{2n-d+1}$. }

\smallskip
Theorem \wi6 was proved directly in [PWZ61, Ha61, Ze62, Ir65, Hu69] for
{\it homotopically} $d$-connected manifolds and using the deleted product
method (\Wu) in [Ha62'', We67, Ad93, \S7, Sk97, Sk02] for {\it homologically}
$d$-connected manifolds.
Theorem \wi6 follows from Theorem \wi8 below.
%The proof of PL case of Theorem \wi6.a for {\it homotopically} $d$-connected
%manifolds is given in \Eng.
Note that if $n\le2d+1$, then every closed homologically
$d$-connected $n$-manifold is a homology sphere, so the PL case of
Theorem \wi6 gives nothing more than Theorem \wi5.

For generalizations of Theorem \wi6 see Theorem \wi8 below or
[Hu67, Ha68', Hu72, Go72, Ke79].
We shall use of them (the simplest case $2m\ge3n+3$ of) the following
relative version of Theorem \wi6.a.

\smallskip
{\bf The Penrose-Whitehead-Zeeman-Irwin Embedding Theorem \wi6.c.}
{\it If $m-n\ge3$, then any proper map from a $(2n-m)$-connected
PL $n$-manifold with boundary to a $(2n-m+1)$-connected PL
$m$-manifold with boundary, whose restriction to the boundary is an embedding,
is homotopic (relatively to the boundary) to a PL embedding} [PWZ61, Ir65].

\smallskip
The dimension assumption in Theorem \wi6.b is sharp:

%\lc7
\smallskip
{\bf The Hudson Torus Example \wi6.c.} {\it For each $p\le q$ there exists a
non-trivial embedding $S^p\times S^q\to\R^{p+2q+1}$} [Hu63].

\smallskip
{\it Simplified construction [Sk06].} (This construction is
interesting even for $p=q=1$!) Take the standard embedding
$2D^{p+q+1}\times S^q\subset\R^{p+2q+1}$. The Hudson torus is the
(linked!) connected sum of the $(p+q)$-sphere $2\partial
D^{p+q+1}\times x$ with the standard embedding $\partial
D^{p+1}\times S^q\subset D^{p+q+1}\times S^q\subset\R^{p+2q+1}$.

\smallskip
The Hudson torus can be distinguished from the standard embedding using the
{\it Whitney invariant} defined in the subsection under the same name below
[Sk06] or the {\it Haefliger-Wu invariant} defined in \Wu\ [Sk02].

The rest of this subsection can be omitted for the first reading.

\bigskip
\centerline{\epsffile{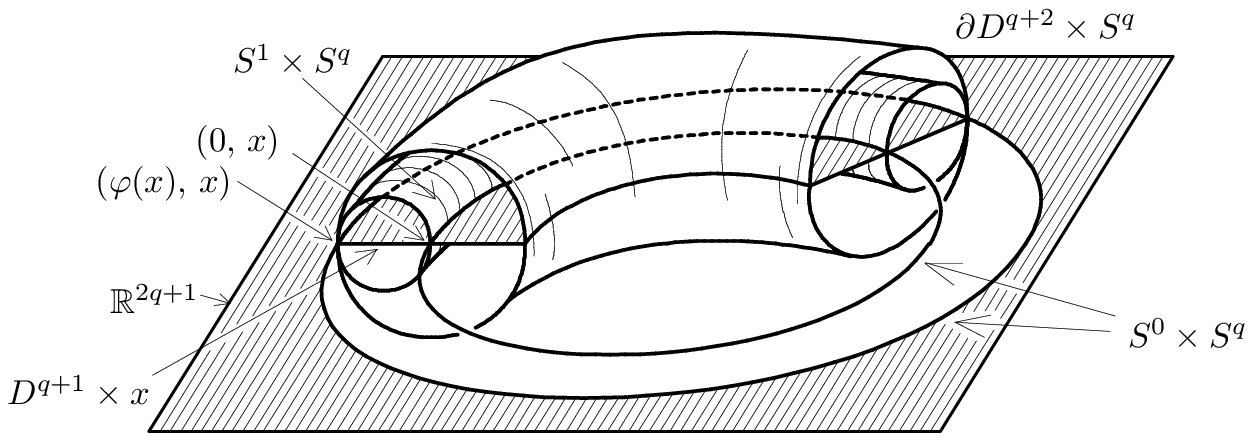}}
\centerline{\it Figure \wi2.a} %approximately here}
\bigskip
\centerline{\epsffile{3-11.eps}}
\centerline{\it Figure \wi2.b} %approximately here}
\bigskip

{\it An alternative simplified construction of the Hudson torus
[Sk02, Sk].} Define a map $S^0\times S^q\to D^{q+1}$ to be the
constant $0\in D^{q+1}$ on one component and the standard
embedding $\varphi$ on the other component. This map gives an {\it
embedding}
$$S^0\times S^q\to D^{q+1}\times S^q\subset D^{q+2}\times S^q\subset\R^{2q+2}$$
(Figure \wi2.a).
Each disk $D^{q+2}\times x$ intersects the image of this embedding at two
points lying in $D^{q+1}\times x$.
Extend this embedding $S^0\to D^{q+1}\times x$ for each $x\in S^q$ to an embedding
$S^1\to D^{q+2}\times x$ (Figure \wi2.b).
Thus we obtain a non-standard embedding
$$S^1\times S^q\to D^{q+2}\times S^q\subset\R^{2q+2}.$$
Taking spheres of dimensions $p\le q$ we obtain analogously
an embedding
$$S^p\times S^q\to D^{q+p+1}\times S^q\subset\R^{2q+p+1}.$$
(Taking as $\varphi$ above an arbitrary CAT-map $S^q\to\partial D^{m-q-p}$ we
obtain analogously a CAT embedding
$S^p\times S^q\overset{\varphi\times\pr_2}\to\to D^{m-q}\times S^q\to\R^m$.)

\smallskip
It would be interesting to know whether the smooth case of Theorem \wi6.a
holds for $n=2d+2$, i.e. for  $d$-connected $2(d+1)$-manifolds.
An {\it almost smooth} embedding is a PL embedding which is a smooth
embedding outside a point.

\smallskip
{\bf Theorem \wi7.}
{\it Let $N$ be a closed smooth $(l-1)$-connected $2l$-manifold.

(a) If $l\equiv3,5,7\mod8$ and $l\ne2^s-1$, then $N$ almost smoothly embeds
into $\R^{2l+1}$.
%(or if $\Arf N=0$)

(b) For $l$ even the manifold $N$ almost smoothly embeds into
$\R^{2l+1}$ if and only if $N$ is stably parallelizable.

(c) If $N$ is almost parallelizable, then $N$ almost smoothly embeds
into $\R^{2l+2}$ [MRS03].

(d) For each even $l$ there exists a closed smooth $(l-1)$-connected (even
almost parallelizable) $2l$-manifold which does not smoothly embed into
$\R^{3l-1}$ [MRS03].}

\smallskip
Theorems \wi7.a and \wi7.b are proved in [Mi65, Corollary 1, Theorem 2 and
Addendum (1)], cf. [Sa65].
By Theorem \wi7.c, the manifold from \wi7.d almost smoothly embeds into
$\R^{2l+2}$.

\smallskip
{\it Proof of \wi7.c.}
If $l=2$, then the result holds by [Ma78, Corollary 10.11, CS79].
So assume that $l\ne2$.
Let $N_0$ be a complement in $N$ to some open $2l$-ball.
Then $N_0$ is parallelizable and hence there is an immersion
$f:N_0\to\R^{2l+1}$.
Since $N_0$ is $(l-1)$-connected and $l\ne2$, it follows that it has an
$l$--dimensional spine [Wa64, Ho69].
By general position the restriction of $f$ to this spine is an embedding.
Hence the restriction of $f$ to some neighborhood of this spine is an
embedding.
But this neighborhood is homeomorphic to $N_0$.
So there is an embedding $g:N_0\to\R^{2l+1}$.
Extending the embedding $g|_{\partial N_0}$ as a cone in $\R^{2l+2}$ we
obtain an almost smooth embedding of $N$ into $\R^{2l+2}$.
\qed

\smallskip
{\it Proof of \wi7.d.}
Take the Kervaire-Milnor closed smooth almost parallelizable
$4k$-manifold $N$ whose signature $\sigma(N)\ne0$ [Ma80, MK58].
We can modify this by surgery [Ma80] and assume further that it is
$(l-1)$-connected.
Hence the Pontrjagin class $p_{l/2}(N,\R)\ne0$ by the Hirzebruch formula.
Therefore $\bar p_{l/2}(N,\R)\ne0$ by the duality theorem for real
Pontryagin classes [Wu65, cf. MS74].
Hence $N$ does not smoothly embed into $\R^{3l-1}$ (it does
not even immerse in $\R^{3l-1}$) [Po42, Wu65].
\qed

\smallskip
The dimension $2n-d$ in Theorem \wi6.a can be decreased by 1 for some pairs
$(n,d)$, as Theorem \wi3 shows.
However, we conjecture that the dimension $2n-d$ in Theorem \wi6.a  cannot be
significantly decreased for some $(n,d)$.
This is so for $d=0$ (as the example $N=\R P^{a_1}\times\dots\times\R P^{a_s}$
shows) and for $n=2,4,8$, $d=\frac n2-1$ (take $N=\R P^2$, $\C P^2$, $\H P^2$
or apply Theorem \wi7.d).
Example of a highly-connected but badly embeddable manifolds were also
exhibited in [HS64, Sa65].

%\newpage
\subhead The Whitney obstruction \endsubhead

{\it Definition of the modulo 2 Whitney obstruction.}
Let $N$ be a closed manifold.
We present the definition in the piecewise linear case, the definition in the
(piecewise) smooth case is analogous.
Take any general position map $f:N\to\R^m$.
Recall the definition of the self-intersection set $\Sigma(f)$ from \Intr\ (Figure
\wi3).

\bigskip
\centerline{\epsffile{2-1.eps}}
\centerline{\it Figure \wi3} %approximately here}
\bigskip

Take a triangulation $T$ of $N$ such that $f$ is linear on simplices of $T$.
Then the self-intersection set $\Sigma(f)$ is a subcomplex of $T$.
Denote by $[\Sigma(f)]\in C_{2n-m}(N;\Z_2)$ the sum of the top-dimensional
simplices of $\Sigma(f)$ (Figure \wi1).
Then $[\Sigma(f)]$ is a cycle [Hu69, Lemma 11.4, Hu70', Lemma 1].

(In order to prove this fact it suffices to prove that each $(2n-m-1)$-simplex
$\eta$ of $T$ is in the boundary of an {\it even} number of $(2n-m)$-simplices
$\sigma\subset\Sigma(f)$.
It suffices to consider the case $\eta\in\Sigma(f)$.
By general position, $f^{-1}f\eta$ consists of simplices
$\eta=\eta_1,\dots,\eta_k$.
The link $\lk_T\eta_i$ is a sphere of dimension $n-(2n-m-1)-1=m-n$.
The link $\lk_{\R^m}f\eta$ is a sphere of dimension $m-(2n-m-1)-1=2(m-n)$.
The intersection of two $f$-images $f(\lk_T\eta_i)$ of $(m-n)$-spheres
in the $2(m-n)$-sphere $\lk_{\R^m}f\eta$ consists of an even number of points.
These intersection points are in 1--1 correspondence with $(2n-m)$-simplices
$\sigma\subset\Sigma(f)$ containing $\eta$ in their boundaries.)

{\it The modulo 2 Whitney obstruction} is the homology class
$$\bar w_{m-n}(N):=[\Sigma(f)]\in H_{2n-m}(N;\Z_2).$$
The class $\bar w_i$ is called {\it the normal Stiefel-Whitney class}.
This definition of the normal Stiefel-Whitney classes is equivalent to other
definitions [MS74], up to Poincar\'e duality.

The independence of $\bar w_{m-n}(N)$ on $f$ follows from the equality
$[\Sigma(f_0)]-[\Sigma(f_1)]=\partial[\Sigma(F)]$ for a general position
homotopy $F:N\times I\to\R^m\times I$ between general position maps
$f_0,f_1:N\to\R^m$.

Hence these classes are obstructions to the embeddability of $N$ into
$\R^m$:

{\it if $N$ embeds into $\R^m$, then $\bar w_i(N)=0$ for $i\ge m-n$} [Wh35].

\smallskip
{\it Definition of the Whitney obstruction for $N$ orientable and $m-n$ odd.}
Fix in advance any orientation of $N$ and of $\R^m$.
The definition is analogous to the above, only $[\Sigma(f)]$ is the sum of
{\it oriented} simplices $\sigma$ with $\pm$ signs defined as follows.
(For $m-n$ even the signs can also be defined but are not used because
$[\Sigma(f)]$ is not necessarily a cycle with integer coefficients).

By general position there is a unique simplex $\tau$ of $T$ such that
$f(\sigma)=f(\tau)$.
The orientation on $\sigma$ induces an orientation on $f\sigma$ and then on
$\tau$.
The orientations on $\sigma$ and $\tau$ induce orientations on normal spaces in
$N$ to these simplices.
These two orientations (in this order) together with the orientation on
$f\sigma$ induce an orientation on $\R^m$.
If this orientation agrees with the fixed orientation of $\R^m$, then the
coefficient of $\sigma$ is $+1$, otherwise $-1$.
Clearly, the change of orientation of $\sigma$ changes the sign of $\sigma$ in
$[\Sigma(F)]$, so the sign is well-defined.

(An equivalent definition of the signs in $[\Sigma(f)]$ is as follows.
The orientation on $\sigma$ induces an orientation on $f\sigma$ and then on
$\tau$, hence it induces an orientation on their links.
Consider the oriented sphere $\lk_{\R^m}f\sigma$, that is the link of $f\sigma$
in certain triangulation of $\R^m$ 'compatible' with $T$.
The dimension of this sphere is $m-1-(2n-m)=2(m-n)-1$.
This sphere contains disjoint oriented $(m-n-1)$-spheres $f(\lk_T\sigma)$ and
$f(\lk_T\tau)$.
The coefficient of $\sigma$ in $[\Sigma(F)]$ is their linking coefficient,
which equals $\pm1$.)

{\it The Whitney obstruction} is the homology class
$$\bar W_{m-n}(N):=[\Sigma(f)]\in H_{2n-m}(N;\Z).$$

Clearly, $\bar w_i(N)$ are modulo 2 reductions of $\bar W_i(N)$.
The classes $\bar w_i(N)$ are easier to compute, however they are possibly
weaker than $\bar W_i$.

Pontryagin introduced for closed orientable $n$-manifold $N$ {\it the
Pontryagin classes} $\bar p_i\in H_{n-4i}(N;\Z)$ which obstruct to
embeddability into $\R^{n+2i-1}$ [Po42].

Recall the definition of $\Z_{(k)}$ from \Intr.
For a closed orientable $n$-manifold $N$ denote by
$\bar W_{m-n}(N)\in H_{2n-m}(N,\Z_{(m-n-1)})$ the class $\bar W_{m-n}(N)$ for
$m-n$ odd and the class $\bar w_{m-n}(N)$ for $m-n$ even.

\smallskip
{\bf Theorem \wi8.a.}
{\it Let $N$ be a closed $d$-connected $n$-manifold, $d\ge1$.
The manifold $N$ embeds into $\R^{2n-d-1}$ if $\bar W_{n-d-1}(N)=0$, provided
$n\ge d+4$ or $n\ge2d+5$, in the PL or DIFF cases respectively.}

\smallskip
See references to proofs after Theorem \wi8.b.

In Theorem \wi8.a the $d$-connectedness assumption can be weakened to the
{\it homological} $d$-connectedness except when $n=2d+2$ in the PL case.
The PL case of Theorem \wi8.a gives nothing but Theorem \wi5.c for
$d+4\le n\le 2d+1$.
The smooth case of Theorem \wi8.a is true if $d$ is even and $n=2d+3$ [Sk02,
Corollary 1.7].

\subhead The Whitney invariant \endsubhead

\smallskip
{\it Definition of the Whitney invariant.} Let $N$ be a closed
connected orientable $n$-manifold. Let $f_0:N\to\R^m$ be a certain
fixed (\lq standard') embedding and let $f:N\to\R^m$ be an
arbitrary embedding. Take a general position homotopy $F:N\times
I\to\R^m\times I$ between $f$ and $f_0$ (Figure \wi4). 

\bigskip
\centerline{\epsffile{2-2.eps}}
\centerline{\it Figure \wi4} %approximately here}
\bigskip

Analogously
to the above, the self-intersection set $\Sigma(F)$ supports a
$(2n-m+1)$-cycle $[\Sigma(F)]$ in $N\times I\simeq N$ with the
coefficients $\Z_{(m-n-1)}$. The {\it Whitney invariant} of $f$ is
the homology class of this cycle:
$$W(f):=[\Sigma(F)]\in H_{2n-m+1}(N,\Z_{(m-n-1)}).$$
Analogously to the above, $W(f)$ depends only on $f$ and $f_0$ but not on the
choice of $F$ [Hu69, \S11, cf. HH63, Vr77, Sk06].

\proclaim{Theorem \wi8.b}
Let $N$ be a closed orientable homologically $d$-connected $n$-manifold,
$d\ge0$.
Then the Whitney invariant
$$W:\Emb\phantom{}^{2n-d}(N)\to H_{d+1}(N,\Z_{(n-d-1)})$$
is a bijection, provided $n\ge d+3$ or $n\ge2d+4$, in the PL or DIFF cases
respectively.
\endproclaim

Theorems \wi8.a and \wi8.b were proved in [Le62, HH63, Hu69, \S11, BH70, Bo71,
Vr77] directly for {\it homotopically} $d$-connected manifolds (except the
PL case of Theorem \wi8.a) and using the deleted product method (\Wu) in
[Ha62'', We67, Sk97, Sk02] for {\it homologically} $d$-connected manifolds.
(The author is grateful to J. Boechat for indicating that [Bo71, Theorem 4.2]
needs a correction; this does not affect the main result of [Bo71].)

E.g. by Theorem \wi8.b we obtain that {\it the Whitney invariant
$W:\Emb\phantom{}^{p+2q+1}(S^p\times S^q)\to\Z_{(q)}$ is bijective
for $1\le p\le q-2$}, cf. Theorem \lc8 below.
The generator is the Hudson torus.

The PL case of Theorem \wi8.b gives nothing but Theorem \wi5.a for
$d+3\le n\le2d+1$.

Analogously to Theorem \wi8.b it may be proved that
{\it if $N$ is a closed connected non-orientable $n$-manifold, then
$$\Emb\phantom{}^{2n}(N)=\cases H_1(N,\Z_2)& n\text{ odd, }\\
\Z\oplus\Z_2^{s-1}&n\text{ even and }H_1(N,\Z_2)\cong\Z_2^s,\endcases$$
provided $n\ge3$ or $n\ge4$, in the PL or DIFF case respectively}
[Ba75, Vr77].
(There is a mistake in the calculation for the non-orientable case in
[Ha62'', We67, Theorem B].)

Because of the existence of knots the analogues of Theorem \wi8.b for $n=d+2$
in the PL case, and for (most) $n\le2d+3$ in the smooth case are false.
So for the smooth category and $n\le2d+3$ a classification is much harder:
until recently the {\it only} known concrete complete classification results were
for spheres and their disjoint unions.
Recently the following two results were obtained using the Kreck
modification of surgery theory.

\smallskip
{\bf Theorem \wi9.} [Sk06]
{\it Let $N$ be a closed parallelizable homologically $(2k-2)$-connected
$(4k-1)$-manifold.
Then the Whitney invariant $W:\Emb^{6k}_{DIFF}(N)\to H_{2k-1}(N)$ is
surjective and for each $u\in H_{2k-1}(N)$ there is a 1--1 correspondence
$\eta_u:W^{-1}u\to\Z_{d(u)}$, where $d(u)$ is the divisibility of the
projection of $u$ to $H_1(N)/Tors$.}

\smallskip
Recall that the divisibility of zero is zero and the divisibility
of $x\in G-\{0\}$ is $\max\{d\in\Z\ | \ \text{there is }x_1\in G:
\ x=dx_1\}$.
E.g. by Theorem \wi9 we obtain that {\it the Whitney invariant
$W:\Emb\phantom{}^{6k}(S^{2k-1}\times S^{2k})\to\Z$ is surjective and for each
$u\in\Z$ there is a 1--1 correspondence $W^{-1}u\to\Z_u$.}

\smallskip
{\bf Theorem \wi10.} [KS05]
{\it (a) Let $N$ be a closed connected smooth 4-manifold
such that $H_1(N)=0$ and the signature $\sigma(N)$ of $N$ is free of squares
(i.e. is not divisible by a square of an integer $s\ge2$).
Then the Whitney invariant $W:\Emb\phantom{}_{DIFF}^7(N)\to H_2(N)$ is
injective.
There exists $x_0\in H_2(N)$ such that $x_0^2=\sigma(N)$ and $x_0\mod2=w_2(N)$;
then $\im W=\{y\in H_2(N)\ |\ y^2+y\cap x_0=0\}$.

(b) Let $N$ be a closed simply-connected smooth 4-manifold embeddable into
$S^6$.
Take a composition $f:N\to S^6\subset S^7$ of an embedding and the inclusion.
Then $\#W^{-1}W(f)=12$.}

\smallskip
E.g. by Theorem \wi10.a we obtain that $\Emb_{DIFF}^7(\C P^2)$ is in 1--1
correspondence with $\{+1,-1\}\in\Z\cong H_2(\C P^2;\Z)$.
One can check that $+1$ corresponds to the standard embedding [BH70, p. 164]
and $-1$ to its composition with mirror symmetry.

\smallskip
{\bf Conjecture.} {\it Every smooth embedding $S^1\times S^1\to\R^4$ is PL
isotopic to a connected sum of a knot $S^2\to S^4$ either with the standard
embedding, or with the right Hudson torus, or with the left Hudson torus,
or with the composition of Dehn twist along the parallel and the right Hudson 
torus.}

\smallskip
A similar conjecture question can be stated for arbitrary closed 2-manifolds.
Cf. [FKV87, FKV88].

%For a composition of a homeomorphism of torus and the standard embedding
%into $\R^4$ the Whitney invariant cannot assume value (1,1)
%only values $(0,0)$, $(0,1)$
%and $(1,0)$ [Mo83, Hi]; thus the Whitney invariant is not surjective on
%the subset of such compositions contrary to a conjecture in [Sk02, \S6].
%It would be interesting to find a description of the possible values
%of the Whitney invariant for such compositions.

%It would be interesting to know if {\it the Whitney invariant is surjective for
%embeddings of the 2-torus (and even for a closed 2-manifold) into $\R^4$}.
%(Since there is an isomorphism $\pi^3_{eq}(\t N)\cong H_1(N,\Z_2)$ commuting
%with the $\alpha$-invariant defined in \Wu\ and the Whitney invariant, the
%surjectivity of $W$ for orientable surfaces is equivalent to the surjectivity
%of $\alpha^4_{PL}(N)$.)
%{\it If $N$ is a closed orientable 2-surface, $f_0:N^2\to\R^4$ is the
%standard embedding, $\gamma\subset N$ is a circle representing an element
%$[\gamma]\in H_1(N,\Z_2)$, \ $h_\gamma:N\cong N$ is the Dehn twist along
%$\gamma$, then $W(f_0\circ h_\gamma)=[\gamma]$.}[Hud69, \S11]
%\proclaim{Theorem \wi9} The Whitney invariant
%$$W:\Emb\phantom{}^m(N)\to H_{2n-m+1}(N,\Z_{(m-n-1)})$$
%is a bijection, provided $N$ is homologically $(2n-m)$-connected, $m\ge n+3$
%and $m\ge\frac{3n}2+1$ or $m\ge\frac{3n}2+2$, in the PL or DIFF cases
%respectively. \endproclaim

\subhead Low-dimensional manifolds \endsubhead

For relatively low-dimensional manifolds there are the following results
not covered by Theorems \wi2, \wi3, \wi6.a and \wi8.a.
(We need not specify whether PL or DIFF manifolds are under consideration
because every PL manifold of dimension at most 7 is smoothable.)

%Every closed 6-manifold $N$ such that $H_1(N)=0$ PL embeds into $\R^{10}$ [Sko97].

\smallskip
{\bf Theorem \wi11.}
{\it (a) A closed orientable 4-manifold $N$ PL embeds into $\R^6$ if and only
if $\overline w_2(N)=0$ [Ma78, Corollary 10.11, CS79].

(b) A closed orientable 4-manifold smoothly embeds into $\R^6$ if and only if
$\overline w_2(N)=0$ and $p_1(N)=0$ [Ma78, Corollary 10.11, CS79, Ru82].

(c) Every 2-connected closed 6-manifold is a connected sum of $S^3\times S^3$
[Sm62, Theorem B] and therefore embeds into $\R^7$.

(d) Every closed non-orientable 6-manifold $N$ such that
$\bar w_2(N)=0$ and $\bar w_3(N)=0$ PL embeds into $\R^{10}$ [Sk02].

(e) Let $N$ be a closed simply-connected 6-manifold whose homology are torsion
free, and $\bar w_2(N)=0$.

$N$ embeds into $\R^7$ if and only if $N$ is a connected sum of copies of
$S^2\times S^4$ and $S^3\times S^3$;

$N$ smoothly (or PL locally flat) embeds into $\R^8$ if and only if $p_1(N)=0$;

$N$ smoothly embeds into $\R^{10}$ [Wa66', Theorems 12 and 13].

(f) Every closed homologically 2-connected 7-manifold PL embeds into
$\R^{11}$ [Sk97, Sk02].}

%For d we need additionally the remark on the smooth case after Theorem \wi8.a
%and[Mas62, Sko97, discussion after Theorem 1.1].

\smallskip
The embeddability in $\R^{10}$ in Theorem \wi11.e is true also in the PL case,
but this is covered by Theorem \wi3.b.
%Concerning embeddability of 4-manifolds into $\R^5$ see [Co84I, Co84T, CHS02].

Take the Dold 5-manifold $N$ such that
$\bar w_{2,3}(N):=\bar w_2(N)\bar w_3(N)\ne0$ and make surgery killing
$\pi_1(N)$.
We obtain a simply connected 5-manifold $N'$ with
$\bar w_{2,3}(N')\ne0$, therefore $\bar w_3(N')\ne0$ and hence $N'$ does not
embed into $\R^8$.
This remark of Akhmetiev shows that the dimension $2n-d$ is the
minimal in Theorem \wi6.a for $n=5$ and $d=1$.

We conjecture that there exists a 1-connected 6-manifold $N$ with normal
Stiefel-Whitney class $\bar W_3(N)\ne0$ so that $N$ does not embed into $\R^9$,
cf. [Wa66', Zh75, Zh89].

\newpage
\head \lc\ Links and knotted tori \endhead

\subhead The linking coefficient \endsubhead
\smallskip
{\it Definition of the linking coefficient.}
Fix orientations of $S^p$, $S^q$, $S^m$ and $D^{m-p}$.
Assume that $m\ge q+3$ and $f:S^p\sqcup S^q\to S^m$ is an embedding.
Take an embedding $g:D^{m-q}\to S^m$ such that $gD^{m-q}$ intersects $fS^q$
transversally at exactly one point with positive sign (Figure \lc1).
Then the restriction of $g$ to $\partial D^{m-q}$ is an orientation preserving
homotopy equivalence $h:S^{m-q-1}\to S^m-fS^q$.
The induced isomorphism of homotopy groups does not depend on $g$.
{\it The linking coefficient} is
$$\lambda_{12}(f)=
[S^p\overset{f|_{S^p}}\to\to S^m-fS^q\overset{h}\to\to S^{m-q-1}]
\in\pi_p(S^{m-q-1}).$$
%\bigskip
\centerline{\epsffile{3-2.eps}}
\centerline{\it Figure \lc1} %approximately here}
\bigskip
Clearly, $\lambda_{12}(f)$ is indeed independent on $h$.

Analogously we may define $\lambda_{21}(f)\in\pi_q(S^{m-p-1})$ for $m\ge p+3$.
The definition works for $m=q+2$ if the restriction of $f$ to $S^q$
is PL unknotted (this is always so for $m\ge q+3$ by Theorem \wi5.a).
For $m=p+q+1$ there is a simpler alternative definition.

\bigskip
\centerline{\epsffile{3-5.eps}}
\centerline{\it Figure \lc2} %approximately here}
\bigskip

{\it Construction of a link with prescribed linking coefficient for
$p\le q\le m-2$.}
Define $f$ on $S^q$ to be the standard embedding into $\R^m$.
Take any CAT map $\varphi:S^p\to\partial D^{m-q}$.
Define the CAT embedding $f$ on $S^p$ by
$$S^p\overset{\varphi\times i}\to\to\partial D^{m-q}\times S^q
\subset D^{m-q}\times S^q\subset\R^m,$$
where $i:S^p\to S^q$ is the equatorial inclusion and the latter inclusion
is the standard.
See Figure \lc2.
Clearly, $\lambda_{12}(f)=\varphi$.

\smallskip
If $m\ge\frac p2+q+2$, then
$\Sigma^{\infty}:\pi_p(S^{m-q-1})\to\pi^S_{p+q+1-m}$ is an isomorphism.

Consider the \lq connected sum' commutative group structure on
$\Emb^m(S^p\sqcup S^q)$ defined for $m-3\ge p,q$ (Figure \lc3) in [Ha66, Ha66'].

\bigskip
\centerline{\epsffile{3-4.eps}}
\centerline{\it Figure \lc3} %approximately here}
\bigskip

\proclaim{The Haefliger-Zeeman Theorem \lc1} If $1\le p\le q$, then the map
$$\Sigma^{\infty}\lambda_{12}:\Emb\phantom{}^m(S^p\sqcup S^q)\to
\pi^S_{p+q+1-m}$$
is an isomorphism for $m\ge\frac p2+q+2$ and for $m\ge\frac{3q}2+2$, in the
PL and DIFF cases respectively.
\endproclaim

The surjectivity is proved above and does not require the dimension
restrictions.
The injectivity is proved in [Ha62', Ze62], or follows from the Haefliger-Weber
Theorem \wu4 and Deleted Product Lemma \wu3.a below.

By the Haefliger-Zeeman Theorem \lc1 we have the following table for
$m\ge\frac{3q}2+2$.
$$\minCDarrowwidth{1pt} \CD
m                    @=2q+2 @=2q+1 @=2q   @=2q-1 @=2q-2    @=2q-3 @=2q-4 \\
\Emb^m(S^q\sqcup S^q)@=0    @=\Z   @=\Z_2 @=\Z_2 @=\Z_{24} @=0    @=0 \endCD$$
%\subhead The $\alpha$-invariant \endsubhead
The stable suspension
of the linking coefficient can be described alternatively as
follows. For an embedding $f:S^p\sqcup S^q\to S^m$ define a map
$$\t f:S^p\times S^q\to S^{m-1}\quad\text{by}\quad
\t f(x,y)=\frac{fx-fy}{|fx-fy|}.$$
For $p\le q\le m-2$ define the $\alpha$-invariant by
$$\alpha(f)=[\t f]\in[S^p\times S^q,S^{m-1}]\overset{v^*}\to\cong\pi_{p+q}(S^{m-1})
\cong\pi^S_{p+q+1-m}.$$

\bigskip
\centerline{\epsffile{3-3.eps}}
\centerline{\it Figure \lc4} %approximately here}
\bigskip

The second isomorphism in the formula for $\alpha(f)$ is given by the
Freudenthal Suspension Theorem.
The map $v:S^p\times S^q\to\frac{S^p\times S^q}{S^p\vee S^q}\cong S^{p+q}$
is the quotient map (Figure \lc4).
The map $v^*$ is an isomorphism for $m\ge q+2$.

(For $m\ge q+3$ this follows by general position and for $m=q+2$ by the
cofibration Barratt-Puppe exact sequence of the pair
$(S^p\times S^q,S^p\vee S^q)$ and by the existence of a retraction
$\Sigma(S^p\times S^q)\to\Sigma(S^p\vee S^q)$, cf.\ [MR86, \S3].)

By [Ke59, Lemma 5.1] we have $\alpha=\pm\Sigma^{\infty}\lambda_{12}$.

Note that $\alpha$-invariant can be defined in more general situations [Ko88].

%\newpage
\subhead Borromean rings, the Whitehead link and the Trefoil knot \endsubhead
An analogue of the Haefliger-Zeeman Theorem \lc1 holds for links with many components.
However, {\it the collection of pairwise $\alpha$-invariants (or even linking
coefficients) is not injective for $2m<3n+4$ and $n$-dimensional links with
more than two components in $\R^m$.}
This is implied by the following example.

\smallskip
{\bf The Borromean Rings Example \lc2.}
{\it The Borromean rings
$$S^{2l-1}\sqcup S^{2l-1}\sqcup S^{2l-1}\to\R^{3l}$$
form a non-trivial embedding whose restrictions to 2-componented sublinks are
trivial} [Ha62, 4.1, Ha62', cf. Ma90].

\bigskip
\centerline{\epsffile{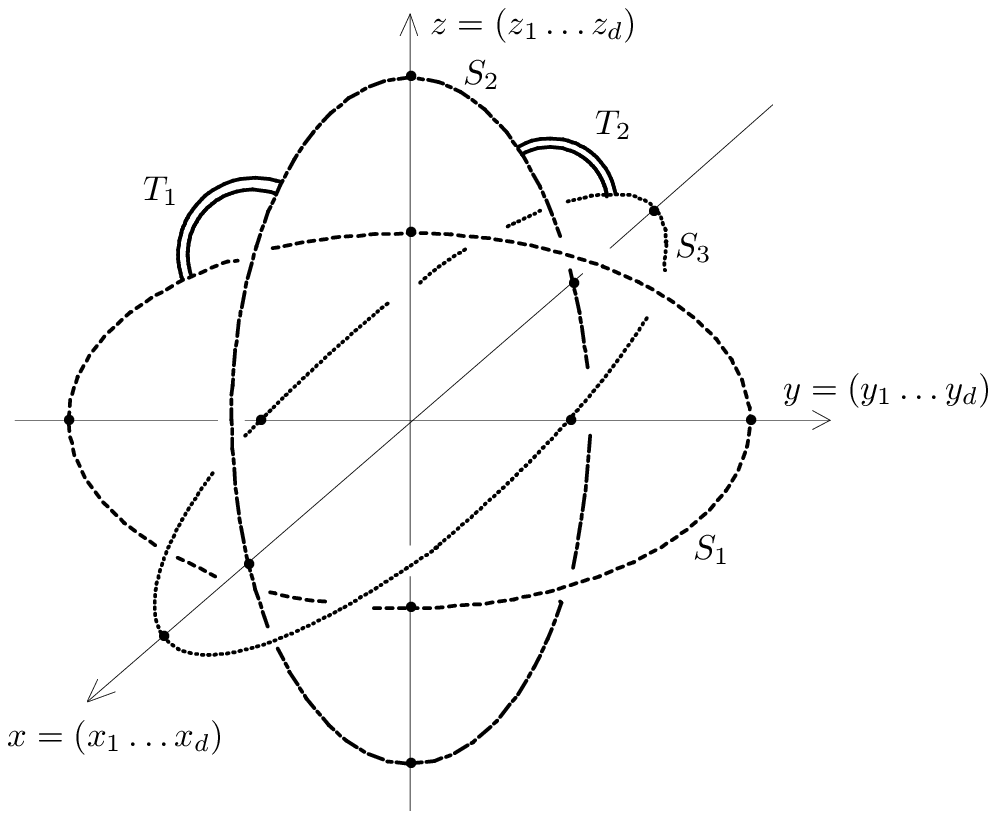}}
\centerline{\it Figure \lc5} %approximately here}
\bigskip

Consider the space $\R^{3l}$ of coordinates
$$(x,y,z)=(x_1,\dots,x_l,y_1,\dots,y_l,z_1,\dots,z_l).$$
The (higher-dimensional) Borromean rings are three embedded spheres
(Figures \lc5, \lc6.b)
given by the equations
$$\cases x=0\\ y^2+2z^2=1\endcases, \qquad \cases y=0\\ z^2+2x^2=1\endcases
\qquad\text{and}\qquad\cases z=0\\ x^2+2y^2=1 \endcases.$$

The following classical example shows that {\it the invariant
$\alpha=\pm\Sigma^\infty\lambda_{12}: Emb^m(S^p\sqcup S^q)\to\pi_{p+q+1-m}^S$ can
be incomplete for $m<\frac p2+q+2$ and links with two components}, i.e.\ that
the dimension restriction in the Haefliger-Zeeman Theorem \lc1 is sharp.

\proclaim{The Whitehead Link Example \lc3}
The Whitehead link $w:S^{2l-1}\sqcup S^{2l-1}\to\R^{3l}$ is non-trivial
although $\alpha(w)=\Sigma^\infty\lambda_{12}(w)=0$.
\endproclaim

%\bigskip
\centerline{\epsffile{3-7.eps}}
\centerline{\it Figure \lc6} %approximately here}
\bigskip

The {\it Whitehead link} is obtained from Borromean rings by joining two
components with a tube (Figure \lc6.w).
We have
$$\alpha(w)=0\quad\text{but}\quad\lambda_{12}(w)=[\iota_l,\iota_l]\ne0\quad
\text{for}\quad l\ne1,3,7.$$
Cf. [Ha62', \S3].
Note that for $l=1,3,7$ the Whitehead link is still non-trivial, although
$\lambda_{12}(w)=\lambda_{21}(w)=0$ [Ha62', \S3].

\proclaim{The Trefoil Knot Example \lc4} The trefoil knot $S^{2l-1}\to\R^{3l}$
is not smoothly trivial (but is PL trivial for $l\ge2$) [Ha62, Ha66'].
\endproclaim

The trefoil knot is obtained by joining
the three Borromean rings by two tubes (Figure \lc6.t).

If we take a cone or a suspension over any codimension 2 knot, then we obtain a
PL embedding $f$ of a ball or of a sphere which is not {\it smoothable}, i.e.
is not PL isotopic to a smooth (not necessarily standard) embedding.
This is so because $f$ is not locally flat.
Recall that an embedding $N\subset\R^m$ of a PL $n$-manifold $N$ is {\it
locally flat} if each point $x\in N$ has a closed neighborhood $U$ such that
$(U,U\cap N)\cong(D^m,D^n)$.
Observe that for $m\ge n+3$ the suspension extension $S^n\to\R^m$ of any knot
$S^{n-1}\to\R^{m-1}$ is PL isotopic to the standard embedding and is therefore
smoothable.

\proclaim{The Haefliger Torus Example \lc5} There is a PL embedding
$S^{2k}\times S^{2k}\to\R^{6k+1}$ which is (locally flat but) not PL isotopic
to a smooth embedding  [Ha62, BH70, p.165, Bo71, 6.2].
\endproclaim

In order to construct the Haefliger torus take the above trefoil knot
$S^{4k-1}\to\R^{6k}$.
Extend this knot to a conical embedding $D^{4k}\to\R^{6k+1}_-$.
By [Ha62], the trefoil knot also extends to a smooth embedding
$S^{2k}\times S^{2k}-\delet D^{4k}\to\R^{6k+1}_+$ (Figure \lc7.a).
These two extensions together form the Haefliger torus (Figure \lc7.b).

\bigskip
\centerline{\epsffile{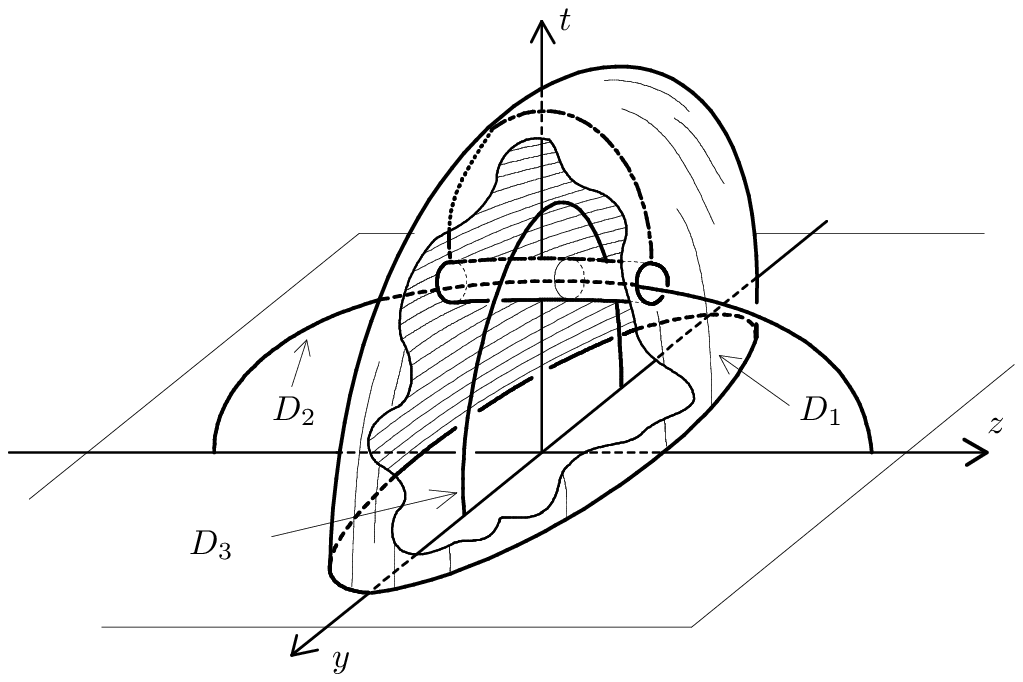}}
\centerline{\it Figure \lc7.a} %approximately here}
\bigskip
\centerline{\epsffile{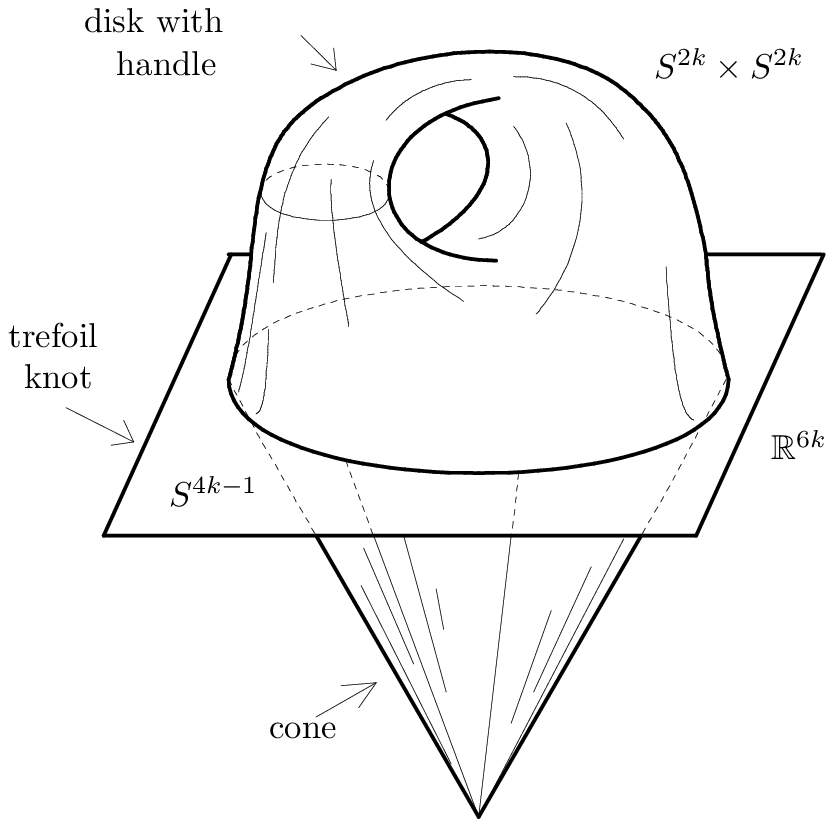}}
\centerline{\it Figure \lc7.b} %approximately here}
%\bigskip

\subhead A classification of knots and links below the metastable range
\endsubhead
Denote
$$C_q^{m-q}:=\Emb\phantom{}^m_{DIFF}(S^q).$$
The 'connected sum' commutative group structure on $C_q^{m-q}$
was defined for $m\ge q+3$ in [Ha66, cf. Ha66'].
Theorem \wi5.b states that $C_q^{m-q}=0$ for $2m\ge3q+4$.
It is known [Ha66, Mi72', KS05] that
$$C^{2k+1}_{4k-1}\cong\Z,\quad C^{2k+2}_{4k+1}\cong\Z_2,\quad
C_4^3\cong\Z_{12},\quad C^{2k}_{4k-2}=0,$$ \quad
$$C^{4s+3}_{8s+4}\cong\Z_4\quad\text{for}\quad s>0
\quad\text{and}\quad C^{4s+1}_{8s}\cong\Z_2\oplus\Z_2.$$

%\smallskip
{\bf Theorem \lc6.}
{\it (a) [Ha66', cf. Ha86] If $p,q\le m-3$, then
$$\Emb\phantom{}^m_{DIFF}(S^p\sqcup S^q)\cong
\Emb\phantom{}^m_{PL}(S^p\sqcup S^q)\oplus C^{m-q}_q\oplus C^{m-p}_p.$$
(b) [Ha66', Theorem 10.7, Sk']
If $p\le q\le m-3$ and $3m\ge2p+2q+6$, then for large enough $M$
$$\Emb\phantom{}^m_{PL}(S^p\sqcup S^q)\cong\pi_p(S^{m-q-1})\oplus
\pi_{p+q+2-m}(V_{M+m-p-1,M}).$$}

%\smallskip
The isomorphism in Theorem \lc6.b is given by the sum of
$\lambda_{12}$-invariant and {\it the $\beta$-invariant} [Sk'].

By Theorem \lc6.b (and its proof) the invariant $\lambda_{12}\oplus\lambda_{21}$
is injective for $l\ge2$ and PL embeddings $S^{2l-1}\sqcup S^{2l-1}\to\R^{3l}$;
its range is isomorphic to $\pi_{2l-1}(S^l)\oplus\Z_{(l)}$.
However, this invariant is not injective in other dimensions.

The set $\Emb^m(S^{n_1}\sqcup\dots\sqcup S^{n_s})$ for $m\ge n_i+3$ has been
described in terms of exact sequences involving homotopy groups of spheres
[Ha66, Ha66', cf. Le65, Ha86].

%\newpage
\subhead Knotted tori \endsubhead
The classification of knotted tori, i.e.\ the description of the isotopy classes of
embeddings $S^p\times S^q\to\R^m$ is an interesting problem because

(1) it has already provided many interesting examples [Al24, Ko62, Hu63, Wa65',
Ti69, BH70, Bo71, MR71, Sk02, Sk],

(2) by the Handle Decomposition Theorem it may be considered as a natural next
step (after the link theory) towards the classification of embeddings for
{\it arbitrary} manifolds, cf. [Sk05].

(3) it generalizes the classical theory of 2-componented
links of the same dimension,

(4) it reveals new interesting relations between algebraic and geometric
topology,

Denote
$$KT^m_{p,q,CAT}:=\Emb\phantom{}^m_{CAT}(S^p\times S^q).$$
Notice the change in the role of $p$ in this subsection compared with the
previous ones.
We omit CAT if a formula holds for both categories.

From the Haefliger-Zeeman Isotopy Theorem \wi6.b it follows that
$KT^m_{p,q}=0$ for $p\le q$ and $m\ge p+2q+2$.
The dimension restriction in this result is sharp by the Hudson Torus Example
\wi6.c.

\smallskip
{\bf Group Structure Theorem \lc7.} {\it The set $KT^m_{p,q}$ has
a commutative group structure for $m\ge2p+q+3$ in the smooth case
and $m\ge\max\{2p+q+2,q+3\}$ in the PL case} [Sk].

\smallskip
{\it Idea of the proof.} See Figure \lc8.
By [Sk] under the dimension assumptions for any embedding
$f:S^p\times S^q\to\R^m$ there is a {\it web}, i.e. an embedding
$$u:D^{p+1}\to\R^m\quad\text{such that}
\quad u(D^{p+1})\cap f(S^p\times S^q)=u(\partial
D^{p+1})=f(S^p\times1).$$ Moreover, a web is unique up to isotopy.

Now take two embeddings
$f_0,f_1:S^p\times S^q\to\R^m$ and their webs $D^{p+1}_0$ and $D^{p+1}_1$.
Join the centers of $D^{p+1}_0$ and $D^{p+1}_1$ by an arc.
Construct an embedding $\partial D^{p+1}\times I\to\R^m$ \lq along this arc'
so that
$$\partial D^{p+1}\times I\cap f_i(S^p\times S^q)=
\partial D^{p+1}\times i=f_i(S^p\times1)\quad\text{for}\quad i=0,1.$$
Take a \lq connected sum' of $f_0$ and $f_1$
\lq along $\partial D^{p+1}\times I$'.
The resulting embedding $S^p\times S^q\to\R^m$ is the sum of $f_0$ and $f_1$.
\qed

\bigskip
\centerline{\epsffile{3-9.eps}}
\centerline{\it Figure \lc8} %approximately here}
\bigskip

{\bf Theorem \lc8.} [HH63, Hu63, Vr77]
$$KT^{p+2q+1}_{p,q,PL}\cong\cases\Z_{(q)}&1\le p<q \\
\Z_{(q)}\oplus\Z_{(q)}&2\le p=q\endcases\qquad\text{and}
\qquad KT^{p+2q+1}_{p,q,DIFF}\cong\Z_{(q)}\quad\text{for}\quad 1\le p\le q-2.$$

Theorem \lc8 follows from Theorem \wi8.b (as well as from Theorem
\lc9 below). In the PL case of Theorem \lc8 for $p=q$ we only have
a 1--1 correspondence of sets (because Group Structure Theorem
\lc7 does not give a group structure for such dimensions). A
description of $KT^{6k}_{2k-1,2k,DIFF}$ is given after Theorem
\wi9.

This result can be generalized as follows.

\smallskip
{\bf Theorem \lc9.} {\it [Sk02, Corollary 1.5.a] If $2m\ge3q+2p+4$
or $2m\ge3q+3p+4$, in the PL or DIFF cases respectively, then}
$$KT^m_{p,q}\cong\pi_q(V_{m-q,p+1})\oplus\pi_p(V_{m-p,q+1}).$$

Note that $\pi_p(V_{m-p,q+1})=0$ for $m\ge2p+q+2$ (which is
automatic for $p\le q$ and $2m\ge3p+3q+4$). Theorem \lc9 follows
from Theorems \wu4, \wu5 and Deleted Product Lemma \wu3.b below.
For $m\ge2p+q+2$ there is an alternative direct proof [Sk], but
for $m<2p+q+2$ (when no group structure exists) no proof of
Theorem \lc9 without referring to the deleted product method is
known.

Let us construct a map $\tau:\pi_q(V_{m-q,p+1})\to KT^m_{p,q}$
giving one summand in Theorem \lc9. Recall that
$\pi_q(V_{m-q,p+1})$ is isomorphic to a group of CAT maps $S^q\to
V_{m-q,p+1}$ up to CAT homotopy. The latter maps can be considered
as CAT maps $\varphi:S^q\times S^p\to\partial D^{m-q}$. Define the
CAT embedding $\tau(\varphi)$ (Figure \lc9) as the composition
$$S^p\times S^q\overset{\varphi\times\pr_2}\to\to\partial D^{m-q}\times S^q
\subset D^{m-q}\times S^q\subset\R^m.$$

\bigskip
\centerline{\epsffile{3-12.eps}}
\centerline{\it Figure \lc9} %approximately here}
\bigskip

Let us present some calculations based on Theorem \lc9 and the
calculation of $\pi_q(V_{a,b})$ [Pa56, DP97]. Recall that
$\pi_q(V_{m-q,2})\cong\pi_q(S^{m-q-1})\oplus\pi_q(S^{m-q-2})$ for
$m-q$ even ($\cong\pi^S_{2q+1-m}\oplus\pi^S_{2q+2-m}$ for $m-q$
even and $2m\ge3q+6$) because the sphere $S^{m-q-1}$ has a
non-zero vector field. In all tables of this subsection $u^v$
means $(\Z_u)^v$.

%$$\cases\text{is adjoint to }(\pi^S_{2q+1-m}\oplus\pi^S_{2q+2-m})\otimes\Z_2
%& m-q\text{ odd, }2m\ge3q+6,\ m\ne2q+1\\
%\cong\Z_2& m=2q+1\text{ and $q\ne2$ is even}\endcases$$
%This essentially well-known fact was proved in [Sko02, \S7].

The Haefliger-Zeeman Theorem \lc1 suggests the following question:
{\it how can we describe $KT^m_{1,q}$}?
We have the following table for $2m\ge3q+6$ and for $2m\ge3q+7$,
in the PL and DIFF cases respectively.
$$\minCDarrowwidth{1pt} \CD
m                         @=2q+2 @=2q+1      @=2q @=2q-1     @=2q-2@=2q-3\\
KT^m_{1,q},\ q\text{ even}@=\Z   @=2         @=2^2@=2^2      @=24  @=0\\
KT^m_{1,q},\ q\text{ odd} @=2    @=\Z\oplus2 @=4  @=2\oplus24@=2   @=0 \endCD$$
Theorem \lc8 and [MR71] suggest the following problem: describe $KT^m_{p,q}$ for $m\le2q+p$.
We have the following table for $q\ge 4$ or $q\ge p+4$, in the PL or
DIFF cases respectively .
$$\minCDarrowwidth{1pt} \CD
p                       @=1         @=2\le p\le q-2 @=q-1         @=q  \\
KT^{p+2q}_{p,q},\ q=4s  @=2         @=0             @=2           @=0    \\
KT^{p+2q}_{p,q},\ q=4s+2@=2         @=2             @=2^2         @=2^2   \\
KT^{p+2q}_{p,q},\ q=4s+1@=\Z\oplus2 @=2^2           @=\Z\oplus2^2 @=2^4 \\
KT^{p+2q}_{p,q},\ q=4s-1@=\Z\oplus2 @=4             @=\Z\oplus4   @=4^2\endCD$$

%$$KT^{2q+p}_{p,q}=\cases 0&q=0(4)\\  \Z_2\oplus\Z_2&q=1(4)\\
%\Z_2&q=2(4)\\ \Z_4&q=3(4)\endcases$$

Classifications of {\it smooth} embeddings $S^p\times S^q\to\R^m$ for $2m\le3p+3q+3$,
as well as PL embeddings for $2m\le2p+3q+3$ is much harder
(because of the existence of smooth knots and the incompleteness of the Haefliger-Wu
invariant).
However, the {\it statements} are simple.

\smallskip
{\bf Theorem \lc10.}
{\it (a) [Sk] If $p\le q$, $m\ge2p+q+3$ and $2m\ge3q+2p+4$, then
$$KT^m_{p,q, DIFF}\ \cong\ KT^m_{p,q}\oplus G_{m,p,q}
\ \cong\ \pi_q(V_{m-q,p+1})\oplus G_{m,p,q},$$
where $G_{m,p,q}$ is certain quotient of $C^{m-p-q}_{p+q}$.

(b) [Sk06] If $1\le p\le2k-2$, then
$$KT^{6k}_{p,4k-1-p}\cong\pi_{4k-1-p}(V_{2k+p+1,p+1})\oplus\Z.$$}

\smallskip
{\bf Theorem \lc11.} [Sk] {\it (aDIFF) For each $k>1$ we have
$$KT^{6k+1}_{1,4k-1,DIFF}\cong
\pi_{2k-2}^S\oplus\pi_{2k-1}^S\oplus\Z\oplus G_k,$$
where $G_k$ is an abelian group of order 1, 2 or 4.
$$(aPL)
\qquad
KT^{6k+1}_{1,4k-1,PL}\cong\pi_{2k-2}^S\oplus\pi_{2k-1}^S\oplus\Z.$$
(b) For each $k>0$ we have
$$KT^{6k+4}_{1,4k+1,PL}\cong KT^{6k+4}_{1,4k+1,DIFF}\cong
\Z_2^a\oplus\Z_4^b\quad
\text{for some integers $a=a(k)$, $b=b(k)$ such that}$$
$$a+2b-\rk(\pi_{2k}^S\otimes\Z_2)-\rk(\pi_{2k-1}^S\otimes\Z_2)=
\cases 0         & k\in\{1,3\}\\
1                & k+1\text{ is not a power of 2} \\
1\text{ or }0    & k+1\ge8\text{ is a power of 2}\endcases.$$}

\smallskip
For a generalization of Theorem \lc11 and its relation to homotopy
groups of Stiefel manifolds see [Sk, MS04].
By Theorem \lc11 and [Pa56, To62, DP97] (see the details in [Sk]) we have
the following table.
$$\minCDarrowwidth{2pt}\CD
l                    @=2      @=3 @=4       @=5   @=6  @=7 @=8        @=9
@=10\\
KT^{3l+1}_{1,2l-1,PL}@=\Z^2\oplus2 @=4 @=\Z\oplus24\oplus2 @=2^2 @=\Z @=2
@=\Z\oplus240\oplus2 @=2\t\times2^2 @=\Z\oplus2^5\endCD$$

The following strong result was proved using a clever generalization of methods
from [Sk].

\smallskip
{\bf Theorem \lc12.}
{\it Assume that
$$p\le q,\quad p+\frac43q+2<m<p+\frac32q+2\quad\text{and}\quad m>2p+q+2.$$
The group $KT^m_{p,q,DIFF}$ is infinite if and only if either $q+1$ or
$p+q+1$ is divisible by 4} [CRS].

\smallskip
We conclude this subsection by some open problems.
It would be interesting to find $\widehat{C^{m-p-q}_{p+q}}$, at least for
particular cases.
It would be interesting to describe $KT^{6k+2}_{2k,2k+1,DIFF}$.
Note that $\# KT^{6k+2}_{2k,2k+1,DIFF}\in\{2,3,4\}$ [Sk06], cf. Theorem \lc9.b.
%\cong\Z_2\oplus\widehat{\Z_2}$ for $k>0$, where $\widehat{\Z_2}$ is either
%0 or $\Z_2$.
It would be interesting to find $KT^{6k+3}_{2k,2k+1,DIFF}$.
For this case the Whitney invariant is a surjection onto $\Z_2$, and both
preimages consist on 1 or 2 elements.
It would be interesting to find $KT^{3q+1}_{q,q,DIFF}$ for $q\ge2$.
For this case the image of the Whitney invariant is $\Z\vee\Z$
for $q$ even and is either $\Z_2\vee\Z_2$ or $\Z_2\oplus\Z_2$ for $q$ odd.
For a group $G$
define $G\vee G=\{(x,y)\in G\oplus G\ |\ \text{either $x=0$ or }y=0\}$.
The non-empty preimages of the Whitney invariant consist of 1 element for $q$
odd, of 1,2,4 elements for $q$ even $\ge4$, and of 1,2,3,4,6,12 elements
for $q=2$.

It would be interesting to find an action of the group of CAT
auto-homeomorphisms of $S^p\times S^q$ on $KT^m_{p,q,CAT}$ (E. Rees).
The above classification of knotted tori could perhaps be applied to solve for
knotted tori the Hirsch problem about the description of possible normal
bundles for embeddings of manifolds into $\R^m$, cf.\ [MR71].
The same remark holds for the following Hirsch-Rourke-Sanderson  problem
[Hi66, RS01, cf. Ti69, Vr89]: {\it which embeddings
$N\to\R^{m+1}$ are isotopic to embeddings $N\to\R^m$?}

%\newpage
\head \vk\ The van Kampen obstruction \endhead
\subhead The embeddability of $n$-complexes in $\R^{2n}$ \endsubhead
By the General Position Theorem \wi1.a, the first non-trivial case of the
Embedding Problem is the investigation of the embeddability of $n$-polyhedra in
the Euclidean space $\R^{2n}$, cf. Example \wi1.b.
For $n=1$ this problem was solved by the Kuratowski criterion [Ku30], see also
[RS96, \S2, Sk05'] and references there.
However for $n>1$, such a simple criterion does not exist [Sa91].
(Note that there are infinitely many closed non-orientable 2-surfaces, which
do not embed into $\R^3$ and these do not contain a common subspace non-embeddable
into $\R^3$.)
In [Ka32] an obstruction to the embeddability of $n$-polyhedra in $\R^{2n}$
was constructed for arbitrary $n$ (see also a historical remark at the end of
\Wu).

\bigskip
\centerline{\epsffile{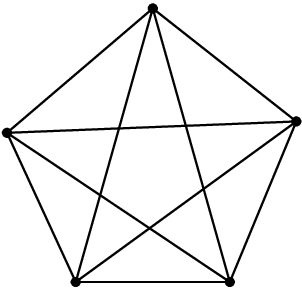}}
\centerline{\it Figure \vk1} %approximately here}
\bigskip

To explain the idea of van Kampen, we sketch a proof of the nonplanarity of
$K_5$ (i.e.\ of the complete graph with 5 vertices, Figure \vk1).
Take any general position map $f:K_5\to\R^2$.
For each two edges $\sigma,\tau$ the intersection $f\sigma\cap f\tau$
consists of a finite number of points.
Let $v_f$ be

{\it the sum mod~2 of the numbers $|f\sigma\cap f\tau|$ over all
non-ordered pairs $\{\sigma,\tau\}$ of disjoint edges of $K_5$.}

For the map $f$, shown in Figure~\vk1, $v_f=1$.
Every general position map $f:K_5\to\R^2$ can be transformed to any other such
map through isotopies of $\R^2$ and \lq the Reidemeister moves' for graphs in
the plane from Figure~\vk2.

\bigskip
\centerline{\epsffile{5-2.eps}}
\centerline{\it Figure \vk2} %approximately here}
\bigskip

This assertion is proved analogously to the Reidemeister theorem for knots.
We will not prove it, since it is needed only for this sketch proof and not
for the rigorous proof.
For each edge of $K_5$ with vertices $a,b$, the graph $K_5-\{a,b\}$, obtained
by deleting from $K_5$ the vertices $a,b$ and the interiors of the edges adjacent to
$a,b$, is a circle (this is the very property of $K_5$ we need for the proof).
Therefore $v_f$ is invariant under the \lq Reidemeister moves'.
Hence $v_f=1$ for {\it each} general position map $f:K_5\to\R^2$.
So $K_5$ is nonplanar.
(For a proof without use of the assertion on the Reidemeister moves see below or
[Sk05''].)

Similarly, one can prove that the graph $K_{3,3}$ (three houses and three
wells) is not embeddable into $\R^2$ and that the 2-skeleton of the 6-simplex
is not embeddable into $\R^4$, cf.\ Example \wi1.b.

\subhead Ramsay link theory \endsubhead
Now let us discuss some generalizations of the above proof which are
interesting in themselves and are used in \Ex.
From that proof one actually gets a stronger assertion.
Let $e$ be an edge of $K_5$ and $\Sigma^1$ the cycle in
$K_5$, formed by the edges of $K_5$ disjoint with $e$.  Then
$K_5-\delet e$ embeds into $\R^2$ (Figure \vk3) and for each
embedding $g:K_5-\delet e\to\R^2$ the $g$-images of the ends of
$e$ (which form a 0-sphere $g\partial e=g\Sigma^0$) lie on
different sides of $g\Sigma^1$.

\bigskip
\centerline{\epsffile{5-3.eps}}
\centerline{\it Figure \vk3} %approximately here}
\bigskip

Moreover, let $e$ be a 2-simplex of the 2-skeleton $\Delta^2_6$
of the 6-simplex  and $P=\Delta^2_6-\delet e$.
Then $P$ embeds into $\R^4$.
Let $\Sigma^1=\partial e$ and $\Sigma^2$ be the sphere formed by 2-faces disjoint
with $e$.
Then for each embedding $P\to\R^4$ the images of these spheres link with a non-zero
(more precisely, with an odd) linking number [Fl34].

Using the above idea one can prove the following.

{\it For any embedding $K_6\to\R^3$ there are two cycles in $K_6$ whose images
are linked with an odd linking coefficient} [Sa81, CG83, for generalizations
see RST95, Ne98, BKK02, Sk03, PS05].

\subhead The van Kampen obstruction mod 2 \endsubhead
Fix a triangulation $T$ of a polyhedron $N$.
The space
$$\t T=\cup\{\sigma\times\tau\in T\times T\ |\ \sigma\cap\tau=\emptyset\}$$
is called the {\it simplicial deleted product} of $N$.
By the Simplicial Deleted Product Lemma \wu3.a below the equivariant homotopy
type of $\t T$ depends only on $N$, so we write $\t N$ instead of $\t T$ in
this section.
Consider the 'exchanging factors' involution $t:\t N\to\t N$ defined by
$t(x,y)=(y,x)$.
Let $N^*=\t N/t$.

\smallskip
{\it Definition of the mod 2 van Kampen obstruction $v(N)$ for $n=\dim N=1$.}
%Throughout this subsection we shall omit $\Z_2$ coefficients
%from the notation of the (co)chain and (co)ho\-mo\-logy groups.
For any general position PL map $f:N\to\R^2$ and disjoint edges $\sigma,\tau$
of $T$, the intersection $f(\sigma)\cap f(\tau)$ consists of a finite number of
points.
Let
$$v_f(\sigma,\tau)=|f(\sigma)\cap f(\tau)|\mod2.$$
Then $v_f$ is an element of the group $C^2(N^*;\Z_2)$ of vectors (which are
called {\it cochains}) with components from $\Z_2$ indexed by 2-simplices of
$N^*$, i.e.\ by non-ordered products of disjoint edges of $N$.

This cochain $v_f$ is invariant under isotopy of $\R^2$ and the first four
\lq Reidemeister moves' (Figure \vk2.I--IV).
The fifth \lq Reidemeister move' (passing of $f\tau$ through $fa$, see
Figure \vk2.V) adds to $v_f$ the vector that assumes value 1 on the class of
2-simplex $\sigma\times\tau$ for $a\in\sigma$, and value 0 on the other
2-simplices of $N^*$.
This vector is denoted by $\delta[a\times\tau]$.
Denote by $B^2(N^*;\Z_2)$ the subgroup of $C^2(N^*;\Z_2)$ generated by
elementary coboundaries.
%of the elementary cochain from $C^1(N^*)$, that assumes value 1 on the class
%of 1-simplex $a\times\tau$, and value 0 on the other 2-simplices
%(the reader can consider this as a {\it definition} of the
%coboundary for this particular case).
{\it The van Kampen obstruction mod2} is the equivalence class
$$v(N):=[v_f]\in H^2(N^*;\Z_2)=C^2(N^*;\Z_2)/B^2(N^*;\Z_2).$$
(Since $\dim N^*=2$, it follows that $C^2=Z^2$.)

\smallskip
{\it Sketch of a proof that $v(N)$ is indeed independent of $f$
(without using the assertion about \lq the Reidemeister moves'
which is not proved here).} We follow [Sh57, Lemma 3.5]. Consider
an arbitrary general position homotopy $F:N\times I\to\R^2\times
I$ between general position maps $f_0,f_1:N\to\R^2$. Colour a
non-ordered pair $a\times\tau$ in red if $F(a\times I)$ intersects
$F(\tau\times I)$ in an odd number of points. Then
$v_{f_0}-v_{f_1}$ is the sum of $\sum\delta(a\times\tau)$ over all
red pairs $a\times\tau$. Hence $v_f$ is independent on $f$. \qed

\smallskip
Clearly, $v(N)=0$ for all planar graphs $N$.

Analogously one defines the mod 2 van Kampen obstruction
$v(N)\in H^{2n}(N^*;\Z_2)$ to the embeddability of an $n$-polyhedron $N$
into $\R^{2n}$.

%$$V_f(\sigma,\tau)=f\sigma\cdot f\tau=
%\sum\limits_{P\in f\sigma\cap f\tau}\sign P,$$
%where $\sign P=+1$, if the positive $n$-bases of $f\sigma$ and $f\tau$
%(in this order???) constitute a
%positive $2n$-basis $\R^{2n}$, and $\sign P=-1$ otherwise.
%Clearly, $V_f(\sigma\times\tau)=(-1)^nV_f(\tau\times\sigma)$.
%Note that $H_s^{2n}(\t N)\cong H^{2n}(N^*)$ for $n$ even [FKT94].

\subhead The integral van Kampen obstruction \endsubhead
Fix a triangulation of $N$ and define $\t N$, $t$ and $N^*$ as above.
Choose an orientation of $\R^{2n}$ and orientations of the $n$-simplices of
$\t N$.
The latter give orientations on $2n$-simplices of $\t N$.
Clearly, $t(\sigma\times\tau)=(-1)^n(\tau\times\sigma)$ (the case $n=1$ helps
to check the sign) [Sh57, p.257, above].
%orientations for the $2n$-simplices of $\t N$ so that the symmetry with
%respect to $\diag N$ preserves these orientations.

For any general position map $f:N\to\R^{2n}$ and any two disjoint $n$-simplices
$\sigma,\tau$ of $N$ the intersection $f(\sigma)\cap f(\tau)$ consists of a
finite number of points.
Define the {\it the intersection cochain}
$$V_f\in C^{2n}(\t N)\quad\text{by the formula}\quad V_f(\sigma,\tau)=
f\sigma\cdot f\tau:=\sum\limits_{P\in f\sigma\cap f\tau}\sign P.$$
Here $\sign P=+1$, if for the positive $n$-bases $s_1,\dots,s_n$ and
$t_1,\dots,t_n$ of $\sigma$ and $\tau$, respectively, we have that
%such that $s_1,\dots,s_n,t_1,\dots,t_n$ is a positive $2n$-base of
%$\sigma\times\tau$ and
$fs_1,\dots,fs_n,ft_1,\dots,ft_n$ is a positive $2n$-base of $\R^{2n}$;
and $\sign P=-1$ otherwise.
Clearly [Sh57],
$$V_f(\sigma\times\tau)=(-1)^nV_f(\tau\times\sigma)=V_f(t(\sigma\times\tau)).$$
So $V_f$ induces a cochain in the group $C^{2n}(N^*)$.
Denote this new cochain by the same notation $V_f$ (we use the old $V_f$
only in the proof of Lemma \vk2, so no confusion will arise).
{\it The van Kampen obstruction} is the equivalence class
$$V(N):=[V_f]\in H^{2n}(N^*)=C^{2n}(N^*)/B^{2n}(N^*).$$
This class is independent on $f$ (this is proved analogously to the proof for
$v(N)$).
Clearly, $V(N)$ is an obstruction to the embeddabiltiy of $N$ into $\R^{2n}$.

One can easily show that $V(N)$ depends on the choice of orientations of
$\R^{2n}$ and of the $n$-simplices of $N$ only up to an automorphism of the
group $H^{2n}(N^*)$.

(The author is grateful to S. Melikhov for indicating that in [FKT94, Kr00,
BKK02] the signs are not accurate and so the Van Kampen obstruction for $n$ odd
erroneously assumes its values in in $H_{\Z_2}^{2n}(\t N;\Z)$, where the
involution acts on $\t N$ by exchanging factors and on $\Z$ by $(-1)^n$.)

\subhead The van Kampen-Shapiro-Wu Theorem \endsubhead

\smallskip
{\bf Theorem \vk1.} {\it If an $n$-polyhedron $N$ embeds into $\R^{2n}$, then
$V(N)=0$. For $n\ne2$ the converse is true, whereas for $n=2$ it is not}
[Ka32, Sh57, Wu58, FKT94].

\smallskip
%Theorem \vk1 implies the PL case of the statement (a) after Theorem \wi4
%(and Theorem \wi2.a) for $n\ge3$.
The necessity of $V(N)=0$ in Theorem \vk1 was actually proved
in the construction of the van Kampen obstruction.
The sufficiency in Theorem \vk1 for $n\ge3$ follows from Lemmas \vk2 and \vk3
below, and for $n=1$ is obtained using the Kuratowski graph planarity
criterion.
A counterexample to the completeness of the van Kampen obstruction for $n=2$
is presented in \Ex.

A map $g:N\to\R^m$ is of a polyhedron $N$ is called a {\it non-degenerate
almost embedding} if there exists a triangulation $T$ of $N$ such that
$g|_\alpha$ is an embedding for each $\alpha\in T$ and
$g\alpha\cap g\beta=\emptyset$ for each $\alpha\times\beta\subset\t T$.

\smallskip
{\bf Lemma \vk2.}
{\it If $N$ is an $n$-polyhedron and $V(N)=0$, then there exists a general
position non-degenerate almost embedding $g:N\to\R^{2n}$} [cf. FKT94, Lemmas
2 and 4].

\smallskip
{\it Sketch of the proof.}
Let $T$ be a triangulation of $N$.
Let $\varphi:N\to\R^{2n}$ be a map linear on the simplices of $T$.
Then $\varphi$ is non-degenerate.
The condition $V(N)=0$ implies that $V_\varphi\in C^{2n}(\t N)$ is an symmetric
coboundary.
Hence $V_\varphi$ equals to the sum of some \lq elementary' symmetric coboundaries
$\delta(\sigma^n\times\nu^{n-1})+\delta t(\sigma^n\times\nu^{n-1})$ over some
$\sigma^n\times\nu^{n-1}\in\t T$.
Applying the van Kampen finger moves (higher-dimensional analogue of Figure \vk2.V)
for all pairs $\sigma^n\times\nu^{n-1}$
from this sum we obtain a general position nondegenerate map $f:N\to\R^{2n}$ such
that $f\alpha\cdot f\beta=0$ for each $\alpha,\beta\in\t T$.

Then by induction on pairs of $n$-simplices of $\t T$ and using the Whitney
trick (see below) in the inductive step we obtain the required map $g$.
See the details in [FKT94].
\qed

\smallskip
Let us illustrate the application of the Whitney trick by the following
argument.

\smallskip
{\it Sketch of the proof of Theorem \wi2.a in the smooth case, stating that
 every smooth $n$-manifold $N$ smoothly embeds into $\R^{2n}$}.
For $n\le2$ the proof is trivial, so assume that $n\ge3$.
Using the higher-dimensional analogue of the first Reidemeister move
(Figure \vk2.I), any smooth general position map $f:N\to\R^{2n}$ can be
modified so that a single self-intersection point with a prescribed sign will
be added.
Hence there exists a general position map $f:N\to\R^{2n}$ whose
self-intersections consist of an even number of isolated points, with algebraic
sum zero.

In order to conclude the proof, we \lq kill' these double points in pairs.
This procedure is analogous to the second Reidemeister move (Figure \vk2.II)
and is called the Whitney trick.
More precisely, take two double points of opposite sign:
$$x_1,y_1,x_2,y_2\in N\quad\text{so that}\quad f(x_1)=f(x_2),
\quad f(y_1)=f(y_2).$$
Join $x_1$ to~$y_1$ and $x_2$ to~$y_2$ by arcs $l_1$ and~$l_2$ so that these
double points have \lq opposite signs' along these arcs (Figure \vk4).

\bigskip
\centerline{\epsffile{8-8.eps}}
\centerline{\it Figure \vk4} %approximately here}
\bigskip

By general position ($n\ge2$), we may assume that the restrictions $f|_{l_1}$
and $f|_{l_2}$ are embeddings and that $l_1$ and $l_2$ do not contain other
double points of $f$.
Since $n\ge3$, by general position we can embed a 2-disk~$C$ into $\R^{2n}$
so that
$$\partial C=f(l_1)\cup f(l_2)\quad\text{and}\quad C\cap f(N)=\partial C.$$
Such a disk~$D$ is called Whitney's disk.
We can move the $f$-image of a regular neighborhood of~$l_1$ in $N$ \lq along' $C$
so that we \lq cancel' the double points $f(x_1)=f(x_2)$ and $f(y_1)=f(y_2)$.
For details see [Ad93, Pr].
\qed

\smallskip
{\bf The Freedman-Krushkal-Teichner Lemma \vk3.} {\it If there is
a non-degenerate almost embedding $g:N\to\R^{2n}$ of an
$n$-polyhedron $N$ and $n\ge3$, then there is an embedding
$f:N\to\R^{2n}$} [FKT94, Lemma 5, cf. We67, \S6, Sk98].

\demo{Proof} Take a triangulation $T$ satisfying to the properties from the
definition of an non-degenerate almost embedding.
We may assume by induction that
$g\alpha\cap g\beta=g(\alpha\cap\beta)$ for each
$\alpha,\beta\in T$ such that $\alpha\cap\beta\ne\emptyset$ except
for $(\alpha,\beta)=(\sigma^n,\tau^n)$.
In addition, we may assume that $g\delet\sigma^n\cap g\delet\tau^n$ is
a point (say, $p$).
Let $v$ be a point of $\sigma^n\cap\tau^n$.
Take PL arcs $l_1\subset v\cup\delet\sigma^p$ and
$l_2\subset v\cup\delet\tau^q$ joining $v$ to $g^{-1}p$ and containing no
self-intersection points of $g$ in their interiors (Figure \vk5).

\bigskip
\centerline{\epsffile{8-9.eps}}
\centerline{\it Figure \vk5} %approximately here}
\bigskip

Then $g(l_1\cup l_2)$ is a circle.
Since $n\ge3$, this circle bounds a PL embedded 2-disk
$D\subset\R^{2n}$.
We have by general position ($n+2<2n$) that
$\delet D\cap gN=\emptyset$.
There is a small neighborhood $D^{2n}$ of $D$ in $\R^{2n}$ $\rel g(v)$ that
is PL homeomorphic to the $2n$-ball and such that $g^{-1}D^{2n}$ is a
neighborhood of $l_1\cup l_2$ in $N$ $\rel v$ and is homeomorphic to the wedge
$D^n\vee D^n$.
By the Unknotting Wedges theorem [Li65], the restriction
$g:\partial g^{-1}D^{2n}\to\partial D^{2n}$ is unknotted.
Hence it can be extended to an embedding $h:\partial g^{-1}D^{2n}\to D^{2n}$.
In order to conclude the proof, set $f$ equal
to $g$ on $N-\partial g^{-1}D^{2n}$ and to $h$ on $\partial g^{-1}D^{2n}$.
\qed\enddemo

Observe that for the embedding $f$ constructed above we have
$\t f\simeq_{eq}\t g$ on $\t T$.

\subhead Generalizations of the van Kampen obstruction \endsubhead
%The constructions above can be generalized in several ways.
The idea of the van Kampen obstruction can be applied to calculate the minimal
$m$ such that a polyhedron, which is a product of graphs, embeds into $\R^m$
[Sk03, cf.\ Ga92, ARS01'].

Analogously, one can construct {\it the van Kampen-Wu invariant}
$U(f)\in H^{2n}(N^*)$ of an embedding $f:N\to\R^{2n+1}$.

\smallskip
{\bf Theorem \vk4.}
{\it If embeddings $f,g:N\to\R^{2n+1}$ of a finite $n$-polyhedron $N$ are
isotopic, then $U(f)=U(g)$.
For $n\ge2$ the converse is true, whereas for $n=1$ it is not} [Wu65].

\smallskip
Note that

{\it embeddings $f,g:N\to\R^3$ of a graph $N$ such that
$U(f)=U(g)$ are {\it homologous} [Ta95]}.

In this paper we shall not present a proof of Theorem \vk4.
For $n\ge2$ it is proved analogously to Theorem \vk1 using the
ideas of [RS99, \S12], and for $n=1$ it is trivial.

As it was pointed out by Shapiro, when $V(N)=0$ (and hence $N$ is embeddable
in $\R^{2n}$ for $n\ge3$), one can construct the `second obstruction' to
the embeddability of $N$ in $\R^{2n-1}$, etc.

For a subpolyhedron $A$ of a polyhedron $N$ one can analogously define the
obstruction to extending given embedding $A\subset\partial B^m$ to an
embedding $N\to B^m$ [FKT94].
The relative van Kampen obstruction is complete for $n\ne2$ (for $n\ge3$ see
[Wu65] and for $n=1$ this follows from a relative version of the Kuratowski
criterion [Sk05']) and is incomplete for $n=2$ (\Ex).

For the van Kampen obstruction for {\it approximability by embeddings} see
[CRS98, \S4, RS98, ARS02, \S4, RS02, Me02, Sk03'].

%\newpage
\head \wu\ The Haefliger-Wu invariant \endhead
\subhead Basic idea \endsubhead
\lq The complement of the diagonal' and \lq the Gauss map' ideas play a great
role in different branches of mathematics [Gl68, Va92].
The Haefliger-Wu invariant is a manifestation of these ideas in the theory of
embeddings.
The complement to the diagonal idea originated from two celebrated
theorems: the Lefschetz Fixed Point Theorem and the Borsuk-Ulam Antipodes
Theorem.
In order to state the latter denote the antipode of a point $x\in S^n$ by $-x$
and recall that a map $f:S^n\to S^m$ between spheres is {\it equivariant}
(or {\it odd}), if $f(-x)=-f(x)$ for each $x\in S^n$.

\proclaim{The Borsuk-Ulam Theorem \wu1}
(a) For any map $f:S^n\to\R^n$ there exists $x\in S^n$ such that
$f(x)=f(-x)$.

(b) There are no equivariant maps $S^n\to S^{n-1}$.

(c) Every equivariant map $S^n\to S^n$ is not homotopic to the constant map.
\endproclaim

\demo{Sketch of a deduction of (a) and (b) from (c)}
Part (c) is non-trivial, see the proof e.g. in [Pr04, 8.8].
%can be proved by induction on $n$.
Part (b) follows from (c) because if $\varphi:S^n\to S^{n-1}$ is an
equivariant map, then the restriction $\varphi|_{S^{n-1}}$ extends to $S^n$
and therefore is null-homotopic.
In order to present the idea of the Gauss map in the simplest case, let us
deduce (a) from (b).
Suppose to the contrary, that there exists a map $f:S^n\to\R^n$ which does not
identify any antipodes.
Then a map
$$\t f:S^n\to S^{n-1}\quad\text{is well-defined by}\quad
\t f(x)=\frac{f(x)-f(-x)}{|f(x)-f(-x)|}.$$
Evidently, $\t f$ is equivariant.
This contradicts to (b).
\qed\enddemo

%Let us show how to use the Borsuk-Ulam theorem in order to construct an
%$n$-polyhedron $N$ non-embeddable into $\R^{2n}$ [cf.\ Sch84, KaRe].

\demo{Construction of Example \wi1.b}
We present a simplified construction invented by Schepin
and the author (and, possibly, others), cf.\ [Sc84, Appendix, RS01'].
Let $T$ be a triod, i.e.\ the graph with four vertices $O$, $A$, $B$, $C$
and three edges $OA$, $OB$ and $OC$.
The product $T^{n+1}$ is a cone over some $n$-polyhedron $N$.

In order to prove that $N$ does not embed into $\R^{2n}$ it suffices to prove
that $T^{n+1}$ does not embed in $\R^{2n+1}$.
Suppose to the contrary that there is an embedding $f:T^{n+1}\to\R^{2n+1}$.
Let $p:D^2\to T$ be a map which does not identify
any antipodes of $S^1=\partial D^2$ (e.g.  the map from Figure  \wu1).
It is easy to check that the map
$p^{n+1}|_{\partial D^{2n+2}}:\partial D^{2n+2}\to T^{n+1}$ also does not
identify any antipodes.
Then the composition of $p^{n+1}$ and $f$ again does not identify antipodes.
This contradicts the Borsuk-Ulam Theorem \wu1.a.
\qed\enddemo

\bigskip
\centerline{\epsffile{4-1.eps}}
\centerline{\it Figure \wu1} %approximately here}
\bigskip

\subhead Definition of the Haefliger-Wu invariant \endsubhead
%To formulate the deleted product condition, we need the following definition.
The {\it deleted product} $\t N$ of a topological space $N$ is the product of
$N$ with itself, minus the diagonal:
$$\t N=\{(x,y)\in N\times N\ |\ x\ne y\}.$$

\bigskip
\centerline{\epsffile{4-2.eps}}
\centerline{\it Figure \wu2} %approximately here}
\bigskip

This is the configuration space of ordered pairs of distinct points of $N$.

Now suppose that $f:N\to\R^m$ is an embedding.
Then the map $\t f:\t N\to S^{m-1}$ is well-defined by the Gauss formula
$$\t f(x,y)=\frac{f(x)-f(y)}{|f(x)-f(y)|}.$$

\bigskip
\centerline{\epsffile{4-3.eps}}
\centerline{\it Figure \wu3} %approximately here}
\bigskip

This map is equivariant with respect to the 'exchanging factors' involution
$t(x,y)=(y,x)$ on $N$ and the antipodal involution $a_{m-1}$ on $S^{m-1}$.
Thus the existence of at least one equivariant map $\t N\to S^{m-1}$ is a
necessary condition for the embeddability of $N$ in $\R^m$.

%The Haefliger-Wu (deleted product) invariant is constructed as follows.
For isotopic embeddings $f_0,f_1:N\to\R^m$ and an isotopy
$f_t:N\times I\to\R^m$ between them the homotopy $\t f_t$ is an equivariant
homotopy between $\t{f_0}$ and $\t{f_1}$.
Hence the existence of  an equivariant homotopy between $\t{f_0}$ and $\t{f_1}$
is necessary for the embeddings $f_0,f_1:N\to\R^m$ to be isotopic.
The {\it Haefliger-Wu invariant} $\alpha(f)$ of the embedding $f$ is the
equivariant homotopy class of the map $\t f$, cf. [Wu59, Ha61, Gr86, 2.1.E].

Let $\pi^{m-1}_{eq}(\t N)=[\t N;S^{m-1}]_{eq}$ be the set of equivariant maps
$\t N\to S^{m-1}$ up to equivariant homotopy.
Thus the Haefliger-Wu invariant is a mapping
$$\alpha=\alpha_{CAT}^m(N):\Emb\phantom{}_{CAT}^m(N)\to \pi^{m-1}_{eq}(\t N)
\quad\text{defined by}\quad\alpha(f)=[\t f]\in\pi^{m-1}_{eq}(\t N).$$

\subhead Calculations of the Haefliger-Wu invariant \endsubhead
It is important that using algebraic topology methods the set
$\pi^{m-1}_{eq}(\t N)$ can be explicitly calculated in many cases
Let us give several examples, cf. [FS59, CF60, beginning of \S2, Ha62'', Ha63,
1.7.1, Ba75, Ad93, 7.1, RS99, Sk02].
We denote by \lq $\cong$' the existence of a 1--1 correspondence between sets.

\proclaim{The Deleted Product Lemma \wu3}
(a) $\pi^{m-1}(\t{S^p\sqcup S^q})\cong\pi_{p+q+1-m}^S$ for $m\ge q+2$.

(b) $\pi^{m-1}(\t{S^p\times S^q})\cong
\pi_q(V_{m-q,p+1})\oplus\pi_p(V_{m-p,q+1})$ for
$2m\ge3q+2p+4$.
\endproclaim

The Deleted Product Lemma \wu3.a and \wu3.b is proved in \Kt\ and follows from
the Torus Lemma \kt1, respectively.
From the proof of the Deleted Product Lemma \wu3.a it follows that the
$\alpha$-invariants of \Lc\ and of \Wu\ for $N=S^p\sqcup S^q$ indeed coincide.

\proclaim{The Simplicial Deleted Product Lemma \wu3.c}
Fix a triangulation $T$ of a polyhedron $N$.
The simplicial deleted product
$$\t T=\cup\{\sigma\times\tau\in T\times T\ |\ \sigma\cap\tau=\emptyset\}$$
is an equivariant strong deformation retract of $\t N$ [Sh57, Lemma 2.1,
Hu60, \S4].
\endproclaim

\demo{Sketch of a proof} %Consider cell decomposition of $\t T$ by cells
Denote
$$E_{\sigma\tau}:=\cup\{U_\sigma\times U_\tau\ |\  U_\sigma,\ U_\tau
\text{ non-empty faces of }\sigma,\ \tau,
\text{ respectively, and }U_\sigma\cap U_\tau=\emptyset\}.$$
Then $\sigma\times\tau\cong E_{\sigma\tau}*\diag(\sigma\cap\tau)$.
So for
$\sigma\cap\tau\ne\emptyset$ there is an equivariant strong deformation
retraction $\sigma\times\tau-\diag(\sigma\cap\tau)\to E_{\sigma\tau}$.
These retractions agree on intersections, so together they form an
equivariant strong deformation retraction $\t N\to\t T$.
\qed\enddemo

Now we sketch how to deduce in a purely algebraic way all the necessary
conditions for embeddability and isotopy presented in \Wi\ and \Vk\ from
the \lq deleted product necessary conditions'.
Let $N$ be a polyhedron (in particular, a smooth manifold).
By the Simplicial Deleted Product Lemma \wu3.c, there exists an equivariant map
$\t N\to S^{m-1}$ if and only if there exists an equivariant map
$\t T\to S^{m-1}$.
Define an $S^{m-1}$-bundle
$$\gamma:\frac{\t T\times S^{m-1}}{(x,y,s)\sim(y,x,-s)}
\overset{S^{m-1}}\to\to\frac{\t T}{(x,y)\sim(y,x)},\quad
\text{by}\quad\gamma[(x,y,s)]:=[(x,y)].$$
The existence of an equivariant map $\t T\to S^{m-1}$ is equivalent to the
existence of a cross-section of $\gamma$; the existence of an equivariant
homotopy between $\t{f_0}$ and $\t{f_1}$ is equivalent to the equivalence of
the corresponding cross-sections of the bundle $\gamma$.
Thus the existence of either can be checked using the methods of obstruction
theory.
In particular, the Whitney and the van Kampen obstructions (\Wi, \Vk) are
the first obstructions to the existence of a cross-section of the bundle
$\gamma$ [Wu65]; the Whitney and the van Kampen-Wu invariants (\Wi, \Vk) are
the first obstructions to an equivariant homotopy of $\t{f_0}$ and $\t{f_1}$
[Wu65].

\subhead The completeness of the Haefliger-Wu invariant \endsubhead
However trivial the \lq deleted product necessary conditions' may seem, the
above shows that they are very useful.
Thus it is very interesting to find out for which cases they are also
{\it sufficient} for embeddability and isotopy, i.e.\ for which cases the
following assertions hold (the converses of which were just proved).

(e) If there exists an equivariant map $\Phi:\t N\to S^{m-1}$, then $N$
embeds into $\R^m$.

(s) If there exists an equivariant map $\Phi:\t N\to S^{m-1}$, then there
exists an embedding $f:N\to\R^m$ such that $\t f\simeq_{eq}\Phi$.

(i) If $f_0,f_1:N\to\R^m$ are two embeddings and $\t{f_0}\simeq_{eq}\t{f_1}$,
then $f_0$ and $f_1$ are isotopic.

Clearly, (s) and (i) are the surjectivity and the injectivity of $\alpha$.
Obviously, (s) implies (e).
From (e) it follows that
{\it if $N$ TOP embeds into $\R^m$, then $N$  PL or DIFF embeds into $\R^m$}
(in particular, that PL or DIFF embeddability of $N$ into $\R^m$ does not
depend on PL or DIFF structures on $N$).
Condition (i) has analogous corollary.

%if $N$ quasi-embeds into $\R^m$, then $N$ embeds into $\R^m$.

Thus the surjectivity and the injectivity of $\alpha$ are directly related
to the Embedding and Knotting Problems.

\proclaim{The Haefliger-Weber Theorem \wu4}
[Ha63, Theorem 1', We67, Theorems 1, 1']
For embeddings $N\to\R^m$ of either an $n$-polyhedron or a smooth $n$-manifold
$N$ the Haefliger-Wu invariant is
$$\text{bijective if}\quad 2m\ge3n+4\quad\text{and surjective if}
\quad 2m=3n+3.$$
\endproclaim

The {\it metastable dimension restrictions} in the Haefliger-Weber Theorem
\wu4 are sharp in the smooth case by the Trefoil Knot Example \lc4 (and
other examples of smooth knots [Ha66]), and by Theorems \wi7.c, \wi7.d
(because the PL embedability $N\to\R^{2l+2}$ implies the existence of an
equivariant map $\t N\to S^{2l+1}\subset S^{3l-2}$).
Such dimension restrictions appeared also in the PL cases of the classical
theorems on embeddings of highly-connected manifolds and of Poincar\'e
complexes (see Theorems \wi5, \wi6 and [Ru73, BM99, BM00]), but were later
weakened to $m\ge n+3$.
So since 1960's it was conjectured by Viro, Dranishnikov, Koschorke, Szucs,
Schepin and others that also in Theorem \wu4 for the PL case and connected $N$
the metastable dimension restrictions can be weakened to $m\ge n+3$ (possibly
at the price of adding the $p$-fold Haefliger-Wu invariant, see the subsection
'the generalized Haefliger-Wu invariant').
This turned out to be false not only for polyhedra (Examples \wu7.d and \wu9.c
below) but even for PL manifolds (Examples \wu7.b, \wu7.c and \wu8.b below).
So it is surprising that the metastable dimension restrictions {\it can be
weakened to $m\ge n+3$} for highly connected PL $n$-manifolds (less highly
connected than in Theorems \wi6 and \wi8).

\proclaim{Theorem \wu5} [Sk02]
For embeddings $N\to\R^m$ of is a closed $d$-connected PL $n$-manifold $N$ and
$m\ge n+3$, the Haefliger-Wu invariant is
$$\text{bijective if}\quad 2m\ge3n+3-d\quad\text{and surjective if}
\quad 2m=3n+2-d.$$
\endproclaim

For $d=1$ we need only {\it homological} 1-connectedness in Theorem \wu5.

Observe that Theorem \wu5 is not quite the result expected in the 1960's,
and that its proof cannot be obtained by direct generalization of the
Haefliger-Weber proof without invention of new ideas.
This follows from the preceding discussion and the following two remarks.

(1) The PL case of Haefliger-Weber Theorem \wu4 holds for polyhedra,
but Theorem \wu5 holds only for highly enough connected PL manifolds.

(2) The same $(3n-2m+2)$-connectedness assumption as in the surjectivity part
of Theorem \wu5 ($2m\ge3n+2-d\Leftrightarrow d\ge3n-2m+2$) appeared in the
Hudson PL version of the Browder-Haefliger-Casson-Sullivan-Wall Embedding
Theorem for closed manifolds (roughly speaking, it states that a homotopy
equivalence between PL manifolds is homotopic to a PL embedding, and it was
proved by engulfing) [Hu67].
This assumption was soon proved to be superfluous (by surgery) [Ha68, Hu70'].
So it is natural to expect that the $(3n-2m+3)$-connectedness assumption in
Theorem \wu5 is superfluous (Theorem \wu5 is proved using generalization
of the engulfing approach).
However, Non-injectivity Examples \wu7.b,c of the next subsection show that
this assumption is {\it sharp}.

In this paper we sketch the proof of the surjectivity in the Haefliger-Weber
Theorem \wu4 in the PL case and present the idea of the proof of Theorem \wu5.

Most of the results of \Wi, \Lc\ and \Vk\ are corollaries of Theorems \wu4 and
\wu5, although some of them were originally proved independently (sometimes in
a weaker form).

\proclaim{Corollary \wu6} If $N$ is a homologically 1-connected closed smooth
$n$-manifold, then $\alpha_{DIFF}^m(N)$ is injective for $2m=3n+2$, $n=4s+2$,
and surjective for $2m=3n+1$, $n=4s+3$.
\endproclaim

\demo{Sketch of the proof}
Theorem \wu5 and smoothing theory [Ha67, 1.6, Ha, 11.1] imply the following.

{\it Let $N$ be a closed $d$-connected
(for $d=1$, just homologically 1-connected) smooth $n$-manifold and $m\ge
n+3$.

If $2m\ge3n+2-d$, then for each $\Phi\in\pi^{m-1}_{eq}(\t N)$
there is a PL embedding $f:N\to\R^m$ smooth outside a point and such that
$\alpha(f)=\Phi$; a complete obstruction to the smoothing of $f$ lies in
$C^{m-n}_{n-1}$.

If $2m\ge3n+3-d$, then any two smooth embeddings $f_0,f_1:N\to\R^m$
such that $\alpha(f_0)=\alpha(f_1)$ can be joined by a PL isotopy, which is
smooth outside a point; a complete obstruction to the smoothing of such a PL
isotopy lies in $C^{m-n}_n$.}

The Corollary follows from these assertions and $C^{2k}_{4k-2}=0$ [Ha66, 8.15,
Mi72', Corollary C]. (There is a misprint in [Ha66, 8.15]: instead of
'$C^{3k}_{4k-2}=0$' it should be '$C^{4k}_{8k-2}=0$'.)
\qed\enddemo

In Theorem \wu5 the surjectivity is not interesting for $m<\frac{5n+6}4$.
Indeed, $\frac{5n+6}4>\frac{3n+2-d}2$ implies that $d>\frac n2-1$ and $n\ge 6$;
hence $N$ is a homotopy sphere, so $N\cong S^n$ and the surjectivity in
Theorem \wu5 is trivial.
But the proof is not simpler for $m\ge\frac{5n+6}4$; the proof can also be
considered as a step towards the analogue of Theorem \wu5 for embeddings into
{\it manifolds}, which is interesting even for $m<\frac{5n+6}4$.
An analogous remark can be made about the injectivity in Theorem \wu5.

The Haefliger-Weber Theorem \wu4 has relative and approximative versions [Ha63,
1.7.2, We67,  Theorems 3 and 7, RS98], which require that constructed embedding
or isotopy extend given one or is close to given one.
But Theorem \wu5 has such versions only under some additional assumptions.

An interesting corollary of [Ku30, Cl34, Cl37] was deduced in [Wu65] for graphs
and in [Sk98] for the general case:

{\it a Peano continuum $N$ embeds into $\R^2$ if and only if there exists an
equivariant map $\t N\to S^1$.}

An interesting corollary of [MA41] was deduced in [Wu65]:

{\it embeddings $f,g:N\to\R^2$ of a Peano continuum $N$ are isotopic if and
only if $\t f\simeq_{eq}\t g$.}

\subhead The incompleteness of the Haefliger-Wu invariant \endsubhead
Clearly, $\t{S^n}\simeq_{eq} S^n$.
Therefore the Haefliger-Wu invariant is not injective in codimension 2 (e.g. for
knots in $\R^3$) and any smoothly non-trivial knot $S^n\to\R^m$
demonstrates the non-injectivity of $\alpha_{DIFF}^m(S^n)$.
The deleted product necessary conditions for embeddability or for isotopy
do not reflect either the ambience of isotopy or the
distinction between the DIFF and PL (or TOP) categories.
The same assertions hold for generalized Haefliger-Wu invariants (see below) or
'isovariant maps invariants' [Ha63, Ad93, \S7].

Let us present another examples.
All the examples in this subsection hold for {\it each} set
of the parameters $k,l,m,n,p$ satisfying to the conditions in the statement.

\proclaim{Non-injectivity Examples \wu7}
The following maps are not injective:
$$(a)\quad\alpha^{3l}(S^{2l-1}\sqcup S^{2l-1}\sqcup S^{2l-1}),\quad
\alpha^{3l}(S^{2l-1}\sqcup S^{2l-1})\quad\text{and}\quad
\alpha^{3l}_{DIFF}(S^{2l-1});$$
$$(b)\quad\alpha^{6k}(S^p\times S^{4k-1})\quad\text{for }p<k\quad\text{[Sk02]};
$$
$$(c)\quad\alpha^{3l+1}(S^1\times S^{2l-1}),\quad
\text{if $l+1$ is not a power of 2 [Sk]};$$
$$(d)\quad\alpha^m_{PL}((S^n\vee S^n)\sqcup S^{2m-2n-3})\quad\text{for}
\quad n+2\le m\le(3n+3)/2\quad\text{[Sk02]}.$$
\endproclaim

%[Ha62, Ha62', Ze62].
%0) $\alpha$ is not injective for each pair $(m,n)$ such that
%$n+3\le m\le\frac{3(n+1)}2$ but $m-n\not\in\{4,8\}$ and a polyhedron
%(for $m=\frac{3(n+1)}2$, non-connected PL manifold) $N=S^n\sqcup S^{2m-2n-3}$

The Non-surjectivity Example \wu7.a follows from the Borromean Rings, the
Whitehead Link and the Trefoil Knot Examples of \Lc.
The Non-injectivity Examples \wu7.b, \wu7.c and \wu7.d are constructed in
\Kt, below and in \Ex, respectively.

%Example \wu7.c shows that $\alpha^m_{PL}(N)$ can fail to be injective if
%$$\text{either}\quad d=-1\text{ and }2m=3n+2,\quad\text{or}\quad 2m=3n-3d,
%\quad\text{or}\quad d=0\text{ and }2m=3n-d+2.$$

The construction (but not the proof [Sk]) of Example \wu7.c is very simple
and explicit.

\smallskip
{\it Construction of Example \wu7.c in the PL case.}
Add a strip to the
Whitehead link $\omega_{0,PL}$, i.e. extend it to an embedding
$$\omega_0':S^0\times S^{2l-1}\bigcup\limits_{S^0\times D^{2l-1}_+=\partial
D^1_+\times D^{2l-1}_+} D^1_+\times D^{2l-1}_+\to\R^{3l}.$$ This
embedding contains connected sum of the components of the
Whitehead link. The union of $\omega_0'$ and the cone over the
connected sum forms an embedding $D^1_+\times
S^{2l-1}\to\R^{3l+1}_+$. This latter embedding can clearly be
shifted to a proper embedding. {\it The PL Whitehead torus}
$\omega_{1,PL}:S^1\times S^{2l-1}\to\R^{3l+1}$ is the union of
this proper embedding and its mirror image with respect to
$\R^{3l}\subset\R^{3l+1}$.
Cf. definition of $\mu'$ in \Kt.

\smallskip
It is easy to prove that $\alpha(\omega_{1,PL})=\alpha(f_0)$, where $f_0$ is
the standard embedding [Sk].
It is non-trivial that $\omega_{1,PL}$ is not PL isotopic to the standard
embedding when $l+1$ is not a power of 2.
The Non-injectivity Example \wu7.a for $\alpha^{3l}(S^{2l-1}\sqcup S^{2l-1})$
is based on the linking coefficient.
The Non-injectivity Example \wu7.c is much more complicated because
$S^1\times S^{2l-1}$
is connected,  so the linking coefficient cannot be defined
(the linking coefficient for the restriction to $S^0\times S^{2l-1}$ gives
the weaker Non-injectivity Example \wu7.b); thus a new invariant [Sk] is
required.

%{\it the Haefliger-Massey invariant} (for $S^{2l-1}\sqcup S^{2l-1}\sqcup S^{2l-1}$),
%{\it the Haefliger invariant} (for $S^{2l-1}$)
%Example \wu7.d shows that the dimension restriction is sharp in
%the polyhedral case of the Haefliger-Weber Theorem \wu4.

\proclaim{Non-surjectivity Examples \wu8}
For $m\ge n+3$ the following maps are not surjective:

(a) $\alpha^m(S^n\sqcup S^n)$, if
$\Sigma^\infty:\pi_n(S^{m-n-1})\to\pi_{2n+1-m}^S$ is not epimorphic, e.g. for
$$\minCDarrowwidth{0pt}\CD
n @=6               @=9               @=12            @=13     @=14  @=21\\
m @=10=\frac{3n+2}2 @=13=\frac{3n-1}2 @=18=\frac{3n}2 @=19=\frac{3n-1}2
            @=22=\frac{3n+2}2 @=31=\frac{3n-1}2
\endCD$$
(b) $\alpha^m(S^1\times S^{n-1})$, if $m-n$ is odd and
$\Sigma^\infty:\pi_{n-1}(S^{m-n})\to\pi^S_{2n-m-1}$ is not epimorphic [Sk02],
e.g. for
%if $n\in\{7,15\}$ and $2m=3n-1$ [Sko02].
$$\minCDarrowwidth{0pt}\CD
n @= 7               @= 10              @= 13              @= 14
@= 15              @= 22\\
m @= 10=\frac{3n-1}2 @= 13=\frac{3n-4}2 @= 18=\frac{3n-3}2 @= 19=\frac{3n-4}2
@= 22=\frac{3n-1}2 @= 31=\frac{3n-4}2
\endCD$$
$$(c)\quad\alpha_{DIFF}^{6k+1}(S^{2k}\times S^{2k}).$$
\endproclaim

For $n\le k+2$ the stable suspension mapping is denoted by
$\Sigma^\infty:\pi_{n+k}(S^n)\to\pi_{2k+2}(S^{k+2})=\pi_k^S$.

The Non-surjectivity Example \wu8.a follows from the formula
$\alpha=\pm\Sigma^\infty\lambda_{12}$ and the construction of a link with the
prescribed linking coefficient of \Lc.
The Non-surjectivity Example \wu8.b is constructed in \Kt.
%It shows that the dimension restriction in the surjectivity part of Theorem
%\wu5 is {\it almost sharp}: $\alpha(N)$ can fail to be surjective for $d=0$,
%$2m\le3n-1$.
The Non-surjectivity Example \wu8.c follows because
$\alpha_{PL}^{6k+1}(S^{2k}\times S^{2k})$ is bijective by Theorem \wu5 (or by
[Bo71, Ha62'']) but by the Haefliger Torus Example \lc5 there exists a PL
embedding $S^{2k}\times S^{2k}\to\R^{6k+1}$ that is not PL isotopic to a smooth
embedding.

Links give many other examples of the non-injectivity and the non-surjectivity
of $\alpha$.
From a link example, by gluing an arc joining connected components we can
obtain a highly connected polyhedral example.

\proclaim{Non-embeddability Examples \wu9} There exists an equivariant map
$\t N\to S^{m-1}$ but $N$ does not CAT embed into $\R^m$ (hence
$\alpha^m_{CAT}(N)$ is not surjective) for

(a) CAT=DIFF, $m=n+3$, $n\in\{8,9,10,16\}$ and a certain homotopy $n$-sphere
$N$;

(b) CAT=DIFF, $m=6k-1$, $n=4k$ and a certain (almost parallelizable
$(2k-1)$-connected) $n$-manifold $N$;

(c) CAT=PL, $\max\{3,n\}\le m\le\frac{3n+2}2$ and a certain $n$-polyhedron $N$
[MS67, SS92, FKT94, SSS98, GS06].
\endproclaim

The Non-Embeddability Example \wu9.a follows from the existence of a homotopy
$n$-sphere $N$ which is non-embeddable in codimension 3 [HLS65, Le65,
cf.\ Re90, \S2, MT95, pp. 407--408] (because $N\cong S^n$ topologically and so
$\t N\simeq_{eq}\t{S^n}$).
The Non-Embeddability Example \wu9.b follows from Theorems \wi7.c and \wi7.d.
The Non-Embeddability Example \wu9.c  is proved in \Ex.

%$\alpha$') $\alpha$ is not surjective for $N=S^1\times S^q$ and
%$m=q+l+2$, if $l\ge2$,
%$\Sigma^\infty_{(l)}:\pi_q(S^l)_{(l)}\to\pi^S_{q-l,(l)}$ is
%not epimorphic and $\Sigma^\infty:\pi_q(S^{l+1})\to\pi^S_{q-l-1}$ is
%monomorphic.
% @= 13 @= 11 @= 12\\
% @= 21 @= 17 @= 18\\

In [Sk98] it is proved that

{\it although the 3-adic solenoid $\Sigma$ (i.e.\ the intersection of infinite
sequence of filled-tori, each inscribed in the previous one with degree 3)
does not embed into  $\R^2$, nevertheless there exists an equivariant map
$\t\Sigma\to S^1$.}

We conjecture that there exists a non-planar {\it tree-like} continuum $N$,
for which there are no equivariant maps $\t N\to S^1$.

\subhead The Generalized Haefliger-Wu invariant \endsubhead
 The Borromean Rings Example \lc2 suggests that one can introduce
an obstruction to embeddability, analogous to the deleted
product obstruction (and the van Kampen obstruction, see \Vk)
but deduced from a triple, quadruple, etc.\ product.
Moreover, the vanishing of this obstruction should be sufficient for
embeddability even when this is not so for the deleted product obstruction.
Such an obstruction can indeed be constructed as follows, cf.\ [Kr00].
Suppose that $G$ is a
subgroup of the group $S_p$ of permutations of $p$ elements and let
$$\t N_G=\{(x_1,\dots,x_p)\in N^p\ |\ x_i\not=x_{\sigma(i)}
\text{ for each non-identity element }\sigma\in G,\ i=1,\dots,p\}.$$
The space $\t N_G$ is called the {\it deleted $G$-product} of $N$.
The group $G$ obviously acts on the space $\t N_G$.
For an embedding $f:N\to\R^m$ the map $\t f_G:\t N_G\to\t{\R^m}_G$ is
well-defined by the formula $\t f_G(x_1,\dots,x_p)=(fx_1,\dots,fx_p)$.
Clearly, the map $\t f_G$ is $G$-equivariant.
Then we can define the $G$-Haefliger-Wu invarant
$$\alpha_G=\alpha\phantom{}^m_G(N):\Emb\phantom{}^m(N)\to
[\t N_G,\t{\R^m}_G]_G \quad\text{by}\quad \alpha_G(f)=[\t f_G].$$
{\it The deleted $G$-product obstruction for the embeddability of $N$ in
$\R^m$} is the existence of a $G$-equivariant map $\Phi:\t N_G\to\t{\R^m}_G$.

This approach works well in link theory (the simplest example is the
classification of \lq higher-dimensional Borromean rings' [Ha62', \S3, Ma90,
Proposition 8.3] by means of $\alpha^m_{S_3}$).
Surprisingly, in contrast to that, the Non-injectivity Examples \wu7 (except
for the example of the Borromean rings of \wu7.a) demonstrate the
non-injectivity of $\alpha_G$ for each $G$: in their formulations $\alpha$ can
be replaced by $\alpha_G$ for each $G$.
This follows by the construction of these examples.
Clearly, if $\alpha$ is not surjective, then neither is $\alpha_G$ for
each $G$.
Under the conditions of the Non-embeddability Examples \wu9 property (e) is
false even if we replace the $\Z_2$-equivariant map $\t N\to S^{m-1}$ by a
$G$-equivariant map $\t N_G\to\t{\R^m}_G$.

\subhead Historical remarks \endsubhead
A particular case of Theorem \wu4 (Theorem \vk1) was discovered by Van Kampen
[Ka32].
But Van Kampen's proof of the sufficiency in Theorem \vk1 contained a mistake.
However, he modified his argument to prove the PL case of Theorem \wi2.a.
Based on the idea of the Whitney trick invented in [Wh44],
Shapiro and Wu completed the proof [Sh57, Wu58].  Subsequently
their argument was generalized by Haefliger and Weber (using some ideas of
Shapiro and Zeeman) in order to prove the Haefliger-Weber Theorem \wu4.
The second part of the Weber proof was simplified in [Sk98] using the
idea of the Freedman-Krushkal-Teichner Lemma \vk3.

The Whitney trick, on which the proof of sufficiency in Theorem \vk1 for
$n\ge3$ is based, cannot be performed for $n=2$ [KM61, La96].
Sarkaria has found a proof of the case $n=1$ of Theorem \vk1 based on
1-dimensional Whitney trick [Sa91'] (the author is grateful to K. Sarkaria and
M. Skopenkov for indicating that the argument in [Sa91'] is incomplete).
Sarkaria also asked whether the sufficiency in Theorem \vk1 holds for the case
$n=2$.
Freedman, Krushkal and Teichner have constructed an example showing that it
does not [FKT94].

The dimension restriction $2m\ge3n+3$ in the Haefliger-Weber
Theorem \wu4 comes from the use of the Freudenthal Suspension
Theorem, the Penrose--Whitehead--Zeeman--Irwin Embedding Theorem
\wi6.c, a relative version of the Zeeman Unknotting Theorem \wi5.a
and  general position arguments (\Web). Toru\'nczyk and Spie\D z
showed that one can try to relax the restriction coming from the
Zeeman Unknotting Theorem (using relative regular neighborhoods)
and those coming from the Freudenthal Suspension Theorem (using
Whitehead's `hard part' of the Freudenthal Suspension Theorem and
the Whitehead higher-dimensional finger moves [FF89, \S10]; note
that application of higher-dimensional finger moves in this
situation was first suggested by Schepin) [Sp90, ST91, see also
DRS91, DRS93]. This was the reason why in 1992 Dranishnikov
conjectured the surjectivity in the Haefliger-Weber Theorem \wu4
for $2m=3n+2$. However, Segal and Spie\D z constructed a
counterexample (using the same higher-dimensional finger moves)
--- a weaker version of Non-Embeddability Example \wu9.c [SS92].
Their construction  used the Adams theorem on Hopf invariant one,
and therefore has exceptions corresponding to the exceptional
values $1,3,7$.  In 1995 the author suggested how to remove the
codimension 2 exceptions.  Subsequently, this idea was generalized
independently by Segal--Spie\D z and the author to obtain a
simplification of [SS92] which did not use the Adams theorem, and
therefore has no exceptions [SSS98].  This simplified construction
leads to Non-Embeddability Example \wu9.c.

%\head II. Proofs of embedding theorems \endhead Clearly, for each $f$ the map
%$\t f:\t N-\t\Delta(f)\to S^{m-1}$ is well-defined.
%\newpage

\head \kt\ On the deleted product of the torus \endhead

\subhead Proof of the Deleted Product Lemma \wu3.a,b \endsubhead

\bigskip
\centerline{\epsffile{6-0.eps}}
\centerline{\it Figure \kt0} %approximately here}
\bigskip

\demo{Proof of the Deleted Product Lemma \wu3.a}
We have
$$\t{S^p\sqcup S^q}\simeq_{eq}
S^p\times S^q\sqcup S^p\times S^q\sqcup S^p\sqcup S^q,$$
where the involution on the right-hand term exchanges antipodes
in $S^p$ and in $S^q$, as well as the corresponding points from the
two copies of $S^p\times S^q$ (Figure \kt0).
Therefore analogously to the definition of the $\alpha$-invariant in \Lc\
we have
$$\pi^{m-1}_{eq}(\t{S^p\sqcup S^q})\cong [S^p\times
S^q,S^{m-1}]\cong\pi_{p+q+1-m}^S\quad\text{for} \quad m-2\ge p,q.\qed$$
\enddemo

Recall that the {\it equivariant Stiefel manifold} $V^{eq}_{mn}$ is the space
of equivariant maps $S^{n-1}\to S^{m-1}$.
Denote by $a_k:S^k\to S^k$ the antipodal map.

\proclaim{The Torus Lemma \kt1}
For $p\le q$ and $m\ge p+q+3$ there exist groups $\Pi^{m-1}_{p,q}$ and
$\Pi^{m-1}_{q,p}$, a group structure on $\pi^{m-1}_{eq}(\t{S^p\times S^q})$,
and homomorphisms $\sigma$, $\gamma$, $\rho$, $\tau$ and $\alpha'$ forming
the diagram below, in which the right-hand square commutes and the
left-hand square either commutes or anticommutes.
The homomorphisms
$\sigma$, $\gamma$, $\rho$ and $\pr_1$ are isomorphisms under the
dimension restriction $m\ge A$, where $A$ is shown near the notation of a
map.  $$\minCDarrowwidth{12pt}\CD \pi_q(V_{m-q,p+1}) @>> \tau\
\frac{3q}2+p+2 > KT^m_{p,q} @>> \alpha > \pi^{m-1}_{eq}(\t{S^p\times
S^q})\\ @VV \rho\ \frac{3q}2+p+2 V @VV \alpha' V @VV \gamma\ p+q+3 V \\
\pi_q(V^{eq}_{m-q,p+1}) @>> \sigma\ \frac{3q+p}2+2 > \Pi^{m-1}_{pq} @<
\pr_1\ q+2p+2 << \Pi^{m-1}_{p,q}\oplus\Pi^{m-1}_{q,p} \endCD$$
\endproclaim

%, the left-hand square (anti)commute for $2m\ge3q+2p+3$ and

\bigskip
\centerline{\epsffile{6-1.eps}}
\centerline{\it Figure \kt1} %approximately here}
\bigskip

\demo{Proof}
There is an equivariant deformation retraction
$$H_1:\t{S^p\times S^q}\to\adiag S^p\times S^q\times S^q
\bigcup\limits_{\adiag S^p\times\adiag S^q}S^p\times S^p\times\adiag S^q,$$
where $\adiag$ is antidagonal (Figure \kt1).
More precisely, for non-antipodal points $x$ and $y$ of $S^k$ and
$t\in[0,1]$, let $[x,y,t]=\frac{(1-t)x+ty}{|(1-t)x+ty|}$.
Define the deformation $H_t:\t{S^p\times S^q}\to\t{S^p\times S^q}$ by
$$H(x,y,x_1,y_1)=\cases
([x,y,t],[y,x,t],[x_1,y_1,2\delta t],[y_1,x_1,2\delta t])
& |x,y|\ge|x_1,y_1|\\
([x,y,2(1-\delta)t],[y,x,2(1-\delta)t],[x_1,y_1,t],[y_1,x_1,t])
& |x_1,y_1|\ge|x,y|\endcases,$$
where $\delta=\frac{|x,y|}{|x_1,y_1|+|x,y|}$.
%\quad\text{and}\quad\delta_1=1-\delta=\frac{|x_1,y_1|}{|x_1,y_1|+|x,y|}.$$

Define $\Pi^{m-1}_{pq}:=\pi^{m-1}_{eq}(S^p\times S^{2q})$, where the
involution on $S^p\times S^{2q}$ is $a_p\times t_q$ and
$t_q:S^{2q}\to S^{2q}$ is the symmetry with respect to $S^q\subset S^{2q}$.
The group structure on $\Pi^{m-1}_{pq}$ is defined as follows.
For equivariant maps $\varphi,\psi:S^p\times S^{2q}\to S^{m-1}$ define the
map $\varphi+\psi:S^p\times S^{2q}\to S^{m-1}$ on $x\times S^{2q}$ to be the
ordinary sum of the restrictions of $\varphi$ and $\psi$ to $x\times S^{2q}$.
Analogously, define the unity and the inverse of $\varphi$ on
$x\times S^{2q}$ to be the ordinary unity and the ordinary
inverse of $\varphi|_{x\times S^{2q}}$.

Let $v_q:S^q\times S^q\to\frac{S^q\times S^q}{S^q\vee S^q}\cong S^{2q}$
be the quotient map, cf. \Lc.
Consider the involution $(s,x,y)\to(-s,y,x)$ on $S^p\times S^q\times S^q$.
For $m\ge p+q+3$ by general position we have a 1--1 correspondence
$$(\id\phantom{}_{S^p}\times v_q)^*:
\pi^{m-1}_{eq}(S^p\times S^q\times S^q)\cong\Pi^{m-1}_{pq}.$$
One can check that the involution on $S^q\times S^q$ exchanging
factors corresponds to $t_q$.

Consider the restrictions of an equivariant map
$$\t{S^p\times S^q}\to S^{m-1}\quad\text{to}\quad
\adiag S^p\times S^q\times S^q\quad\text{ and to}
\quad S^p\times S^p\times\adiag S^q.$$
Define the map $\gamma$ to be the direct sum of the compositions of such
restrictions and the isomorphisms $(\id_{S^p}\times v_q)^*$ and
$(\id_{S^q}\times v_p)^*$.
If
$$\dim(\adiag S^p\times\adiag S^q)=p+q\le(m-1)-2,$$
then $\gamma$ is a 1--1 correspondence by general position and
the Borsuk Homotopy Extension Theorem.
Take the group
structure on $\pi^{m-1}_{eq}(\t{S^p\times S^q})$ induced by $\gamma$.
Then $\gamma$ is an isomorphism.

%(This is a cohomotopy analogue of a partial case of the Mayer--Vietoris
%sequence; if $\dim(S^p\times S^p\times\adiag S^q)=2p+q\le (m-1)-1$, then the
%proof is simplified).

By general position, for $2p+q\le m-2$ we have $\Pi^{m-1}_{qp}=0$, hence
$\pr_1$ is an isomorphism.

Let $\alpha'$ be the map corresponding under the isomorphism
$(\id_{S^p}\times v_q)^*$ to the map
$$S^p\times S^q\times S^q\to S^{m-1}\quad\text{defined by}
\quad(s,x,y)\mapsto\t f((s,x),(-s,y)).$$

\bigskip
\centerline{\epsffile{6-2.eps}}
\centerline{\it Figure \kt2} %approximately here}
\bigskip

Clearly, the right-hand square of the diagram commutes.

The map $\tau$ was defined after Theorem \lc9.

Recall that $\rho$ is the inclusion-induced homomorphism.
By [HH62, Sk], $\rho$ is an isomorphism for $m\ge\frac{3q}2+p+2$.

An element $\varphi\in\pi_q(V^{eq}_{m-q,p+1})$ can be considered
as a map $\varphi:S^p\times S^q\to S^{m-q-1}$ such that
$\varphi(-x,y)=-\varphi(x,y)$ for each $x\in S^p$.
Let $\sigma(\varphi)$ be the $q$-fold $S^p$-fiberwise suspension
of such a map $\varphi$, i.e.\  $\sigma(\varphi)|_{x\times
S^{2q}}=\Sigma^q(\varphi|_{x\times S^q})$.  It is easy to see that
$\sigma$ is a homomorphism.

The (anti)commutativity of the left-hand square is proved for $p=0$ using the
representation $S^{m-1}\cong S^{m-q-1}*S^{q-1}$ and deforming
$\alpha'\tau(\varphi)$ to the $S^p$-fiberwise suspension $\sigma(\varphi)$
of $\varphi$ [Ke59].
For $p>0$ we apply this deformation for each $x\in S^p$ independently.

It remains to prove that $\sigma$ is an isomorphism
for $m\ge\frac{3q+p}2+2\ge p+q+2$.
For this observe that $\sigma$ is a composition
$$\pi_q(V^{eq}_{m-q,p+1})=\pi^{m-q-1}_{eq}(S^p\times S^q)
\overset{\Sigma^q}\to\to\pi^{m-1}_{eq}(\Sigma^q(S^p\times S^q))
\overset{\pr^*}\to\to\Pi^{m-1}_{pq}$$
Here the involution on $S^p\times S^q$ is $a_p\times\id_{S^q}$,
the involutions on
 $\Sigma^qS^q$ and on $\Sigma^q(S^p\times S^q)$ are the \lq suspension'
involutions over $\id_{S^q}$ and $a_p\times\id_{S^q}$;
the map
$$\pr:S^p\times\Sigma^q S^q=
S^p\times\frac{S^q\times D^q}{S^q\times y,\ y\in\partial D^q}\to
\frac{S^p\times S^q\times D^q}{S^p\times S^q\times y,\ y\in\partial D^q}=
\Sigma^q(S^p\times S^q)$$
is a quotient map (Figure \kt3).
The $S^p$-fiberwise group structures on $\pi_{eq}^{m-q-1}(S^p\times S^q)$ and
on $\pi_{eq}^{m-1}(\Sigma^q(S^p\times S^q))$ are defined analogously to that
on $\Pi^{m-1}_{p,q}$.
By the equivariant version of the Freudenthal Suspension Theorem,
it follows that the above $\Sigma^q$ is an isomorphism for $p+q\le2(m-q-1)-2$.
The non-trivial preimages of $\pr$ are $S^p\times[S^q\times y]$,
$y\in\partial D^q$.
Their union is homeomorphic to $S^p\times\partial D^q$.
Since $\dim(S^p\times\partial D^q)=p+q-1$, by general position
it follows that $\pr^*$ is an isomorphism for $p+q-1\le m-3$.
Therefore $\sigma$ is an isomorphism for $m\ge\frac{3q+p}2+2\ge p+q+2$.
\qed\enddemo

\bigskip
\centerline{\epsffile{6-3.eps}}
\centerline{\it Figure \kt3} %approximately here}
\bigskip

Note that the above map $(\id\phantom{}_{S^p}\times v_q)^*$ is an
isomorphism also for $m=p+q+2$.
This is proved using the cofibration exact sequence of the pair
$(S^p\times S^q\times S^q,S^p\times(S^q\vee S^q))$ and a retraction
$$\Sigma(S^p\times S^q\times S^q)\to\Sigma(S^p\times(S^q\vee S^q))$$
obtained from the retraction $$\id\phantom{}_{S^p}\times
r_q:S^p\times\Sigma(S^q\times S^q)\to S^p\times\Sigma(S^q\vee S^q)$$ by
shrinking the product of $S^p$ with the vertex of the
suspension to a point.

%By Theorems \wu2 and \wu5 the maps $\alpha_{DIFF}$ and $\alpha_{PL}$ are
%bijective for $m\ge\frac{3(q+p)}2+2$ and $m\ge\frac{3q}2+p+2$, respectively.
%From the existence of $\tau$ and by Torus Lemma \kt1 it follows that
%$\alpha_{DIFF}$ is surjective for $m\ge\max(\frac{3q}2+p+2,q+2p+2)$ [cf.
%Boe71, BH70].
%Note that some maps of Torus Lemma \kt1 are {\it epimorphisms}
%under the weaker by one dimension restrictions than stated there.
%The {\it PL Stiefel manifold} $V_{mn}^{PL}$ is the space of
%PL embeddings $S^{n-1}\to S^{m-1}$.
%Replacing $V\to V^{PL}$, we can define analogously a map
%$\tau^{PL}:\pi_q(V^{PL}_{m-q,p+1})\to\Emb^m_{PL}(S^p\times S^q)$.
%$\alpha'$ and $\tau$ are bijections for $m\ge\max(\frac{3q}2+p+2,q+2p+2)$
%\{$\frac{3(q+p)}2+2$\},

\proclaim{The Generalized Torus Lemma \kt2} [Sk02] (a) If
$$s\ge3, \quad p_1\le\dots\le p_s,\quad n=p_1+\dots+p_s\quad\text{and}
\quad N=S^{p_1}\times\dots\times S^{p_s},$$
then the same assertion  as in the Torus Lemma \kt1 holds for the following
diagram:
$$\minCDarrowwidth{12pt}\CD \pi_{n-p_1}(V_{m-n+p_1,p_1+1}) @>> \tau >
\Emb^m(N) @>> \alpha >
\pi^{m-1}_{eq}(\t N)\\
@VV \rho\ \frac{3n-p_1}2+2 V @VV \alpha' V @VV \gamma\ 2n-p_1-p_2+3 V \\
\pi_{n-p_1}(V^{eq}_{m-n+p_1,p_1+1}) @>> \sigma\ \frac{3n}2-p_1+2 >
\Pi^{m-1}_{p_1,n-p_1} @< \pr_1\ 2n-p_2+2 <<
\oplus_i\Pi^{m-1}_{p_i,n-p_i}
\endCD$$

(b) Suppose that
$$s\ge3,\quad n=p_1+\dots+p_s,\quad p_1\le\dots\le p_s\quad\text{and}
\quad m\ge2n-p_1-p_2+3$$
(for $s=3$ and CAT=DIFF assume also that $m\ge\frac{3n}2+2$). Then
$$\Emb\phantom{}^m(S^{p_1}\times\dots\times S^{p_s})=
\oplus_{i=1}^s\pi_{n-p_i}(V_{m-n+p_i,p_i+1}).$$
\endproclaim

\demo{Proof} Part (b) follows from part (a) and Theorems \wu4 and \wu5.
The proof of (a) is analogous to the proof of the Torus Lemma \kt1.
We shall only define $\tau$ and $\sigma$ and omit the details.

The map $\tau$ is defined as follows.  An
element $\varphi\in\pi_{n-p_1}(V_{m-n+p_1,p_1+1})$ is
represented by a map $S^{n-p_1}\times S^{p_1}\to S^{m-n+p_1-1}$.
Consider the projections
$$\pr\phantom{}_1:N\to S^{p_1}\times S^{p_2+\dots+p_s}=S^{p_1}\times S^{n-p_1}
\quad\text{and}\quad\pr\phantom{}_2:N\to S^{p_2}\times\dots\times S^{p_s}.$$
Analogously to the case $s=2$, define an embedding $\tau(\varphi)$
as
the composition
$$S^{p_1}\times S^{p_2}\times\dots\times S^{p_s}
\overset{(\varphi\circ\pr_1)\times\pr_2}\to\to
\partial D^{m-n+p_1}\times S^{p_2}\times\dots\times S^{p_s}\subset\R^m.$$
The map $\sigma$ is defined analogously to the case $s=2$ as the
$S^{p_1}$-fiberwise $(n-p_1)$-fold suspension.

Denote $q_1:=n-p_1$.
Equivalently, $\sigma$ is a composition
$$\pi_{q_1}(V^{eq}_{m-q_1,p_1+1})=
\pi^{m-q_1-1}_{eq}(S^{p_1}\times S^{q_1})\overset{\Sigma^{q_1}}\to\to
\pi^{m-1}_{eq}(\Sigma^{q_1}(S^{p_1}\times S^{q_1}))\overset{\pr^*}\to\to
\Pi^{m-1}_{{p_1,q_1}}.$$
Here the maps $\Sigma^{q_1}$ and $\pr^*$ are isomorphisms for
$2m\ge3n-2p_1+4$ and $m\ge n+2$, respectively.
Therefore $\sigma$ is an isomorphism for $m\ge3n/2-p_1+2$.
\qed\enddemo

Note that under the assumptions of the Generalized Torus Lemma \kt2.b we have
$m\ge\frac{3n}2+2$ for $s\ge4$.

\newpage
\subhead Proof of the Non-surjectivity Example \wu7.b  and
the Non-injectivity Example \wu8.b \endsubhead

\proclaim{The Decomposition Lemma \kt3}
[Sk02] For $m\ge2p+q+1\ge q+3$ there is the following (anti)commutative
diagram, in which the first and the third lines are exact.
The map $\nu$ is epimorphic for $m-q$ even and $\im\nu$ is the subgroup of
elements of order 2 for $m-q$ odd.
$$\minCDarrowwidth{12pt}\CD
\pi_q(V_{m-q-1,p}) @>>\mu''> \pi_q(V_{m-q,p+1}) @>>\nu''> \pi_q(V_{m-q,1})\\
@VV \tau_{p-1} V @VV \tau V @VV = V \\
\Emb^{m-1}_{PL}(S^{p-1}\times S^q) @>> \mu' >\Emb^m_{PL}(S^p\times S^q) @>>
\nu' > \pi_q(S^{m-q-1})\\
@VV \alpha'_{p-1} V @VV \alpha' V @VV \Sigma^\infty V \\
\Pi^{m-2}_{p-1,q} @>> \mu > \Pi^{m-1}_{p,q} @>> \nu > \pi^S_{2q+1-m}
\endCD.$$
\endproclaim

Of this Lemma only the right squares and the exactness at $\Pi^{m-1}_{pq}$ are
used for the examples.
The left squares are interesting in themselves and are useful elsewhere [Sk].
The definition of $\mu'$ here is simpler than that in [Sk02].

\demo{Definition of the maps from the diagram}
Let $\nu''$ and $\mu''$ be the homomorphisms induced by the
\lq forgetting the first $p$ vectors' bundle
$V_{m-q,p+1}\overset{V_{m-q-1,p}}\to\to V_{m-q,1}$.

For an embedding
$f:S^p\times S^q\to\R^m$ let $\nu'(f)$ be the linking coefficient of
$f|_{x\times S^q}$ and $f|_{-x\times S^q}$ in $\R^m$.

Define the map
$\nu:\Pi^{m-1}_{p,q}\to\Pi^{m-1}_{0,q}\cong\pi^S_{2q-m+1}$
as \lq the restriction over $*\times S^{2q}$'.

In order to define $\mu'$ denote by
$S^p=D^p_+\bigcup\limits_{\partial D^p_+=S^{p-1}=\partial D^p_-}D^p_-$
the standard decomposition of $S^p$.
Define analogously $\R^m_\pm$ and $\R^{m-1}$.
Add a strip to an embedding $f:S^{p-1}\times S^q\to\R^{m-1}$, i.e. extend it to
an embedding
$$f':S^{p-1}\times S^q\bigcup\limits_{S^{p-1}\times D^q_+=\partial
D^p_+\times D^q_+} D^p_+\times D^q_+\to\R^m.$$
Since $m\ge2p+q+1$, it follows that this extension is unique up to isotopy.
The union of $f'$ and the cone over the restriction of $f'$ to the boundary
forms an embedding $D^p_+\times S^q\to\R^m_+$.
This latter embedding can clearly be shifted to a proper embedding.
Define $\mu'f:S^p\times S^q\to\R^m$ to be the union of this proper embedding
and its mirror image with respect to $\R^{m-1}\subset\R^m$.

Let us define the map $\mu$ first for the case $p=1$.
For a map $\varphi:S^{2q}\to S^{m-2}$ define the map $\mu\varphi$ to
be the equivariant extension of the composition
$D^1\times S^{2q}\overset{\pr}\to\to\Sigma S^{2q}\overset{\Sigma\varphi}
\to\to S^{m-1}$.
In order to define the map $\mu$ for arbitrary $p$, replace
$$\Pi^{m-2}_{p-1,q}\quad\text{and}\quad\Pi^{m-1}_{pq}\quad\text{by}
\quad\pi^{m-2}_{eq}(\Sigma^q(S^{p-1}\times S^q))\quad\text{and}\quad
\pi^{m-1}_{eq}(\Sigma^q(S^p\times S^q)),$$ respectively (see the the proof of
Torus Lemma \kt1). For an equivariant map
$\varphi:\Sigma^q(S^{p-1}\times S^q)\to S^{m-2}$ let $\mu\varphi$ be the
composition
$$\Sigma^q(S^p\times S^q)=\Sigma^q(\Sigma S^{p-1}\times S^q)
\overset{\Sigma^q\pr}\to\to\Sigma^{q+1}(S^{p-1}\times S^q)
\overset{\Sigma\varphi}\to\to\Sigma S^{m-2},$$
where $\pr$ is the map from the proof of the Torus Lemma \kt1.
Clearly, the definition for arbitrary $p$ agrees with that for $p=1$.
\enddemo

Note that all the maps of the Decomposition Lemma \kt3 except $\mu'$
are defined for $m\ge p+q+3$.

The embedding $\mu'f$ can also be defined by the
Penrose-Whitehead-Zeeman-Irwin Embedding Theorem \wi6.c and its
isotopy analogue. For $m\ge2p+q+2$ any embedding $f:S^{p-1}\times
S^q\to S^{m-1}$ can be extended to a PL embedding
$f_\pm:D^p_\pm\times S^q\to\R^m_\pm$, uniquely up to isotopy. Then
two embeddings $f_+$ and $f_-$ define an embedding
$\mu'(f):S^p\times S^q\to\R^m$.

\demo{Proof of the Decomposition Lemma \kt3}
It is easy to check that both $\nu$ and $\mu$ are homomorphisms.

Clearly, the left-upper square of the diagram commutes.

Clearly, the right-upper square of the diagram commutes, see details in [Ti69,
Lemma 3].

The right-bottom square of the diagram (anti)commutes by [Ke59, Lemma 5.1].

We prove the commutativity of the left-bottom square for $p=1$; the proof is
analogous for the general case.
Take an embedding $f:S^0\times S^q\to\R^{m-1}$.
Then $\alpha\mu'f=\mu\alpha'_0f$ on $S^0\times S^{2q}$.
Also, for each $y\in S^1\times S^q\times S^q$ the points $(\alpha'\mu'f)y$ and
$(\mu\alpha_0'f)y$ are either both in the upper or both in the lower hemisphere
of $S^{m-1}$.
Hence $\alpha'\mu'f\simeq_{eq}\mu\alpha_0'f$.

Let us prove the exactness at $\Pi^{m-1}_{pq}$.
Clearly, $\nu\mu=0$.
On the other hand, if $\Phi:S^p\times S^{2q}\to S^{m-1}$ is an
equivariant map such that $\Phi|_{*\times S^{2q}}$ is null-homotopic,
then by the Borsuk Homotopy Extension Theorem, $\Phi$ is equivariantly
homotopic to a map which maps $*\times S^{2q}$ and
$a_p(*)\times S^{2q}$ to antipodal points of $S^{m-1}$.
By the Equivariant Suspension Theorem, the latter map is in $\im\mu$, since
$p-1+2q\le2(m-2)-1$. So $\ker\nu=\im\mu$.

Clearly, $\im\nu$ consists of homotopy classes $\varphi\in\Pi_{0,q}^{m-1}$
extendable to a map $D^1\times S^{2q}\to S^{m-1}$.
Such maps $\varphi$, considered as maps $\varphi:S^{2q}\to S^{m-1}$ are
exactly those which satisfy $a_{m-1}\circ\varphi\circ t_q\simeq\varphi$.
The latter condition is equivalent to $(-1)^m\varphi=(-1)^q\varphi$ (for $m$
odd this follows by [Po85, complement to lecture 6, (10), p.264],
since $h_0:\pi_{2q}(S^{m-1})\to\pi_{2q}(S^{2m-3})$ and $2q<2m-3$).
So $\im\nu=\ker(1-(-1)^{m-q})$.
\qed\enddemo

\demo{Proof of the Non-surjectivity Example \wu8.b} Set $q=n-1\le m-4$.
Look at the right-bottom square of the diagram from the Decomposition Lemma
\kt3 and use the surjectivity of $\nu$ for $m-q$ even.
The specific examples can be found using [To64, \S14]
(set $l=m-n=m-q-1$ and $k=2q+1-m$).
\qed\enddemo

\demo{Proof of the Non-injectivity Example \wu7.b}
Since $p<k$, we have $m\ge2p+q+2$.
Look at the right squares of the diagram from the Decomposition Lemma \kt3 and
use Lemma \kt4 below.
\qed\enddemo

\proclaim{Lemma \kt4} (a) $\Pi_{pq}^{m-1}$ is finite when $p+q+2\le m\le2q$.
%or when $m=2q+p+1$, $q$ odd $\ge3$ and $p\ge1$.

(b) The image of the restriction-induced homomorphism
$\nu_p'':\pi_{4k-1}(V_{2k+1,p+1})\to\pi_{4k-1}(S^{2k})$ is infinite for
$p<2k$.
\endproclaim

\demo{Proof}
Let us prove (a) by induction on $p$.
%For the case $p+q+2\le m\le2q$
The  base of the induction is $p=0$, when
$\Pi_{0q}^{m-1}\cong\pi_{2q}(S^{m-1})$ is indeed finite.
%For the case $m=2q+p+1$, $q$ odd $\ge3$ and $p\ge1$ the induction
%base is $p=1$, when $\Pi_{1q}^{2q+1}\cong\Z_{(q+2)}$ is indeed
%finite by Corollary \kt4.
The inductive step of (a) follows by the induction hypothesis and
the exactness of the bottom line from Decomposition Lemma \kt3.

In order to prove (b) for $p=0$ observe that the map $\nu_0$ is an
isomorphism and $\pi_{4k-1}(S^{2k})$ is infinite.
Suppose that $p\ge1$ and there is an infinite set
$\{x_i\}\in\pi_{4k-1}(V_{2k+1,p})$ with distinct $\nu_{p-1}''$-images.
Consider the Serre fibration $S^{2k+1-p}\to V_{2k+1,p+1}\overset\psi\to
\to V_{2k+1,p}$ and
the following segment of its exact sequence:
$$\pi_{4k-1}(V_{2k+1,p+1})\overset{\psi_*}\to\to\pi_{4k-1}(V_{2k+1,p})\to
\pi_{4k-2}(S^{2k-p}).$$
Since $\pi_{4k-2}(S^{2k-p})$ is finite, by exactness it follows that the
number of congruence classes of $\pi_{4k-1}(V_{2k+1,p})$ modulo $\im\psi_*$
is finite.
Therefore an infinite number of the $x_i$ (we may assume all the $x_i$) lie in
the same congruence class.
By passing from $\{x_i\}$ to $\{x_i-x_1\}$ we may assume that this congruence
class is the subgroup $\im\psi_*$ itself.
Hence the inductive step follows from $\nu_p''=\nu_{p-1}''\psi_*$.
\qed\enddemo

%For a group $G$ let $G_{(k)}=G$ for $k$ even and let $G_{(k)}=G/2G$ for $k$
%odd. If $G$ is finite abelian, then
%$G_{(k)}\cong G_{[k]}\cong G\otimes\Z_{(k)}$.
For a group $G$ let $G_{[k]}=G$ for $k$ even and let $G_{[k]}$ be the
subgroup of $G$ formed by elements of order 2 for $k$ odd.
For $m-q$ even, from the existence of a section
$s:\pi_q(V_{m-q,1})\to\pi_q(V_{m-q,2})$ it follows that $\nu''$ is
epimorphic and hence $\nu'$ is epimorphic.
%In the proof of assertion (*) above it was essentially shown that
We also have $\im\nu=\im\nu'=\im\nu''=\pi_{2q-m+1,[m-q]}$ for $2m\ge3q+4$.
Note that $\im\nu=\pi_{2q-m+1,[m-q]}$ even for $2m\le3q+3$ but by Lemma
\kt4.b $\im\nu''\ne\pi_q(S^{m-q-1})_{[m-q]}$ for $2m\le3q+3$.

We conjecture that if $m-q$ is even $\ge4$ and $2m\ge3q+4$, then
$\Pi^{m-1}_{1,q}\cong\pi^S_{2q-m+2}\oplus\pi^S_{2q-m+1}$.
We conjecture that in general $\Pi^{m-1}_{1q}$ is adjoint to
$\pi^S_{2q-m+2,[m-q]}\oplus\pi^S_{2q-m+1,[m-q]}$, unless $m=2q+1$
and $q=2l$ is even, when $\Pi^{4l}_{1,2l}\cong\Z_2$ (cf.\ the formula for
$\pi_q(V_{m-q,2})$ before Theorem \lc10).
Since $\im\nu=\pi^S_{2q-m+1,[m-q]}$, by the Decomposition Lemma \kt3 the
conjecture would follow from $\coim\mu\cong\pi^S_{2q-m+2,[m-q]}$
(recall that we identify $\Pi^{m-2}_{0q}=\pi^S_{2q-m+2}$).

\comment

\demo{Proof} Recall that $V_{m-q,2}\cong TS^{m-q-1}$.
For $m-q$ even we even have a section
$s:S^{m-q-1}\to V_{m-q,2}$ such that $\nu''s_*=\id$ (note that the map
$\tau s_*$ is a generalization of [Zee62, Example] for $q-r=1$).
If $2m\ge3q+4$, then $\Sigma^\infty:\pi_q(S^{m-q-1})\to\pi^S_{2q+1-m}$ is an
isomorphism.
Hence there is a section
$$\alpha\tau s_*(\Sigma^\infty)^{-1}:\pi^S_{2q+1-m}\to\Pi^{m-1}_{1q},\quad
\text{so}\quad\Pi^{m-1}_{1q}\cong\pi^S_{2q-m+2}\oplus\pi^S_{2q-m+1}.\qed$$

???We need a section from the other side; we need $\ker\mu=0$.
\enddemo

Many other concrete examples of non-surjectivity of $\alpha$ can be constructed
using [Tod62] and the following example.

\proclaim{Non-surjectivity Example \kt5}
$\alpha^m(S^1\times S^q)$ is not surjective if $m-q$ is odd $\ge5$ and
$\Sigma^{\infty}_{[1]}:\pi_q(S^{m-q-1})_{[1]}\to\pi^S_{2q+1-m,[1]}$
is not epimorphic.???
\endproclaim

\demo{Proof} It suffices to prove that
$$(*)\quad\im\nu=\pi^S_{2q-m+1,[1]}\quad\text{and}\quad
\ker2\subset\im\nu''\subset\im\nu'?\subset?\ker(2\Sigma^\infty)\quad\text{for }
m-q\text{ odd }\ge5.$$
Here by 2 is denoted the multiplication by 2 in the group $\pi_q(S^{m-q-1})$.

From the exact sequence of the bundle
$S^{m-q-2}\to V_{m-q,2}\to S^{m-q-1}$ it follows that
$$\ker(1-(-1)^{m-q})\subset\im\nu''\subset\ker[(1-(-1)^{m-q})\Sigma^\infty].$$
Since $\nu''=\nu'\tau$, it follows that $\im\nu''\subset\im\nu'$.
If $f:S^1\times S^q\to\R^m$ is an embedding, then the link $f|_{S^0\times S^q}$
and the link obtained by interchanging the components are isotopic.
Then
$$\Sigma^\infty(-1)^{m-q}\nu'(f)=\Sigma^\infty\nu'(f)\quad\text{i.e.}\quad
\nu'(f)\in\ker[(1-(-1)^{m-q})\Sigma^\infty],$$
and the second part of assertion (*) follows.
\qed\enddemo

In order to calculate $\ker\mu$ denote by $M$ the composition
$$S^1\times\Sigma^qS^q
\overset{\pr}\to\to\Sigma^q(S^1\times S^q)\cong\Sigma^q(\Sigma S^0\times S^q)
\overset{\Sigma^q\pr}\to\to\Sigma^{q+1}(S^0\times S^q),$$
where $\pr$ and $\pr$ are maps from the proof of Torus Lemma \kt1
(we use the same notation $\pr$ for two distinct maps).
Each map $\psi:S^{2q}\to S^{m-2}$ can be identified with an
equivariant map $\psi:\Sigma^q(S^0\times S^q)\to S^{m-2}$.
For each map $\psi:S^{2q}\to S^{m-2}$ we can construct an equivariant
map $h:\partial I^2\times S^{2q}\to S^{m-1}$ such that
$h|_{\partial I\times I\times S^{2q}}$ \lq represents' $M\circ*$,
\ $h|_{I\times0\times S^{2q}}$ \lq represents' $\Sigma\psi\circ M$ and
$$h(s,1,x,y)=-h(-s,0,y,x)\quad\text{for}\quad
(s,1,x,y)\in I\times1\times S^{2q}.$$
Then $h|_{I\times1\times S^{2q}}
$ \lq represents'
$((-1)^{m+q+1}\Sigma\psi)\circ M$. Hence
$$h=\Sigma((1-(-1)^{m-q})\psi)\circ M,\quad\text{so}
\quad
\mu((1-(-1)^{m-q})\pi^S_{2q-m+2})=0.$$
It is easy to see that the above
construction describes the entire $\ker\mu$.

\proclaim{Non-injectivity Example \kt7} If $l\ne3,7$ is odd and
$\Sigma^3:\pi_{2l-1}(S^{l-1})\to\pi^S_l$ is epimorphic, then
$\alpha^{3l}(S^1\times S^{2l-1})$ is not injective.
$$\minCDarrowwidth{12pt}\CD
\pi_{2l-1}(V_{l,1})  @>> \mu'' > \pi_{2l-1}(V_{l+1,2}) @>> \nu'' >
\pi_{2l-1}(S^l) \\
@VV = V @VV \tau V @VV = V \\
\pi_{2l-1}(S^{l-1}) @>> \tau\mu'' > \Emb^{3l}_{PL}(S^1\times S^{2l-1})
@>> \nu' > \pi_{2l-1}(S^l)\\
@VV \Sigma^3 V @VV \alpha' V @VV \Sigma V \\
\pi^S_l @>> \mu > \Pi^{3l-1}_{1,2l-1} @>> \nu > \pi^S_{l+1} \endCD.$$
\endproclaim

\demo{Proof}
The case $l=1$ is obvious, so suppose that $l>1$.
By Decomposition Lemma \kt3 there exists (anti)commutative diagram as above.
Since $3l\ge2+(2l-1)+2$, it suffices to
construct embeddings $f,g:S^1\times S^{2l-1}\to\R^{3l}$ such that
$\nu'(f)\ne\nu'(g)$ but $\alpha' g=\alpha' f$.
Let $\varphi=[\iota_l,\iota_l]\in\pi_{2l-1}(S^l)$.
Recall that $\varphi\ne0$ for $l\ne1,3,7$ but $\Sigma\varphi=0$.
Since $l$ is odd, it follows that there is a section
$s:S^l\to V_{l+1,2}$ such that $\nu''s_*=\id$.
Let $f=\tau s_*\varphi$.  We have
$\alpha'f\in\ker\nu=\im\mu$, hence there is $y\in\pi^S_l$ such that
$\mu y=\alpha' f$.  Since $\Sigma^3$ is epimorphic, it follows that there is
$y'\in\pi_{2l-1}(S^{l-1})$ such that $\Sigma^\infty y'=y$.  Let
$g=\tau\mu''y'$.
Now the example follows from
$$\alpha' f=\mu\Sigma^3y'=\alpha'\tau\mu''y'=\alpha' g\quad\text{and}
\quad\nu'f=\varphi\ne0=\nu'\tau\mu''y'=\nu'g.\qed$$
\enddemo

$\Sigma^3:\pi_{2n-1}(S^{n-1})\to\pi_{2n+2}(S^{n+2})=\pi^S_n$ is
epimorphic for each integer $n\not\in\{1,2,3,7\}$

\endcomment

%\newpage
\head \ex\ The Borromean rings and the Haefliger-Wu invariant \endhead
%\subhead Borromean rings \endsubhead
All examples illustrating that the metastable dimension restrictions in
embedding theorems are sharp have their origin in the Borromean Rings Example
\lc2.
So let us illustrate the idea of the Non-injectivity Example \wu7.d and the
Non-embeddability Example \wu9.c by an alternative (comparatively to Example
\lc2) construction of {\it three circles embedded into $\R^3$ so that every
pair of them is unlinked but all three are linked together.}
%We shall call such triples of circles in $\R^3$ {\it generalized
%Borromean rings}.
Our construction is based on the fact that the fundamental group is not
always commutative.

Take two unknotted circles $\Sigma$ and
$\bar\Sigma$ in $\R^3$ far away from each other.
Embed in $\R^3-(\Sigma\sqcup\bar\Sigma)$ the Figure Eight,
i.e.\ the wedge $C$ of two oriented circles so that the
inclusion $C\subset\R^3-(\Sigma\sqcup\bar\Sigma)$ induces an
isomorphism of the fundamental groups.  Take generators $a$ and
$\bar a$ of $\pi_1(C)\cong\pi_1(\R^3-(\Sigma\sqcup\bar\Sigma))$
represented by the two arbitrarily oriented circles of the
Figure Eight.  Consider a map $S^1\to C\subset\R^3$ representing
the element $a\bar aa^{-1}\bar a^{-1}$.  By general position,
there is an embedding $f:S^1\to\R^3$, very close to this map.
It is easy to choose $f$ so that $\Sigma$ and $f(S^1)$,
$\bar\Sigma$ and $f(S^1)$ are unlinked (because $f$ induces the
zero homomorphism of the 1-dimensional homology groups).  Then
$\Sigma$, $\bar\Sigma$ and $f(S^1)$ are as required (cf. Figure \ex1).
 Indeed, $\Sigma$ and $\bar\Sigma$ are unlinked by their
definition.  But $f$ induces a nonzero homomorphism of the
fundamental groups.  Therefore the three circles $\Sigma$,
$\bar\Sigma$ and $f(S^1)$ are linked together.

Higher-dimensional Boromean rings can be constructed analogously using
Whitehead products instead of commutators.

\demo{Sketch of a counterexample to the relative versions of Theorem \vk1 for
$n=2$ and to the surjectivity in the Haefliger-Weber Theorem \wu4 for $m=2n=4$}
Let
$$N=D^2\sqcup D^2\sqcup D^2\quad\text{and}
\quad A=\partial D^2\sqcup\partial D^2\sqcup\partial D^2.$$
Let $A\subset S^3\cong\partial D^4$ be (generalized) Borromean rings.

Since all the three rings are linked, it follows that the
embedding $A\to\partial D^4$ cannot be extended to an embedding
$N\to D^4$. But since each pair of Borromean rings is unlinked, it
follows that the corresponding relative Haefliger-Wu or van Kampen
obstruction to this extension vanishes. This is so because the
Haefliger-Wu obstruction or the van Kampen obstruction involves
2-fold products and double intersections but does not involve
3-fold products and triple intersections. \qed\enddemo

\bigskip
\centerline{\epsffile{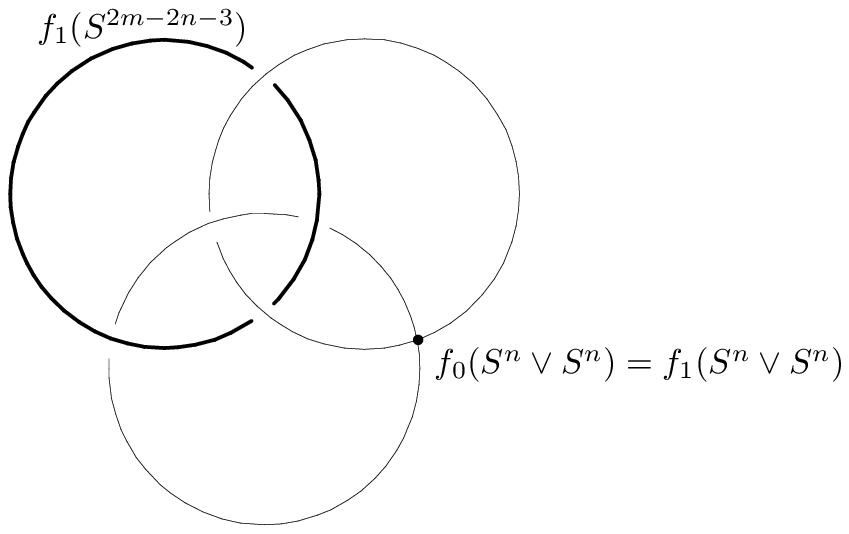}}
\centerline{\it Figure \ex1} %approximately here}
\bigskip

\demo{Proof of the Non-injectivity Example \wu7.d} The reader is
recommended to read this proof first for $n=1$ and $m=3$. Let
$N=S^n\vee S^n\sqcup S^{2m-2n-3}$. Take the standard embedding
$f_0:N\to S^m$. Then $S^m-f_0(S^n\vee S^n)\simeq S^{m-n-1}\vee
S^{m-n-1}$. Take a map $\varphi:S^{2m-2n-3}\to S^m-f_0(S^n\vee
S^n)$ representing the Whitehead product (for $m-n=2$, a
commutator) of generators. If $n=1$ and $m=3$, then $\varphi$ is
homotopic to an embedding by general position. If $n>1$, then
$2m\le3n+3$ implies that $m\le2n$, i.e. $m-(2m-2n-3)\ge3$. Since
also $2(2m-2n-3)-m+1\le m-n-2$ by the
Penrose-Whitehead-Zeeman-Irwin Embedding Theorem \wi6.c, it
follows that $\varphi$ is homotopic to an embedding. Define
$f_1:N\to S^m$ on $S^n\vee S^n$ as $f_0$ and on $S^{2m-2n-3}$ as
such an embedding (Figure \ex1). Since the homotopy class of
$\varphi$ is non-trivial, it follows that $f_1:N\to S^m$ is not
isotopic to the standard embedding $f_0:N\to S^m$.

Let us prove that $\alpha(f_0)=\alpha(f_1)$.
Using \lq finger moves' (analogously to the construction of
Example \ex1 below) 
%; for $n=1$ and $m=3$ see Figure \ex2) 
we construct a map $F:N\times I\to\R^m\times I$ such that
$$F(x,0)=(f_0(x),0),\quad F(x,1)=(f_1(x),1)\quad\text{and}
\quad F((S^n\vee S^n)\times I)\cap F(S^{2m-2n-3}\times I)=\emptyset.$$

%\bigskip
%\centerline{\epsffile{7-2'.ps}}
%\centerline{\it Figure \ex2} %approximately here}
%\bigskip

Then there is a triangulation $T$ of $N$ such that no images of disjoint
simplices intersect throughout $F_t$.
Then the map $\t{F_t}$ is well-defined on the simplicial deleted product $\t T$.
Since $\t T$ is an equivariant deformation retract of $\t N$, it follows that
$\alpha(f_0)=\alpha(f_1)$.

(One can also check that in general $\alpha^m_G(N)f_0=\alpha^m_G(N)f_1$ for
each $G$, so $\alpha^m_G(N)$ is not injective.)
\qed\enddemo

\bigskip
\centerline{\epsffile{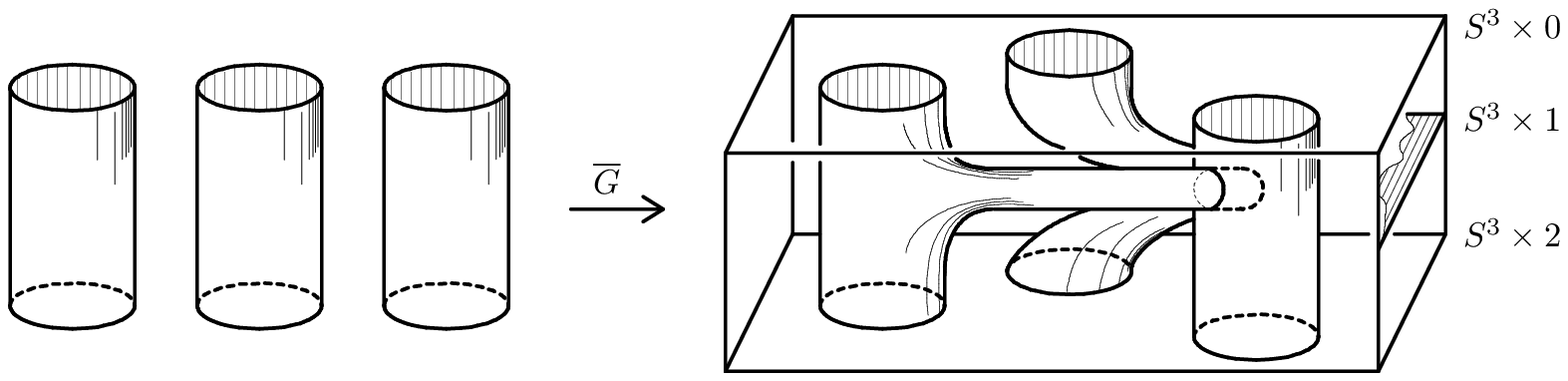}}
\centerline{\it Figure \ex3} %approximately here}
\bigskip

The above gives us {\it a 3-dimensional visualization of the celebrated
Casson finger moves}.  Combine the homotopy $F$ constructed in the above
proof for $n=1$ and $m=3$ with the \lq reverse' homotopy.  We get a homotopy
$G:N\times I\to\R^3\times I$ between standard embeddings $N\to\R^3$ (Figure
\ex3).  This homotopy is obtained from the identity isotopy by Casson finger
moves.

We conjecture that the non-trivial embedding $f$ of the Non-injectivity Example \wu7.b can be
obtained from explicitly defined Borromean rings
$S^n\sqcup S^n\sqcup S^{2m-2n-3}\subset\R^m$ [Ha62', Ma90, Proposition 8.3] by
\lq wedging' $S^n\sqcup S^n$.
We also conjecture that by joining the two $n$-spheres of the above linking
by a tube we obtain a non-trivial embedding $S^n\sqcup S^{2m-2n-3}\to\R^m$
with trivial $\alpha$-invariant (although this is harder to prove: either we
assume that $m-n\not\in\{2,4,8\}$ and need to check that the linking
coefficient of the obtained link is $[\iota_{m-n-1},\iota_{m-n-1}]\ne0$, or we
need to apply the $\beta$-invariant [cf.\ Ha62', \S3, Sk]).

%\subhead The Freedman-Krushkal-Teichner example and its modification
%\endsubhead

\smallskip
{\bf Example \ex1.} {\it There exist a 2-polyhedron $N$ non-embeddable into
$\R^4$ but for which there exist a PL non-degenerate almost-embedding
$f:N\to\R^4$.}

\smallskip
The definition of a non-degenerate almost embedding is given before Lemma
\vk2.
The polyhedron $N$ from Example \ex1 is even topologically non-embeddable into
$\R^4$.

Before constructing Example \ex1 let us explain its meaning.
By the definition of the van Kampen obstruction (\Vk), Example \ex1 implies the
Freedman-Krushkal-Teichner example, i.e. Theorem \vk1 for $n=2$.
Take a triangulation $T$ of $N$ from the definition of an almost embedding
$f:N\to\R^4$.
Then the Gauss map $\t f:\t T\to S^3$ is well-defined on $\t T$.
Since $\t T$ is an equivariant deformation retract of $\t N$, it follows that
there exists an equivariant mapping $\t N\to S^3$.
So Example \ex1 implies the Non-Embeddability Example \wu9.c for $m=2n=4$.
But Example \ex1 gives even more and shows the non-surjectivity not only of the
Haefliger-Wu invariant, but also of the generalized Haefliger-Wu invariants
(\Wu).

\demo{Preliminary construction for Example \ex1}
Let $Q$ be the 2-skeleton of the 6-simplex minus the
interior of one
2-simplex from this 2-skeleton.
Recall from the subsection \lq Ramsay link theory' of \Vk\ that

{\it (*) $Q$ contains two disjoint spheres $\Sigma^2$ and $\Sigma^1$ such
that for each embedding $Q\to\R^4$ these spheres link with an odd linking
number.}

An alternative proof of this fact is presented in the proof
of the Linking
Lemma \ex2 below.

In this section $\bar Q$ denotes a copy of the space $Q$ (for a subset
$A\subset Q$ its copy is denoted by $\bar A\subset\bar Q$).
Embed $Q\sqcup\bar Q$ into $\R^4$ in the standard way, i.e.\ so that

(a) the copies $Q$ and $\bar Q$ are far away from one another;

(b) both $\Sigma^2$ and $\bar\Sigma^2$ are unknotted.

Then $\Sigma^2$ and $\bar\Sigma^2$ are unlinked.
Take any point $x\in\Sigma^1$.
Join the points $x$ and $\bar x$ by an arc in $\R^4$ and
pull small neighborhoods in $Q$ and $\bar Q$ of these points to each other
along this arc.
We obtain an embedding $Q\vee\bar Q\subset\R^4$ (Figure \ex4).

\bigskip
\centerline{\epsffile{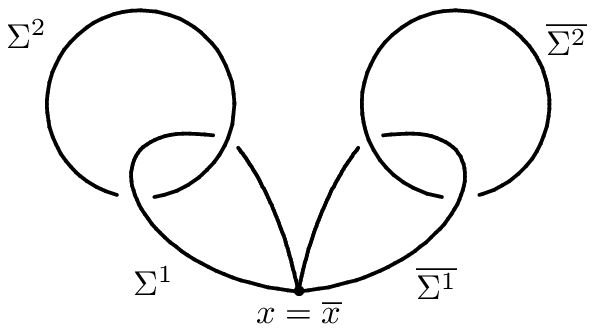}}
\centerline{\it Figure \ex4} %approximately here}
\bigskip

In Figure \ex4 each sphere $\Sigma^1$, $\bar\Sigma^1$, $\Sigma^2$ and
$\bar\Sigma^2$ of dimensions 1, 1, 2, and 2, respectively, is shown as
1-dimensional.
Consider the wedge $\Sigma^1\vee\bar\Sigma^1$
%be a Figure Eight
with the base point $x=\bar x$.
Then the inclusion
$\Sigma^1\vee\bar\Sigma^1\subset\R^4-(\Sigma^2\sqcup\bar\Sigma^2)$ induces an
isomorphism of the fundamental groups.
Take generators $a$ and $\bar a$ of the group $\pi_1(\Sigma^1\vee\bar\Sigma^1)$
represented by the
two (arbitrarily oriented) circles of the Figure Eight.
\enddemo

\demo{Sketch of the Freedman-Krushkal-Teichner construction}
Here we sketch a construction of an example for Theorem \vk1 for $n=2$.
This construction is a bit simpler than that of Example \ex1, but it makes
vanishing of the obstruction less clear, and it gives
neither Example \ex1 nor the Non-Embeddability Example \wu9.c.

Take a map $S^1\to \Sigma^1\vee\bar\Sigma^1$ representing the
element $[a,\bar a]=a\bar aa^{-1}\bar a^{-1}$. Let $N'$ be the
mapping cone of the composition of this map with the inclusion
$\Sigma^1\vee\bar\Sigma^1\subset Q\vee\bar Q$, i.e.
$$N'=B^2\bigcup\limits_{[a,\bar a]:\partial B^2\to \Sigma^1\vee\bar\Sigma^1}
(Q\vee\bar Q).$$

Let us sketch the proof of the nonembeddability of $N'$ into $\R^4$.
Suppose to the contrary that there exists an embedding $h:N'\to\R^4$.
The non-trivial element $[a,\bar a]$ of $\pi_1(\Sigma^1\vee\bar\Sigma^1)$ goes
under $h$ to a loop in $\R^4-h(\Sigma^2\sqcup\bar\Sigma^2)$, which extends to
$hD^2$ and hence is null-homotopic.
This is a contradiction because

{\it $h$ induces a monomorphism
$\pi_1(\Sigma^1\vee\bar\Sigma^1)\to\pi_1(\R^4-h(\Sigma^2\sqcup\bar\Sigma^2))$.}

If both $h\Sigma^2$ and $h\bar\Sigma^2$ are unknotted in $\R^4$, then
the above assertion is clear.
In general (i.e.\ when the spheres are knotted) the above assertion is proved
using the Stallings theorem on the lower central series of groups
[St65, FKT94], cf. below.

We have $V(N')=0$ because the van Kampen obstruction can only detects the
{\it homology} property that the loop $[a,\bar a]$ is null-homologous (for a
detailed proof see [FKT94]).

Since the van Kampen obstruction is a complete obstruction for the
existence of an equivariant map $\t{N'}\to S^3$ (for these dimensions),
we obtain Non-embeddability Example \wu9.c for $m=2n=4$.
\qed\enddemo

\demo{Construction of Example \ex1} Take an embedding
$Q\sqcup\bar Q\subset\R^4$ with the properties (a) and (b) from the
Preliminary Construction.
Take any points $x\in\Sigma^1$ and $y\in\Sigma^2$.
Join points $x$ to $\bar x$ and $y$ to $\bar y$ by two arcs in $\R^4$.
Pull small neighborhoods in $Q$
and $\bar Q$ of these points to each other along this arc (Figure \ex5).
We obtain an embedding
$K:=Q\bigcup\limits_{x=\bar x,y=\bar y}\bar Q\subset\R^4$.

\bigskip
\centerline{\epsffile{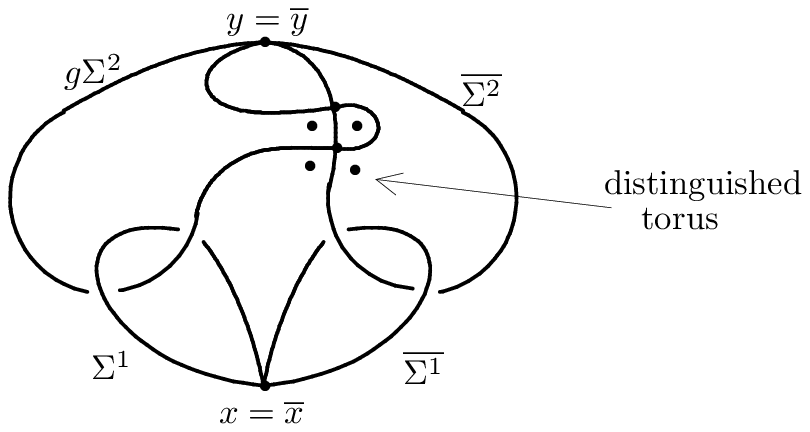}}
\centerline{\it Figure \ex5} %approximately here}
\bigskip

Push an 2-dimensional finger from a small disk $D^2\subset\Sigma^2$ near
$y=\bar y$ intersecting $\bar\Sigma^2$ near $y=\bar y$.
We get a new PL map $g:K\to\R^4$ which has transversal self-intersection
points (Figure \ex5).

We can represent a disk neighborhood $B^4$ of an arbitrary intersection
point $c\in\R^4$ as the product $B^2\times B^2$ of balls, where $0\times 0$
corresponds to the intersection while $B^2\times 0$ and $0\times B^2$
correspond to the images of $\Sigma^2$ and $\bar\Sigma^2$ (we denote by 0
the center of $B^2$).
In a neighborhood of the point $c$ we have the {\it distinguished} or {\it
characteristic} torus $\partial B^2\times\partial B^2$ [cf.\ Ca86, Ki89, FQ90].
In Figure \ex5 the 2-dimensional distinguished torus is shown as
0-dimensional.  By (b) we have
$\pi_1(\R^4-\Sigma^2\vee\bar\Sigma^2)\cong\pi_1( S^1\vee S^1)$.  Denote by
$a$ and $\bar a$ the elements of this group represented by homeomorphisms
$S^1\to z\vee S^1$ and $S^1\to S^1\vee z$ (for some point $z\in S^1$),
respectively (with some orientations).
With appropriate orientations the inclusions of
$\partial B^2\times z$ and $z\times\partial B^2$
into $\R^4-\Sigma^2\vee\bar\Sigma^2$ are homotopic to $a$ and
$\bar a$, respectively.
Since the map
$$a\bar a a^{-1}\bar a^{-1}:S^1\to S^1\vee S^1\cong
(z\times\partial B^2)\vee(\partial B^2\times z)$$
extends to a map $B^2\to \partial B^2\times\partial B^2$,
it follows that $a\bar a a^{-1}\bar a^{-1}$ is null-homotopic
in $\R^4-g(\Sigma^2\vee\bar\Sigma^2)$.
Then there exists a PL map $r:B^2\to\R^4-g(\Sigma^2\vee\bar\Sigma^2)$  such
that $r|_{\partial B^2}:\partial B^2\to\Sigma^1\vee\bar\Sigma^1$
represents the commutator of the inclusions
$\Sigma^1\subset\Sigma^1\vee\bar\Sigma^1$ and
$\bar\Sigma^1\subset\Sigma^1\vee\bar\Sigma^1$.
Roughly speaking, $r(B^2)$ is a torus $0\times\partial B^2\times\partial B^2$.
Set
$$N=B^2\bigcup\limits_{\partial B^2=\partial D^2}(K-\delet D^2)\cup r(B^2).$$
Analogously to the proof of the Freedman-Krushkal-Teichner example above, $N$
does not embed into $\R^4$ (the details are analogous to Construction of
Non-Embeddability Example \wu9.c below).

We have $N\supset(K-\delet D^2)\cup B^2\cong K$.  Define a map
$f:N\to\R^4$ on $(K-\delet D^2)\cup B^2$ as the composition of a
homeomorphism with $K$ and $g$, and on $r(B^2)$ as the identity.  Then
$\Sigma(f)\subset B^2\cup\bar D^2$.  By the construction of the balls $D^2$
and $\bar D^2$ it follows that the balls $B^2$ and $\bar D^2$ are contained
in the interiors of some adjacent 2-simplices of some triangulation $T$ of $N$.
Hence $f$ is a non-degenerate almost embedding (whose image is
$g(K)\cup r(B^2)$).
\qed\enddemo

\demo{Construction of the Non-Embeddability Example \wu9.c}
The case $m=3$ is proved in [GS06] using different ideas.
The case $4\le m\le n+1$ can either be proved analogously to
[GS06] or is covered by the case $4\le m\ge n+2$ of the
Non-Embeddability Example \wu9.c.
So we present the proof for $m\ge n+2$.
This is a higher-dimensional generalization of Example \ex1.

Let $l=m-n-1\ge1$.
Denote by $\Delta^k_{a_0\dots a_s}$ the $k$-skeleton of the
$s$-simplex with vertices $a_0\dots a_s$.

(The definition of $\Delta^n_{012\dots m+2}$ makes sense even
for $n=l+1$, which case is outside the dimension range of the
Non-Embeddability Example \wu9.c. If $n=l+1=1$, then
$\Delta^n_{012\dots m+2}$ is one of the Kuratowski non-planar
graphs, namely $K_5$, and if $n=l+1>1$ then $\Delta^n_{012\dots
m+2}$ is an $n$-dimensional polyhedron non-embeddable in
$\R^{2n}$.)

Set
$$Q=\Delta^n_{12\dots m+2}\cup\Con(\Delta^l_{12\dots m+2}-
\Int\Delta^l_{12\dots l+1},0),
\quad K=Q\bigcup\limits_{0=\bar 0, m=\bar m}\bar Q,$$
$$\Sigma^l=\partial\Delta^{l+1}_{01\dots l+1}\quad\text{and}
\quad\Sigma^n=\partial\Delta^{n+1}_{l+2\dots m+2}.$$
The polyhedron $Q$ embeds in $\R^m$ (this was actually proved in the
first two paragraphs of the proof of [SS92, Lemma 1.1]).
Embed into $\R^m$ two copies of $Q$ which are far apart.
Since $m\ge n+2$, we can join two points of $\Sigma^n$ and
$\bar\Sigma^n$ by an arc and pull the points of the spheres
together along this arc.
Making the same construction for $\Sigma^l$ and $\bar\Sigma^l$ we obtain an
embedding $K\to\R^m$; so we assume that $K$ is a subset of $\R^m$.
We may assume that the wedge
$\Sigma^n\vee\bar\Sigma^n$ is unknotted in $\R^m$.
(If $m>n+2$, then this holds for {\it any} embedding $K\subset\R^m$ [Li65,
Theorem 8]; if  $m=n+2$, then for our embedding $Q\to\R^m$ the sphere
$\Sigma^n$ is unknotted in $\R^m$, and we can choose an embedding $K\to\R^m$ so
that $\Sigma^n\vee\bar\Sigma^n$ is unknotted in $\R^m$.)

Take a triangulation $T$ of $K$.
Let $D^n\subset\Sigma^n$ and $\bar D^n\subset\bar\Sigma^n$ be PL disks
each in the interiors of those $n$-simplices of $T$ that contain the common
point $m=\bar m$ of $\Sigma^n$ and $\bar\Sigma^n$.
Take points $a\in\delet D^n$, $\bar a\in\delet{\bar D^n}$ and a small
arc $s\subset\R^m$ joining $a$ to $\bar a$.
By general position $s\cap K=\{a,\bar a\}$.
Construct a new embedding $g:D^n\to\R^m$ obtained from the old one by
pushing an $n$-dimensional finger from $D^n$ along the arc $s$.
Let $g|_{K-\delet D^n}$ be the inclusion.
We get a new PL map $g:K\to\R^m$ such that $g|_{K-\delet{\bar D^n}}$ is
an embedding but $g(D^n)\cap g(\bar D^n)\ne\emptyset$ (Figure \ex1.b).

By general position $\dim(g(D^n)\cap\bar D^n)\le2n-m$ and $g(D^n)$ intersects
$\bar D^n$ transversally.
Denote by 0 the center of $B^k$.
We can represent a disk neighborhood $B^m$ of an arbitrary point $c$ of this
intersection as the product
$B^{2n-m}\times B^{l+1}\times B^{l+1}$ of balls, where
$B^{2n-m}\times 0\times 0$ corresponds to the intersection, while
$$B^{2n-m}\times B^{l+1}\times 0\quad\text{and}
\quad B^{2n-m}\times 0\times B^{l+1}$$
correspond to $g(D^n)$ and $\t D^n$, respectively.
In a neighborhood of the point $c$ we have the {\it distinguished}  or {\it
characteristic} torus $0\times\partial B^{l+1}\times\partial B^{l+1}$.

Since $\Sigma^n\vee\bar\Sigma^n$ is unknotted in $\R^m\subset S^m$, it
follows that $\pi_l(\R^m-\Sigma^n\vee\bar\Sigma^n)\cong\pi_l(S^l\vee S^l)$.
Denote by $\alpha$ and $\bar\alpha$ the elements of this group represented by
the inclusions of components of the wedge (with some orientations).
Take a point $y\in\partial B^{l+1}$.
With appropriate orientations the inclusions of
$$0\times\partial B^{l+1}\times y\quad\text{and}
\quad0\times y\times\partial B^{l+1}\quad\text{into}
\quad\R^m-\Sigma^n\vee\bar\Sigma^n$$
are homotopic to $\alpha$ and $\bar\alpha$, respectively.
Since the Whitehead product
$$[\alpha,\bar\alpha]:S^{2l-1}\to S^l\vee S^l\cong
(0\times y\times\partial B^{l+1})\vee(0\times\partial B^{l+1}\times y)$$
extends to a map $B^{2l}\to 0\times\partial B^{l+1}\times\partial B^{l+1}$
[Ca86, Ki89, FQ90], it follows that $[\alpha,\bar\alpha]$ is null-homotopic in
$\R^m-\Sigma^n\vee\bar\Sigma^n$.

%By (\ex1.a) below, both $p$ and $\bar p$ are non-zero.

Denote the linking coefficient by $\link(\cdot,\cdot)$
Let $p=\link(\Sigma^l,\Sigma^n)$ and $\bar p=\link(\bar\Sigma^l,\bar\Sigma^n)$.
The inclusions of $\Sigma^l$ and $\bar\Sigma^l$ into
$\R^m-\Sigma^n\vee\bar\Sigma^n$ represent the elements $p\alpha$ and
$\bar p\bar\alpha$ of the group $\pi_l(\R^m-\Sigma^n\vee\bar\Sigma^n)$,
respectively.
Since
$$[p\alpha,\bar p\bar\alpha]=p\bar p[\alpha,\bar\alpha]=0\in
\pi_{2l-1}(\R^m-\Sigma^n\vee\bar\Sigma^n),$$
it follows that the Whitehead product of the (arbitrarily oriented) inclusions
of $\Sigma^l$ and $\bar\Sigma^l$ into $\R^m-\Sigma^n\vee\bar\Sigma^n$ is
null-homotopic.
Hence there exists a PL map $r:B^{2l}\to\R^m-\Sigma^n\vee\bar\Sigma^n$
whose restriction to $\partial B^{2l}$ represents this Whitehead product.
Set
$$N=B^n\bigcup\limits_{\partial B^n=\partial D^n}(K-\delet D^n)\cup r(B^{2l}).$$
Since $m\le\frac{3n}2+1$, it follows that $2l\le n$ and hence $\dim N=n$.
Analogously to the proof of Example \ex1, there exists an almost embedding
$N\to\R^m$.
So it remains to prove that $N$ does not embed into $\R^m$.
\enddemo

\demo{Proof of the PL non-embeddability of $N$ into $\R^m$}
Suppose to the contrary that there is a PL embedding $h:N\to S^m$.
Let
$$\Sigma_1^n=(\Sigma^n-\delet D^n)
\bigcup\limits_{\partial B^n=\partial D^n}B^n \subset N.$$
The map $h\circ r|_{\partial B^{2l}}$ can be extended to map
$$h\circ r:B^{2l}\to\R^m-h(\Sigma_1^n\vee\bar\Sigma^n).$$
Hence $h\circ r|_{\partial B^{2l}}$ is homotopically trivial in
$\R^m-h(\Sigma_1^n\vee\bar\Sigma^n)$.
Now we shall show the contrary and get a contradiction.
Let
$$q=\link(h\Sigma^l,h\Sigma_1^n)\quad\text{and}
\quad\bar q=\link(h\bar\Sigma^l,h\bar\Sigma^n).$$

In the case $m>n+2$ by [Li65, Theorem 8] (cf.\ the construction of $N$ above)
we have $S^m-h(\Sigma_1^n\vee\bar\Sigma^n)\simeq S^l\vee S^l$.
Denote by $\beta$ and $\bar\beta$ the elements of the group
$\pi_l(\R^m-h(\Sigma_1^n\vee\bar\Sigma^n))$ represented by the
homeomorphisms $S^l\to y\vee S^l$ and $S^l\to S^l\vee y$ (where $y\in S^l$),
respectively (with chosen orientations).
Hence the homotopy class of the map
$$h\circ r|_{\partial B^{2l}}:\partial B^{2l}
\to\R^m-h(\Sigma_1^n\vee\bar\Sigma^n)\quad\text{is}
\quad q\bar q[\beta,\bar\beta]\ \in
\ \pi_{2l-1}(\R^m-h(\Sigma_1^n\vee\bar\Sigma^n)).$$
By the Hilton theorem [Po85, complement to Lectures 5 and 6, pp. 231, 257, or
Hu59, p. 511], the map
$$\varphi:\pi_{2l-1}(S^{2l-1})\to\pi_{2l-1}(S^l\vee S^l)\quad
\text{defined by}\quad\varphi(\gamma)=[\beta,\bar\beta]\circ\gamma$$
is injective (this can also be proved by using the homotopy exact
sequence [Hu59, V.3]).  Hence $[\beta,\bar\beta]$ has infinite order. This
implies that the element $q\bar q[\beta,\bar\beta]$ is non-trivial because
both $q$ and $\bar q$ are nonzero analogously to the property (*) in the
preliminary construction for Example \ex1 (see the Linking Lemma \ex2 below).

In the case $m=n+2$ we have $l=1$.
Consider the compositions
$$\Sigma^1\subset\Sigma^1\vee\bar\Sigma^1\to
\R^{n+2}-h(\Sigma^n\vee\bar\Sigma^n)\quad\text{and}\quad
\bar\Sigma^1\subset\Sigma^1\vee\bar\Sigma^1\to
\R^{n+2}-h(\Sigma^n\vee\bar\Sigma^n).$$
They are {\it homologous} to $q\beta$ and $\bar q\bar\beta$, respectively.
The commutator of the homotopy classes of the above compositions is
non-zero because the inclusion
$\Sigma^1\vee\bar\Sigma^1\subset\R^{n+2}-(\Sigma^n\vee\bar\Sigma^n)$
induces a monomorphism of the fundamental groups.
The latter is proved analogously to [FKT94, proof of Lemmas 7 and 8] using
the Stallings theorem [St65] and that by the Linking Lemma 7.2 below
$$\link(\Sigma^n,\Sigma^1)\equiv\link(\bar\Sigma^n,\bar\Sigma^1)
\equiv1\mod2\quad\text{and}\quad\link(\Sigma^n,\bar\Sigma^1)=
\link(\bar\Sigma^n,\Sigma^1)=0.\quad\qed$$
\enddemo

\smallskip
{\bf The Linking Lemma \ex2.} {\it For any PL %or TOP
embedding $K\subset\R^m$ of the above polyhedron $K$ the pairs
$(\Sigma^n,\bar\Sigma^l)$ and $(\bar\Sigma^n,\Sigma^l)$ are
unlinked, and $\link(\Sigma^n,\Sigma^l)$ is odd}
[cf.\ SS92, Lemma 1.4, FKT94, Lemmas 6, 7 and 8].

\demo{Sketch of the proof} The unlinking part  follows because
$\Sigma^l$ (resp.  $\bar\Sigma^l$) bounds a disk $\Delta^{l+1}_{01\dots
l+1}$ (resp.  $\bar\Delta^{l+1}_{01\dots l+1}$) in $K-\bar\Sigma^n$ (resp.
in $K-\Sigma^n$).

We illustrate the idea of proof of the linking part by proving its
particular case $m=2n=4l=4$ (for which, however, there exists a simpler
proof). Recall the formulation for this case:

{\it Let $Q=\Delta^2_{01\dots6}-\delet\Delta^2_{012}$.
Let $\Sigma^1=\partial\Delta^2_{012}$ and let $\Sigma^2$ be the union of
2-simplices disjoint from $\Delta^2_{012}$.
Then for each embedding $Q\to\R^4$ these spheres link with an odd linking
number.}

First we prove a simpler result: that $\link(f\Sigma^2,f\Sigma^1)\ne0$
for each embedding $f:Q\to\R^4$.
(The higher-dimensional analogue of this simpler result is sufficient for
the proof of the non-embeddability for $m\ge n+3$.)
If $\link(f\Sigma^2,f\Sigma^1)=0$, then $f\Sigma^1$ spans a 2-disk outside
$f\Sigma^2$.  Hence we can construct an almost embedding
$\Delta_{01\dots6}^2\to\R^4$.  Therefore there is an equivariant map
$\t{\Delta_{01\dots6}^2}\to S^3$.
Then by [Ya54], any equivariant map
$\t{\Delta_{01\dots6}^2}\to S^4$ should induce a {\it trivial}
homomorphism in $H_4^{eq}(\cdot,\Z_2)$.  But

{\it there is an embedding $g:\Delta_{01\dots6}^2\to\R^5$ such that
$\t g:\t{\Delta_{01\dots6}^2}\to S^4$ induces a {\it non-trivial}
homomorphism in $H_4^{eq}(\cdot,\Z_2)$.}

Indeed,
$\Delta_{01\dots6}^2\cong\Delta_{01\dots5}^2\cup\Con\Delta_{01\dots5}^1$,
where $\Delta^5_{01\dots5}$ is a regular 5-simplex inscribed into the standard unit
4-sphere in $\R^5$ and the vertex of the cone is 0.  This homeomorphism
defines an embedding $g:\Delta_{01\dots6}^2\to\R^5$.  The union of 4-cells
in the simplicial deleted product $\t{\Delta_{01\dots6}^2}$ is a $\Z_2$-equivariant
4-cycle.  Let $p$ be a vertex of $\Delta_{01\dots5}^5$.  The map $\t g$
maps a small neighborhood in $\t{\Delta_{01\dots6}^2}$ of the point
$(p,0)$ homeomorphically onto a neighborhood in $S^4$ of $p$.  Hence the
$\t g$-image of the above 4-cycle is non-trivial, and so $\t g$ induces a
non-trivial homomorphism in $H_4^{eq}(\cdot,\Z_2)$.

In order to prove the full strength of the particular case $m=2n=4l=4$,
assume to the contrary that $\link(f\Sigma^n,f\Sigma^l)$ is even.
Let $\Delta'$ be a polyhedron obtained from $\Delta_{012}^2$ by removing
the interiors of an even number $2r$ of disjoint 2-disks in
$\delet\Delta_{012}^2$ and running $r$ pairwise disjoint tubes between the
holes thus formed.
Then $f\Sigma^1$ spans $\Delta'$ outside $f\Sigma^2$.
Let $X=Q\bigcup\limits_{\partial\Delta^2_{012}=\partial\Delta'}\Delta'$.
Then there is an equivariant map $\t X\to S^3$.
Therefore by [Ya54] any equivariant map $\t X\to S^4$ should induce a
trivial homomorphism in $H_4^{eq}(\cdot,\Z_2)$.

Let $p:\Delta\to\Delta_{01\dots6}^2$ be a map which is the identity on
$Q$ and such that $p(\delet\Delta')=\delet\Delta^2_{012}$.
Then $\t p:\t X\to\t{\Delta_{01\dots6}^2}$ induces an epimorphism in
$H_4^{eq}(\cdot,\Z_2)$ [SS92, p. 278].
Then $\t f\circ\t p$ induces a non-trivial homomorphism in $H_4^{eq}$,
which is a contradiction.
\qed\enddemo

\demo{Proof of the TOP non-embeddability of $N$ into $\R^m$ (sketch)}
For $m\ge n+3$ by [Br72] TOP non-embeddability follows from PL
non-embeddability.
For $m=n+2$ we do the following.
There are arbitrarily close approximations $h':N\to\R^m$ to the
embedding $h$ by PL almost embeddings (for certain triangulations of $N$).
By general position we may assume that $h'|_{\Sigma^1}$ and
$h'|_{\bar\Sigma^1}$ are PL embeddings.
The rest of the proof is analogous to the PL case (in which
we need to replace $h$ by $h'$), because we only used the Linking Lemma \ex2
but not the fact that $h$ is an embedding.  \qed\enddemo

\subhead Appendix: Borromean rings and the Boy immersion \endsubhead
The following was proved in [Ak96'', cf.\ Fr87, ARS02].
Let $h:\R P^2\to\R^3$ be the Boy immersion.
Fix any orientation  on the sphere $S^2$ and on the double
points circle $\Delta(h)$.  Take a small closed ball
$D^3\subset\R^3$ containing the triple point of $h$.
Let $r:S^2\to\R P^2$ be the standard double covering.
% and $h_1=h\circ r$.
Denote by $\pi_1:\R^3\times\R\to\R^3$ and $\pi_2:\R^3\times\R\to\R$ the
projections.
Take  a general position smooth map
$\bar h:(S^2-h_1^{-1}\delet D^3)\to\R^3\times\R$ such that
$\pi_1\circ\bar h=h\circ r$ and for each points
$x,y\in S^2-(h\circ r)^{-1}D^3$ such that $hrx=hry$ and
$\pi_2\bar hx>\pi_2\bar hy$ the following three vectors form a
{\it positive} basis of $\R^3$ at the point $hrx=hry$: the vector
along the orientation of $\Delta$, the normal vector to the
small sheet of $h(\R P^2)$ containing $x$ and the normal vector
to the small sheet of $h(\R P^2)$ containing $y$.
Then $\bar h|_{(h\circ r)^{-1}\partial D^3}\to\partial D^3\times\R$ forms the
Borromean rings linking (after the identification $\partial
D^3\times\R\cong\R^3$).

%\newpage
\head \web\ The Disjunction method \endhead

In this section we illustrate the disjunction method by sketching some ideas
of the proof of the surjectivity in the Haefliger-Weber Theorem \wu4 for the
PL case, and of Theorem \wu5.
The first of these results can be restated as follows.

\proclaim{The Weber Theorem \web1} [We67] Suppose that $N$ is an
$n$-polyhedron, $2m\ge3n+3$ and $\Phi:\t N\to S^{m-1}$ is an equivariant map.
Then there is a PL embedding
$$f:N\to\R^m\quad\text{such that}\quad\t f\simeq_{eq}\Phi\quad\text{on}
\quad\t N.$$
\endproclaim

Simplices of any triangulation $T$ are assumed to be linearly ordered with
respect to increasing dimension.
The lexicographical ordering on $T\times T$ is used.

\subhead Proof of the Weber Theorem \web1 for $m=2n+1$ \endsubhead
We present the proof for $m=3$ and $n=1$ (the general case $m=2n+1$ is proved
analogously).
Take a general position map $f:N\to\R^3$.
Then it is an embedding.
So we only need to modify it in order to obtain the property
$\t f\simeq_{eq}\Phi$. This property does not follow by general position.
We obtain it by applying {\it the van Kampen finger moves}, i.e.\ by winding
edges of the graph $N$ around images of other edges.
Note that van Kampen invented his finger moves for the proof of Lemma \vk2,
and here we present {\it a generalization} of the van Kampen finger moves.

\bigskip
\centerline{\epsffile{8-1.eps}}
\centerline{\it Figure \web1} %approximately here}
\bigskip

\proclaim{Proposition \web2} Let $T$ be a triangulation of a 1-polyhedron $N$
(i.e. $T$ is a graph representing this 1-polyhedron).
For each pair of edges $\sigma,\tau\in T$ such that $\sigma\le\tau$ there exists a PL
embedding
$$f:N\to\R^3\quad\text{such that}\quad\t f\simeq_{eq}\Phi\quad\text{on}
\quad J_{\sigma\tau}=
\bigcup\limits_{\sigma\times\tau>\alpha\times\beta\in\t T}\alpha\times\beta.$$
\endproclaim

Theorem \web1 for $m=2n+1=3$ follows from Proposition \web2 by taking $\sigma$
and $\tau$ to be the last simplices of $T$.

\demo{Proof of Proposition \web2. 1st step: construction of balls}
By induction on $\sigma\times\tau$.
If both $\sigma$ and $\tau$ are the first edge of $T$, then
$\dim J_{\sigma\tau}=1$,
hence Proposition \web2 is true by general position.
Now suppose as inductive hypothesis that
$\t f\simeq_{eq}\Phi$ on $J_{\sigma\tau}$ for an embedding $f:N\to\R^3$.
We need to prove that for $\sigma\cap\tau=\emptyset$ there exists a map
$$f^+:N\to\R^3\quad\text{such that}\quad\t f^+\simeq_{eq}\Phi\quad\text{on}
\quad J_{\sigma\tau}\cup\sigma\times\tau\cup\tau\times\sigma.$$
%If $\sigma\cap\tau\ne\emptyset$, then there is nothing to prove, so suppose 

Take points $C_\sigma\in\delet\sigma$ and $C_\tau\in\delet\tau$ (Figure \web2).
Join their images by an arc
$$C^1\subset\R^3\quad
\text{such that}\quad C^1\cap fN=\{fC_\sigma,fC_\tau\}.$$
Let $D^3$ be a small ball neighborhood of $C^1$ in $\R^3$ such that
$f^{-1}D^3$ is disjoint union of arcs $D_\sigma\subset\delet\sigma$ and
$D_\tau\subset\delet\tau$ containing points $C_\sigma$ and $C_\tau$,
respectively.
We may assume that $fD_\sigma$ is unknotted in $D^3$.
Hence a homotopy
equivalence $h:D^3-fD_\sigma\to S^1$ is constructed analogously to the homotopy
equivalence $S^m-fS^q\to S^{m-q-1}$ from the definition of the linking
coefficients at the beginning of \S3.
\enddemo

\bigskip
\centerline{\epsffile{8-2.eps}}
\centerline{\it Figure \web2} %approximately here}
\bigskip

\demo{Proof of Proposition \web2. 2nd step: the van Kampen finger move}
In order to construct such an $f^+$ we shall wind the arc $f|_{D_\tau}$
around $fD_\sigma$ in $D^3-fD_\sigma$ (Figure \web3).

\bigskip
\centerline{\epsffile{8-3.eps}}
\centerline{\it Figure \web3} %approximately here}
\bigskip

Take any embedding
$$f^+:D_\sigma\sqcup D_\tau\to D^3\quad\text{such that}\quad f^+=f\quad\text{on}
\quad D_\sigma\sqcup\partial D_\tau.$$
Let $D_+$ be a copy of $D_\tau$.
Identify $S^1$ with $D_\tau\bigcup\limits_{\partial D_\tau=\partial D_+}D_+$.
Define a map
$$h_{ff^+}:S^1\to D^3-fD_\sigma\overset{h}\to\to S^1\quad\text{by}
\quad h_{ff^+}(x)=\cases h(f(x)),&x\in D_\tau\\ h(f^+(x)),&x\in D_+\endcases.$$
Since $f^+=f$ on $D_\sigma\sqcup\partial D_\tau$, it follows that there is a
homotopy $f_t:D_\sigma\sqcup D_\tau\to D^3$  from $f$ to $f^+$ fixed on
$D_\sigma\sqcup\partial D_\tau$.
Denote by $\t f$, $\Phi$ and $\t{f^+}$ the restrictions of these maps
to $D_\sigma\times D_\tau$, and by $\t f_t$ the restriction of this map to
$\partial(D_\sigma\times D_\tau)$.
Define a map
$$H_{\t f\t f^+}:\partial(D_\sigma\times D_\tau\times I)\to S^2\quad\text{by}$$
$$H_{\t f\t f^+}|_{D_\sigma\times D_\tau\times0}=\t f,
\qquad H_{\t f\t f^+}|_{D_\sigma\times D_\tau\times1}=\t f^+,
\qquad H_{\t f\t f^+}|_{\partial(D_\sigma\times D_\tau)\times I}=\t f_t.$$

\bigskip
\centerline{\epsffile{8-4.eps}}
\centerline{\it Figure \web4} %approximately here}
\bigskip

By [We67, Lemma 1] we have
$$[H_{\t f\t f^+}]=\Sigma[h_{ff^+}]\in\pi_2(S^2).$$
By the equivariant analogue of the Borsuk Homotopy Extension Theorem,
there is an equivariant extension $\Psi:\t K\to S^2$ of
$\t f|_{J_{\sigma\tau}\cup(\sigma\times\tau-\delet D_\sigma\times\delet D_\tau)}$ such that
$\Psi\simeq_{eq}\Phi$.
We may assume that
$\Psi=\Phi$, so that
$\Phi=\bar f$ on $\partial(D_\sigma\times D_\tau)$.
Therefore $H_{\Phi\t f^+}$ can be defined analogously to the above.
Also we can define $H_{\Phi\t f}$ analogously to the above using the constant homotopy
between $\Phi$ and $\t f$ on $\partial(D_\sigma\times D_\tau)$.
Then
$$[H_{\Phi\t f^+}]=[H_{\Phi\t f}]+[H_{\t f\t f^+}]=
[H_{\Phi\t f}]+\Sigma[h_{ff^+}]\in\pi_2(S^2).$$
For {\it every} element $\beta\in\pi_1(S^1)$ there is an embedding
$f^+:D_\tau\to D^3-fD_\sigma$ such that $[h_{ff^+}]=\beta$.
Therefore by the Suspension Theorem, there exists a map
$f^+:D_\tau\to D^3-fD_\sigma$ such that $[H_{\Phi\t f^+}]=0$.
Extend this $f^+$ to all $N$ by $f$.
Then $f^+$ is as required.
\qed\enddemo

\subhead Proof of the Weber Theorem \web1 for $m=2n\ge6$ \endsubhead
Let us introduce some natural and useful definitions.
Fix a triangulation $T$ of an $n$-polyhedron $N$.
A map $f:N\to\R^m$ is an embedding if and only if the following conditions hold:

$f$ is {\it ($T$-)nondegenerate}, i.e.\ $f|_\alpha$ is an embedding for each
$\alpha\in T$;

$f$ is a {\it ($T$-)almost embedding}, i.e.\ $f\alpha\cap f\beta=\emptyset$ for
each $\alpha\times\beta\subset\t T$;

$f$ is a {\it $T$-immersion}, i.e.\ $f\alpha\cap f\beta=f(\alpha\cap\beta)$ for each
$\alpha,\beta\in T$ such that $\alpha\cap\beta\ne\emptyset$.

\demo{Plan of the proof of the Weber Theorem \web1 for $m=2n\ge6$}
Take a triangulation $T$ of $N$ and a general position map $f:N\to\R^m$ that
is linear on simplices of $T$.
Hence $f$ is nondegenerate.
By general position, the properties of almost embedding and
$T$-immersion hold unless $\dim\alpha=\dim\beta=n$.

\bigskip
\centerline{\epsffile{8-5.eps}}
\centerline{\it Figure \web5} %approximately here}
\bigskip

{\it Step 1.} Wind $n$-simplices of $T$ around $(n-1)$-simplices (analogously
to the case $m=2n+1$) and thus modify $f$ to obtain additionally the condition
$f\simeq_{eq}\Phi$ on the $(2n-1)$-skeleton of $\t T$ (Figure \web5).
Such a {\it van Kampen finger move} is a higher-dimensional generalization of
the fifth Reidemeister move (Figure \vk2.V).

{\it Step 2.} Now it is possible to remove intersections of disjoint
$n$-simplices of $T$ and thus modify $f$ so as to make it additionally an
almost embedding.
This construction is a higher-dimensional generalization of the first and the
second Reidemeister moves (Figure \vk2.I and \vk2.II), which are called {\it
the Penrose-Whitehead-Zeeman trick} and {\it the Whitney trick} (cf.\ \Vk),
respectively.
The details for Step 2 are given as Proposition \web3 below.

Step 1 and Step 2 together are analogues of Lemma \vk2.

{\it Step 3.} Wind $n$-simplices of $T$ around $n$-simplices (analogously to
the case $m=2n+1$) and thus modify
$f$ to obtain additionally the condition $f\simeq_{eq}\Phi$ on $\t T$.
These are the van Kampen finger moves in other dimension.

{\it Step 4.} Remove unnecessary intersections of $n$-simplices of $T$ having
a common face and thus modify $f$ so as to make it additionally a
$T$-immersion (and hence an embedding).
This construction is a higher-dimensional generalization of the fourth
Reidemeister move (Figure \vk2.IV), an analogue of the
Freedman-Krushkal-Teichner Lemma \vk3, and is called {\it the
Freedman-Krushkal-Teichner trick}.
\qed\enddemo

\proclaim{Proposition \web3} (cf.\ Proposition \web2)
Let $T$ be a triangulation of an $n$-polyhedron $N$.
For each pair of $n$-simplices $\sigma,\tau\in T$ such that $\sigma<\tau$
there exists a nondegenerate PL map
$$f:N\to\R^{2n}\quad\text{such that}\quad f\alpha\cap f\beta=\emptyset\quad
\text{if}\quad\sigma\times\tau>\alpha\times\beta\in\t T.$$
\endproclaim

\bigskip
\centerline{\epsffile{8-7.eps}}
\centerline{\it Figure \web6} %approximately here}
\bigskip

\demo{Proof} As inductive hypothesis, assume that we have
such $f$ for a pair $\sigma\times\tau\subset\t T$.
We need to prove that there exists a map
$$f^+:N\to\R^{2n}\quad\text{such that}\quad f\alpha\cap f\beta=\emptyset
\quad\text{if}\quad\sigma\times\tau\ge\alpha\times\beta\in\t T.$$
By the inductive hypothesis $f\sigma\cap f\tau=f\delet\sigma\cap f\delet\tau$
is a finite set of points.
Let $C_\sigma\subset\delet\sigma$ be an arc containing the points of
$\sigma\cap f^{-1}\tau$ (Figure \web6).
Let $C_\tau\subset\delet\tau$ be an arc containing the points of
$\tau\cap f^{-1}\sigma$ \lq in the same order' as $C_\sigma$.
Let $C^2\subset\R^{2n}$ be a union of disks such that
$C^2\cap f\sigma=fC_\sigma$ and $C^2\cap f\tau=fC_\tau$.

Let $D^{2n}$ be a small neighborhood of $C^2$ in $\R^{2n}$ such that
$f^{-1}D^{2n}$ is disjoint union of $n$-balls $D_\sigma\subset\delet\sigma$ and
$D_\tau\subset\delet\tau$, which are small neighborhoods of the arcs $C_\sigma$
and $C_\tau$ in $N$.
We have

$f|_{D_\sigma}$ and $f|_{D_\tau}$ are proper embeddings into $D^{2n}$;

$f\sigma\cap f\tau\subset\delet D^{2n}$;

$D_\sigma=\sigma\cap f^{-1}D^{2n}$ and $D_\tau=\tau\cap f^{-1}D^{2n}$; and

$f\partial D_\sigma\cap fD_\tau=f\partial D_\tau\cap fD_\sigma=\emptyset$.

Since $n\ge3$, it follows that a homotopy equivalence
$h:D^{2n}-fD_\sigma\to S^{n-1}$ is constructed analogously to the homotopy
equivalence $S^m-fS^q\to S^{m-q-1}$ from the definition of the linking
coefficients at the beginning of \S3.
The {\it coefficient of the intersection} of $fD_\sigma$ and $fD_\tau$ is
the homotopy class
$$I(fD_\sigma,fD_\tau)\ =\ [f|_{\partial D_\tau}:\partial D_\tau
\to D^{2n}-fD_\sigma\overset{h}\to\to S^{n-1}]\ \in\ \pi_{n-1}(S^{n-1}).$$
We have
$$\pm\Sigma^n I(fD_\sigma,fD_\tau)=[\t f|_{\partial(D_\sigma\times D_\tau)}]
=[\Phi|_{\partial(D_\sigma\times D_\tau)}]=0\in\pi_{2n-1}(S^{2n-1}).$$
Here the first equality is [We67, Proposition 1], the second equality holds
because $\t f\simeq\Phi$ on $\partial(D_\sigma\times D_\tau)$ by the inductive
hypothesis, and the third equality holds since $\Phi$ is defined over
$\t T\supset D_\sigma\times D_\tau$.
By the Freudenthal Suspension Theorem $I(fD_\sigma,fD_\tau)=0$.
Hence the embedding $f|_{\partial D_\tau}$
extends to a {\it map} $f^+:D_\tau\to D^{2n}-fD_\sigma$.
Using the Penrose-Whitehead-Zeeman trick we modify $f'$ to an {\it
embedding}.
Extend $f^+$ over the entire $N$ by $f$.
By general position $C^2$ (and hence $D^{2n}$ and $f^+D_\tau$) is disjoint
with $f(N-\st\sigma-\st\tau)$.
Therefore $f^+$ is as required.
\qed\enddemo

\subhead Definition of a regular neighborhood \endsubhead
The notions of {\it collapse} and {\it regular neighborhood} are used in our
proofs and are generally important in topology.

A polyhedron $Y$ is said to be obtained from a polyhedron $K$ by an
{\it elementary collapse}, if $K=Y\cup B^n$ and $Y\cap B^n=B^{n-1}$, where
$B^{n-1}$ is a face of the ball $B^n$.
This elementary collapse is said to be made {\it from}
$\Cl(\partial B^n-B^{n-1})$ {\it along} $B^n$ {\it to} $B^{n-1}$.
A polyhedron $K$ {\it collapses to} $Y$ (notation: $K\searrow Y$)
if there exists a sequence of elementary collapses
$K=K_0\searrow K_1\searrow \dots\searrow K_n=Y$.
A polyhedron $K$ is {\it collapsible}, if it collapses to a point.

Clearly, the ball $B^n$ is collapsible, since it is collapsible to its face
$B^{n-1}$ and so on by induction.
Moreover, a cone $cK$ on a compact polyhedron $K$ is collapsible (to its
vertex).
Indeed, note that for each simplex $A\subset K$, the cone $cA$ collapses from
$A$ to $c(\partial A)$, hence $cK$ collapses to a point inductively by
simplices of decreasing dimension.

A collapsing $K\searrow Y$ generates a deformation retraction $r:K\to Y$, given
by deformations of each ball $B^n$ to its face $B^{n-1}$.
Consider a homotopy $H_t$ between the identity map $K\to K$ and the deformation
retraction $r:K\to Y$, given by the collapse $K\searrow Y$.
The {\it trace} of a subpolyhedron $S$ of $K$ under the collapse $K\searrow Y$
is the union of $H_t(S)$ over $t\in[0,1]$ (this in fact depend on the homotopy
not only on the collapse).

Suppose that $K$ is a subpolyhedron of a PL manifold $M$.
A neighborhood $N$ of $K$ in $M$ is called {\it regular}, if $N$ is a compact bounded manifold and $N\searrow K$.
The same polyhedron can have distinct regular neighborhoods.
But the regular neighborhood is unique up to homeomorphism and even
up to isotopy fixed on $K$ [RS72].
The regular neighborhood of a collapsible polyhedron is a ball [RS72].
(The inverse statement, i.e.\ that $K\subset B^n$ and $K\searrow *$
imply $B^n\searrow K$, is true only in codimension $\ge3$ [Hu69].)

Let us fix the following convention.
The notation $R_M(K)$ means \lq a sufficiently small regular neighborhood of
$K$ in $M$', when it first appears, and \lq the regular neighborhood of
$K$ in $M$', after the first appearance.

\subhead Proof of the Weber Theorem for the general case $2m\ge3n+3$
\endsubhead
The proof of the Weber Theorem \web1 for the general case consists of two
steps:

(1) the construction of a nondegenerate almost embedding (an analogue of Lemma
\vk2, the details for this step are given as Proposition \web4 below); and

(2) the construction of an embedding from a non-degenerate almost embedding
(a generalization of Freedman-Krushkal-Teichner Lemma \vk3, the details are
given in [Sk98, RS99]).

\proclaim{Proposition \web4}
Suppose that $N$ is an $n$-polyhedron with a triangulation $T$,
$2m\ge3n+3$ and $\Phi:\t N\to S^{m-1}$ is an equivariant map.
Then for each $\sigma\times\tau\in\t T$ such that $\sigma\le\tau$ there exists
a nondegenerate PL map $f:N\to\R^m$ such that

(*) $f\alpha\cap f\beta=\emptyset$ for each
$\alpha\times\beta<\sigma\times\tau$, and

(**) $\t f\simeq_{eq}\Phi$ on
$J_{\sigma\tau}:=\cup\{\alpha\times\beta\cup\beta\times\alpha\subset\t T\ |
\ \alpha\times\beta<\sigma\times\tau\}$.
\endproclaim

The Weber Theorem \web1 follows from Proposition \web4 by taking
$\sigma$ and $\tau$ to be the last simplex of $T$ and then applying a
generalization of the Freedman-Krushkal-Teichner Lemma \vk3 [Sk98].
The Weber Theorem \web1 can also be proved by first constructing the immersion
and then modifying it to an embedding [Sk02].

\demo{Proof of Proposition \web4}
Take a general position map $f:N\to\R^m$ that is linear on
the simplices of the triangulation $T$.
The map $f$ is already non-degenerate.
By the induction hypothesis on $\sigma\times\tau$ we may assume that $f$ is
non-degenerate and the properties (*) and (**) hold.
Suppose that $p+q\ge m-1$ (otherwise the inductive step
holds by general position).

The first part of the proof (a generalization of Proposition \web3
and the Whitney trick)  is getting the property $f\sigma\cap
f\tau=\emptyset$. The second part of the proof (a generalization
of Proposition \web2 and the van Kampen finger moves) is getting
the property $\t f\simeq_{eq}\Phi$ on
$J_{\sigma\tau}\cup\sigma\times\tau\cup\tau\times\sigma$.
\enddemo

\demo{Constructions of balls $D_\sigma$, $D_\tau$ and $D^m$}
The first step in the proof of Proposition \web4
generalizes the construction of the arcs $l_1,l_2$ and the disk
$D$ in Whitney trick, or
the construction of $D_\sigma$, $D_\tau$ and $D^{2n}$ in the case $m=2n$ above.
Let $\Sigma=f\sigma\cap f\tau$.
By (*) we have
$$f\sigma\cap f\partial\tau=f\partial\sigma\cap f\tau=\emptyset.$$
Hence $\Sigma=f\delet\sigma\cap f\delet\tau$.
By general position, $\dim\Sigma\le p+q-m$.
Let $C_\sigma\subset\delet\sigma$ be the trace of the polyhedron
$\sigma\cap f^{-1}\tau$ under some collapse
$\sigma\searrow(\text{a point in }\delet\tau)$.
Define analogously $C_\tau\subset\delet\tau$.
The polyhedra $C_\sigma,C_\tau$ are generalizations of the arcs $l_1,l_2$ from
the Whitney trick, and of $C_\sigma,C_\tau$ above.
They are collapsible,
$$\Sigma\subset fC_\sigma\cap fC_\tau\quad\text{and}\quad
\dim C_\sigma,\dim C_\tau\le p+q-m+1.$$

Consider a collapse from some PL $m$-ball $J^m$ in $\R^m$,
containing $\Sigma$ in its interior, to a point in $\delet J^m$.
Let $C$ be the trail of $C_\sigma\cup C_\tau$ under this collapse.
The polyhedron $C$ is a generalization of the disk $C$ from the Whitney
trick.
It is collapsible, it contains $C_\sigma\cup C_\tau$ and $\dim C\le p+q-m+2$.
Hence by general position, $C\cap f\sigma=C_\sigma$ and $C\cap f\tau=C_\tau$.

Take the regular neighborhoods of polyhedra $C_\sigma,C_\tau$ and $C$ in some
sufficiently fine (agreeing) triangulations of $\sigma,\tau$ and
$\R^m$, respectively.
They are PL balls
$$D^p_\sigma\subset\delet\sigma,\quad D^q_\tau\subset\delet\tau\quad
\text{and}\quad D^m\subset\R^m\quad\text{such that}$$

(a) $f|_{D_\sigma}$ and $f|_{D_\tau}$ are proper embeddings into $D^m$;

(b) $f\sigma\cap f\tau\subset\delet D^m$;

(c) $D_\sigma=\sigma\cap f^{-1}D^m$ and $D_\tau=\tau\cap f^{-1}D^m$;

(d) $D^m\cap fP=\emptyset$, where $P=N-\st\sigma-\st\tau$.

Only the last property needs a proof.
By (*) we have $C_\sigma\cap P=\emptyset$.
By general position,
$$\dim(fP\cap f\tau)\le n+q-m,\quad\text{hence}
\quad\dim(fP\cap f\tau)+\dim C_\tau<q\quad\text{so}
\quad C_\tau\cap P=\emptyset.$$
Therefore $C\cap fP=\emptyset$, which implies (d).
\enddemo

\demo{A generalization of the Whitney trick} Take PL balls $D^m$,
$D_\sigma$ and $D_\tau$ as above. Since $f$ is non-degenerate and
(*) holds, it follows that $f\partial D_\sigma\cap
fD_\tau=f\partial D_\tau\cap fD_\sigma=\emptyset$. Since
$m-p\ge3$, it follows that a homotopy equivalence
$h:D^m-fD_\sigma\to S^{m-p-1}$ is constructed analogously to the
homotopy equivalence $S^m-fS^q\to S^{m-q-1}$ from the definition
of the linking coefficients at the beginning of \S3. The {\it
coefficient of the intersection} of $f|_{\partial D_\sigma}$ and
$f|_{\partial D_\tau}$ is the homotopy class
$$I(f|_{D_\tau},f|_{D_\tau}):=
[f|_{\partial D_\tau}:\partial D_\tau\to D^m-fD_\sigma\overset{h}\to
\to S^{m-p-1}]\in\pi_{q-1}(S^{m-p-1}).$$
We have
$$\Sigma^p I(f|_{D_\sigma},f|_{D_\tau})=(-1)^{m-p}[\t f|_{\partial(D_\sigma\times D_\tau)}]=
[\Phi|_{\partial(D_\sigma\times D_\tau)}]=0.$$
Here the first equality holds by [We67, Proposition 1].
The second equality holds since $\t f\simeq\Phi$ on $\partial(D_\sigma\times D_\tau)$
by the inductive hypothesis.
The third equality holds since $\Phi$ is defined over
$\t T\supset D_\sigma\times D_\tau$.
Since
$$2p+q\le2m-3,\quad\text{we have}\quad q-1\le2(m-p-1)-2.$$
So by the Freudenthal Suspension Theorem the homomorphism $\Sigma^p$ above is a monomorphism.
Hence the embedding $f|_{\partial D_\tau}$ extends to a {\it map}
$f':D_\tau\to D^m-fD_\sigma$.

Since $2q-m+1\leq m-p-2$, by the Penroze-Whitehead-Zeeman-Irwin
Embedding Theorem \wi6.c it follows that $f'$ is homotopic
$\rel\partial D_\tau$ to an {\it embedding} $f^+:D_\tau\to
D^m-fD_\sigma$. Here we again use the inequality $p+2q\leq2m-3$.
Since $m-q\ge3$, by the relative version of Theorem \wi5.a [Ze63,
Corollary 1 to Theorem 9] it follows that there is an ambient
isotopy $h_t:D^m\to D^m\rel\partial D^m$ carrying $f|_{D_\tau}$ to
$f^+$. Extend $f^+$ over $N$ by the formula
$$f^+(x)=\cases h_1(f(x)),&\text{if}\quad f(x)\in D^m\quad\text{and}
\quad x\in\gamma\quad\text{for some}\quad\gamma\supset\sigma^q\\
f(x),&\text{otherwise}\endcases.$$
It is easy to check that $f^+$ is a non-degenerate PL map satisfying  to the properties (*), (**) and such that
$f^+\sigma\cap f^+\tau=\emptyset$.
\enddemo

\demo{A generalization of the van Kampen finger moves}
We begin with the analogous construction of PL balls.
By general position we can take points
$C_\sigma\in\delet\sigma$ and $C_\tau\in\delet\tau$ such that the restrictions of $f$ to some
small neighborhoods of $C_\sigma$ and $C_\tau$ are embeddings.
Since $x,y\le m-2$, we can join points $fC_\sigma$ and $fC_\tau$ by an arc $C\subset\R^m$
such that $C\cap fN=\{fC_\sigma,fC_\tau\}$.
Let $D^m=R_{\R^m}(C)$.
Then $f^{-1}D^m$ is the disjoint union of PL disks
$D_\sigma\subset\delet\sigma$ and $D_\tau\subset\delet\tau$, which are
regular neighborhoods in $N$ of $C_\sigma$ and $C_\tau$,  respectively.

By the Borsuk Homotopy Extension Theorem, there is an extension
$\Psi:\t T\to S^{m-1}$ of the map
$\t f|_{J_{\sigma\tau}\cup(\sigma\times\tau-\delet D_\sigma\times\delet D_\tau)}$ such that
$\Psi\simeq\Phi$. So $\Psi=\t f$ on $\partial(D_\sigma\times D_\tau)$.
We may assume $\Psi=\Phi$.
By [We67, Lemma 1], for each map
$$f':D_\sigma\sqcup D_\tau\to D^m\quad\text{such that}\quad f'=f\quad\text{on}
\quad D_\sigma\sqcup\partial D_\tau\quad\text{and}
\quad f'D_\sigma\cap f'D_\tau=\emptyset$$
and homotopy $f_t\rel D_\sigma\sqcup\partial D_\tau$ from $f$ to $f'$ we have
$$[H_{\Phi\t f_t\t f'}]=[H_{\Phi\t f}]+[H_{\t f\t f_t\t f'}]=
[H_{\Phi\t f}]+(-1)^{m-p}\Sigma^p[h_{ff'}]\in\pi_{p+q}(S^{m-1}).$$
Here $\Phi$, $\t f$ and $\t f^+$ denote the restrictions of
these maps to $D_\sigma\times D_\tau$; $\t f_t$ denotes the
restriction of this map to $\partial(D_\sigma\times D_\tau)$, and the maps $H$
are defined as in the second step of the proof of Proposition \web2.
Since
$$2p+q\le2m-3,\quad\text{we have}\quad q\le2(m-p-1)-1.$$
So by the Freudenthal Suspension Theorem $\Sigma^p$ is an epimorphism.
Since for {\it every} element $\beta\in\pi_q(S^{m-p-1})$ there is a map (not
necessarily an embedding)
$$f':D_\tau\to D^m-fD_\sigma\quad\text{such that}\quad[h_{ff'}]=\beta\quad\text{
and}\quad f'=f\quad\text{on}\quad D_\sigma\sqcup\partial D_\tau,$$
it follows that we can take $f'$ so that $[H_{\Phi\t f_t\t f'}]=0$.

The rest of the proof is the same as in the generalization of the Whitney
trick.
\qed\enddemo

%\newpage
\subhead Generalization of the Weber Theorem \endsubhead
We illustrate some of the ideas of the proof of Theorem \wu5 by proving the
following weaker result for $d\in\{0,1\}$.

\proclaim{Theorem \web5} Suppose that $N$ is a $d$-connected closed PL
$n$-manifold, \ $2m\ge3n+2-d$, \ $m\ge n+3$ \ and $\Phi:\t N\to S^{m-1}$
is an equivariant map. Then there is a PL embedding $f:N\to\R^m$ [Sk97].
\endproclaim

Recall some classical results and their generalizations required to prove
Theorem \web5.

\proclaim{The Engulfing Lemma \web6} Suppose that $N$ is a $(2k+2-n)$-connected closed $n$-manifold and $K\subset N$
is a $k$-polyhedron such that $n-k\ge3$ and the inclusion
$K\subset N$ is null-homotopic.
Then $K$ can be {\it engulfed} in $N$, i.e.\ is contained in an $n$-ball
$B\subset N$ [PWZ61, Ze63].
\endproclaim

\proclaim{Theorem \web7}
Suppose that $N$ is a closed homologically $(3n-2m+2)$-connected PL $n$-manifold,
$m-n\ge3$ and $g:N\to\R^m$ is a map such that $\Sigma(g)$ is contained in some
PL $n$-ball $B\subset N$. Then there is an embedding $f:N\to\R^m$ such that
$f=g$ on $N-\delet B$.
\endproclaim

\demo{Proof}
The theorem is essentially proved in [Hi65, cf.\ Sk97, Theorem 2.1.2].
Let $M=\R^m-\Int R(g(N-\delet B),g\partial B)$.
Since $N$ is homologically $(3n-2m+2)$-connected, we have by Alexander duality
$$H_i(M)\cong H^{m-1-i}(\R^m-M)\cong H^{m-1-i}(N-\delet B)
\cong H_{n-m+1+i}(N)=0$$ for $i\le2n-m+1$. Since $m-n\ge3$, it
follows that $M$ is simply connected. Therefore by the Hurewicz
Isomorphism Theorem we have that $M$ is $(2n-m+1)$-connected.
Hence by the Penrose-Whitehead-Zeeman-Irwin Embedding Theorem
\wi6.c the embedding $g:\partial B\to\partial M$ extends to an
embedding $f:B\to M$. Extending $f$ by $g$ outside $B$ we complete
the proof. \qed\enddemo

\proclaim{Theorem \web8}
Let $N$ be an $n$-polyhedron with triangulation~$T$.
If $m-n\ge3$ and there exists an equivariant map $\Phi:\t N\to S^{m-1}$,
then there exists a general position nondegenerate PL map $f:N\to\R^m$ such that
$$f\sigma\cap f\tau=f(\sigma\cap\tau)\quad\text{if}
\quad p=\dim\sigma\le\dim\tau=q\quad\text{and}\quad p+q+n\le2m-3.$$
\endproclaim

\demo{Sketch of the proof}
Analogous to the proof of the Weber Theorem \web1.
Recall that the Weber Theorem \web1 is proved by induction on
$\sigma\times\tau\in\t T$.
If $\dim\sigma=p$ and $\dim\tau=q$, then we need the
following dimensional
restrictions:

$p+2q\le2m-3$ to apply the Freudenthal Suspension Theorem (twice) and the
Penrose-Whitehead-Zeeman trick;

$p+q+n\le2m-3$ to get the property
$D^m\cap f(N-\st\sigma-\st\tau)=\emptyset$.
\qed\enddemo

\demo{A reduction of Theorem \web5} It suffices to prove that

{\it (***)\quad For some fine triangulation $T$ of $N$ and a map $f:N\to\R^m$
as is given by Theorem \web8, the inclusion $\Sigma(f)\subset N$ is homotopic to
a map to some $d$-dimensional subpolyhedron of $N$.}

Indeed, since $N$ is $d$-connected, (***) implies that the inclusion
$\Sigma(f)\subset N$ is null-homotopic.
Since $N$ is $(2(2n-m)-n+2)=d$-connected, by the Engulfing Lemma \web6 it follows
that $\Sigma(f)$ is contained in some PL $n$-ball in $N$.
Then $N$ embeds into $\R^m$ by Theorem \web7.
\qed\enddemo

Note that if a map $f$ given by Theorem \web8 is a PL
immersion (i.e.\ a local
embedding), then $\Sigma(f)$ does not intersect the $(2m-2n-3)$-skeleton of
$T$.
Hence it retracts to the $(n-1-(2m-2n-3))=d$-skeleton of a triangulation,
dual to $T$, i.e.\ (***) holds.

\demo{Proof of (***) for $d=0$} We may assume that $2m=3n+2$.
Hence
$$f\sigma\cap f\tau=f(\sigma\cap\tau)\quad\text{unless}\quad
\dim\sigma=\dim\tau=n.$$
For $n$-simplices $\alpha$ and $\beta$ let
$$S_{\alpha\beta}=\cases\alpha\cap f^{-1}f\beta&\alpha\cap\beta=\emptyset\\
f^{-1}\Cl[(f\alpha\cap f\beta)-f(\alpha\cap\beta)]&
\alpha\cap\beta\ne\emptyset\endcases.$$

\bigskip
\centerline{\epsffile{8-6.eps}}
\centerline{\it Figure \web7} %approximately here}
\bigskip

$$\text{Then}\quad\Sigma(f)=\bigcup\limits_{\alpha\ne\beta}S_{\alpha\beta}
\quad\text{and}\quad S_{\alpha\beta}\cap S_{\gamma\delta}=\emptyset\quad
\text{when}\quad\alpha\beta\ne\gamma\delta.$$
Here $\alpha\beta$ is the ordered pair $(\alpha,\beta)$
when $\alpha\cap\beta=\emptyset$ and the non-ordered pair
$\{\alpha,\beta\}$ when $\alpha\cap\beta\ne\emptyset$.
Therefore the contractibility of $\alpha$ (and of $\alpha\cup\beta$ for
$\alpha\cap\beta\ne\emptyset$) implies (***).
%that the inclusion $\Sigma(f)\subset N$ is null-homotopic.
%By general position $\dim\Sigma(f)=\frac n2-1$.
%Since $n-\frac n2+1\ge3$ and $N$ is connected, by Engulfing Lemma \web6 below
%$\Sigma(f)$ is contained in some PL $n$-ball in $N$.
%Then we can construct an embedding $N\to\R^m$ by Theorem \web7 below.
\qed\enddemo

\demo{Proof of (***) for $d=1$}
We may assume that $2m=3n+1$.
Define $\alpha\beta$ and $S_{\alpha\beta}$ as in the case $d=0$.
We denote such pairs $\alpha\beta$ by Latin letters $i$,
$j$, $k$, $l$.
First we prove that

{\it (a) $S_i\cap S_j\cap S_k=\emptyset$ for distinct $i,j,k=1,\dots,s$;

(b) For each $i=1,\dots,s$ there is a contractible polyhedron $A_i\subset N$, containing~$S_i$. }

If $S_i\cap S_j\ne\emptyset$, then there is a contractible
polyhedron~$A_{ij}\subset N$, containing~$A_i\cup A_j$.

Indeed, we may assume that a triangulation~$T$ of $N$ is
such that, for each $x\in N$, the star $\st^2x=\st\st x$ is contractible.
By Theorem \web8 $S_{\alpha\beta}\ne\emptyset$ is possible only when
$$\text{either}\quad\dim\alpha=\dim\beta=n\quad\text{or}\quad
\{\dim\alpha,\dim\beta\}=\{n,n-1\}.$$
By general position, $f$ has no triple points. Therefore each non-empty
intersection of any three of $S_1,\dots,S_s$ can be only of the form
$$S_{\alpha_1\beta}\cap S_{\alpha_2\beta}\cap S_{\alpha_3\beta}=
S_{\alpha\beta}\quad
\text{(or }S_{\beta\alpha _1}\cap S_{\alpha_2\beta}\cap S_{\beta\alpha_3}=S_{\beta\alpha})$$
$$\text{for some}\quad
\alpha_1^n,\alpha_2^n,\alpha_3^n,\beta^n,\alpha\in T,
\quad\alpha=\alpha_1\cap\alpha_2\cap\alpha_3.$$
Since
$S_{\alpha\beta}\ne\emptyset$, it follows that $\dim\alpha=n-1$.
Since $N$ is a closed manifold, then no three distinct $n$-simplices of
$T$ intersect by an $(n-1)$-simplex of $T$. This contradiction
shows that (a) is true.

Let
$$A_{\alpha\beta}=\cases\alpha&\alpha\cap\beta=\emptyset,\\
\alpha\cup\beta&\alpha\cap\beta\ne\emptyset.\endcases$$
If $S_i\cap S_j\ne\emptyset$, then take a point
$a_{ij}\in S_i\cap S_j$ and let $A_{ij}=\st^2a_{ij}$.
From the definition of $S_i$ and $A_i$ it follows that $S\subset A_i$ and
$A_i\cup A_j \subset A_{ij}$.
By the choice of~$T$, $A_i$ and $A_{ij}$ are contractible.

Now we construct a homotopy of $\Sigma(f)$ onto its `reduced' nerve.
From (a) it follows that the sets $S_i\cap S_j$ are
disjoint for distinct non-ordered pairs $i,j=1,\dots,s$.
Take disjoint regular neighborhoods $U_{ij}$ of $S_i\cap S_j$ in $\bigcup\limits_{i=1}^sS_i$.
Since $A_i$ is contractible, it follows that there is a homotopy
$F_i:\Cl\Big(S_i-\bigcup\limits_{j\ne i}U_{ij}\Big)\times I\to A_i$ between the inclusion and a constant map to some point~$a_i\in A_i$.

Suppose that $S_i\cap S_j\ne \emptyset$.
Since $A_{ij}$ is contractible, it follows that there is an arc $l_{ij}\subset A_{ij}$ joining $a_i$ and $a_j$.
Also we can extend homotopies $F_i$ and $F_j$ over~$U_{ij}$ to a homotopy $F_{ij}\:U_{ij}\times I\to N$ between the inclusion and a map of~$U_{ij}$ to~$l_{ij}$ (we can do this first for vertices and then for edges).
Since all~$U_{ij}$ are disjoint and all~$F_{ij}$ are extensions of
$F_i$~and~$F_j$, then all the constructed homotopies define a homotopy
$F:\Big(\bigcup\limits^s_{i=1}S_i\Big)\times I\to N$ between the inclusion
and a map onto the following subgraph of $N$:
$$\Big(\bigcup\limits^s_{i=1}a_i\Big)\cup\Big(
\bigcup\{l_{ij}|1\leq i<j\leq s\text{ and }S_i\cap S_j\ne\emptyset\}\Big).\qed$$
\enddemo

From our proof it follows that Theorem \web5 for $2m\ge3n+1$ is
true even if there only exists an equivariant map to $S^{m-1}$
from the $(\big[\frac{4m}3\big]-2)$-skeleton of $\t T$.
%(equivariant maps are homotopic only on $\big[\frac{4m-5}3\big]$-skeleton)

It would be interesting to know if Theorem \web5 remains true when
$N$ has singularities of dimension at most $m-n-2$.  This is not
clear, contrary to what is written in [Sk97] (because e.g. in the
suspension of a homology sphere the complement to vertices is not
simply-connected so we cannot apply Engulfing Lemma \web6).

\comment

\demo{Proof of (***) for $d=2$}
We may assume that $2m=3n$.
Define $\alpha\beta$ and $S_{\alpha\beta}$ as in the case $d=0$.
By Latin letters $i$, $j$, $k$, $l$ we denote such pairs $\alpha\beta$.
First we prove that

{\it (a) The intersections $S_i\cap S_j\cap S_k$ are either disjoint or the
same for distinct sets $\{i,j,k\}\subset\{1,\dots,s\}$.

(b) For each $i=1,\dots,s$ there is a contractible polyhedron $A_i\subset N$
containing~$S_i$.

If $S_i\cap S_j\ne\emptyset$, then there is a
contractible polyhedron~$A_{ij}\subset N$ containing~$A_i\cup A_j$.

If $S_i\cap S_j\cap S_k\ne\emptyset$, then there is a contractible
polyhedron~$A_{ijk}\subset N$, containing~$A_{ij}\cup A_{jk}\cup A_{ki}$.}

Indeed, we may assume that the triangulation~$T$ of $N$ is such that for each
$x\in N$, the star $\st^4 x$ is contractible.
We shall prove an assertion, equivalent to (a):

{\it if $S_{i_1}\cap\dots\cap S_{i_k}\ne\emptyset$ for some $k\geq4$, then
every three of $S_{i_1}\dots S_{i_k}$ have the same intersection.}

By Theorem \web8 $S_{\alpha\beta}\ne\emptyset$ is possible only when
$$\{\dim\alpha,\dim\beta\}\quad\text{is either}\quad\{n,n\}\quad
\text{or}\quad\{n,n-1\}\quad\text{or}\quad\{n,n-2\}.$$
By general position, $f$ has no quadruple points and each
of its triple points is the intersection of $f$-images of three
$n$-dimensional open simplices of $T$.
Therefore each intersection of $k\geq4$ of $S_1,\dots,S_s$ is of the form
$$S_{\alpha_1\beta}\cap\dots\cap S_{\alpha_k\beta}=
S_{\alpha\beta}\quad\text{(or }
S_{\beta\alpha _1}\cap\dots\cap S_{\beta\alpha_k}=S_{\beta\alpha})$$
$$\text{for some}\quad\alpha_1^n,\dots\alpha_k^n,\beta^n,\alpha\in T,
\quad\alpha=\alpha_1\cap\dots\cap\alpha_k.$$
Consider only the first of the above possibilities, the second is proved
analogously.
Since $S_{\alpha\beta}\ne\emptyset$, it follows that
$\dim\alpha\ge n-2$.
If $\dim\alpha=n-1$, then $k=2$, which is a contradiction.
Therefore $\dim\alpha=n-2$.
The intersection of every three distinct simplices $\alpha_p,\alpha_q,\alpha_r$
from $\{\alpha_1,\dots,\alpha_k\}$ contains $\alpha$ and cannot be
$(n-1)$-dimensional.
Hence $\alpha_p\cap\alpha_q\cap\alpha_r=\alpha$.
Since $f$ has no quadruple points, it follows that
$S_{\alpha _p\beta}\cap S_{\alpha_q\beta}\cap S_{\alpha_r\beta}=
S_{\alpha\beta}$ and we are done.

The polyhedra $A_i$ and $A_{ij}$ are defined as in the case $d=1$.
If $S_i\cap S_j\cap S_k\ne\emptyset$, then take a point
$a_{ijk}\in S_i\cap S_j\cap S_k$ and let $A_{ijk}=\st^4a_{ijk}$.
By the choice of~$T$, $A_i, A_{ij}$ and $A_{ijk}$ are
contractible. Then (b) is proved similarly to the case $d=1$. A
homotopy of $\Sigma(f)$ onto its `reduced' nerve is constructed as
in the case $d=1$. \qed\enddemo

\endcomment

%\newpage


\Refs
\widestnumber\key{ARS0}
\ref \key Ad93 \by M. Adachi
\book Embeddings and Immersions, {\rm Transl. of Math. Monographs \bf 124}
\yr 1993 \bookinfo \publ Amer. Math. Soc. \publaddr
\endref

\ref \key Ak69 \by E. Akin
\paper Manifold phenomena in the theory of polyhedra
\jour Trans. Amer. Math. Soc. \vol 143 \yr 1969 \pages 413--473
\endref

\ref \key Ak96 \by P.~M.~Akhmetiev
\paper On isotopic and discrete realization of mappings from $n$-dimensional
sphere to Euclidean space (in Russian)
\jour Mat.~Sbornik \vol 187:7 \yr 1996 \pages 3--34 \moreref English transl.
\jour Sb. Math.  \vol 187:7 \yr 1996 \pages 951--980
\endref

\ref \key Ak96' \by P.~M.~Akhmetiev \paper Mapping an $n$-sphere in a
$2n$-Euclidean space: its realization (in Russian)
\jour Trudy Mat. Inst. im. Steklova \vol 212 \yr 1996 \pages 37--45
\moreref English transl. in Proc. Math. Steklov Inst., {\bf 212} (1996), 32--39
\endref

\ref \key Ak96'' \by P.~M.~Akhmetiev \paper Prem-mappings, triple points of
orientable surface and Rohlin Signature Theorem (in Russian)
\jour Mat.~Zametki \vol 59:6 \yr 1996 \pages 803--810
\moreref English transl.~in Math.~Notes {\bf 59}:6 (1996), 581--585
\endref

\ref \key Ak00 \by P. M. Akhmetiev \paper Embeddings of compacta, stable
homotopy groups of spheres and singularity theory  (in Russian)
\jour Uspekhi Mat. Nauk \vol 55:3 \yr 2000 \pages 3--62
\moreref English transl.~in Russ. Math.~Surveys {\bf 55}:3 (2000)
\endref

\ref \key Al24 \by J.~W.~Alexander
\paper On the subdivision of 3-space by polyhedron
\jour Proc. Nat. Acad. Sci. USA \vol 10 \yr 1924 \pages 6--8
\endref

\ref \key ARS01 \by P. Akhmetiev, D. Repov\v s and A. Skopenkov
\paper Embedding products of low--dimensional manifolds in $\R^m$
\jour Topol. Appl. \vol 113 \yr 2001 \pages 7--12
\endref

\ref \key ARS02 \by P. Akhmetiev, D. Repov\v s and A. Skopenkov \paper
Obstructions to approximating maps of $n$-surfaces to $\R^{2n}$ by embeddings
\jour Topol. Appl. \vol 123 \yr 2002 \pages 3--14
\endref

\ref \key Ba75 \by D. R. Bausum
\paper Embeddings and immersions of manifolds in Euclidean space
\jour Trans. AMS \vol 213 \yr 1975 \pages 263--303
\endref

\ref \key BH70 \by J. Boechat and A. Haefliger \pages 156--166
\paper Plongements differentiables de varietes de dimension 4 dans $\R^7$
\yr 1970 \vol  \jour Essays on topology and related topics (Springer,1970)
\endref

\ref \key BKK02 \by M. Bestvina, M. Kapovich and B. Kleiner
\paper Van Kampen's embedding obstruction for discrete groups
\jour Invent. Math. \vol 150:2 \yr 2002 \pages 219--235
\endref

\ref \key BM99 \by R. L. Bryant and W. Mio
\paper Embeddings of homology manifolds in codimension $\ge3$
\yr 1999 \jour Topology \pages 811--821 \vol 38:4
\endref

\ref \key BM00 \by R. L. Bryant and W. Mio
\paper Embeddings in generalized manifolds
\yr 2000 \jour Trans. Amer. Math. Soc. \pages 1131--1137 \vol 352:3
\endref

\ref \key Bo33 \by K. Borsuk
\paper \"Uber stetige Abbildungen der euklidischen R\" aume
\jour Fund. Math. \vol 21 \yr 1933 \pages 236--246
\endref

\ref \key Bo71 \by J. Boechat \pages 141--161
\paper Plongements differentiables de varietes de dimension $4k$ dans
$\R^{6k+1}$
\yr 1971 \vol 46:2 \jour Comment. Math. Helv.
\endref

\ref \key Br68 \by W. Browder
\paper Embedding smooth manifolds
\yr 1968 \jour Proc. Int. Congr. Math. Moscow 1966
\pages 712--719 \endref

\ref \key Br71 \by R. L. Brown
\paper Immersions and embeddings up to cobordism
\yr 1971 \jour Canad. J. Math \pages 1102--1115\vol 23 \endref

\ref \key Br72 \by J. L. Bryant
\paper Approximating embeddings of polyhedra in codimension 3
\jour Trans. Amer. Math. Soc. \vol 170 \yr 1972 \pages 85--95
\endref

\ref \key Ca86 \by A. Casson \pages 201-244
\paper Three lectures on new infinite constructions in 4-dimensional manifolds
\inbook  A la Recherche de la Topologie Perdue, {\rm Progress in Mathematics,
62} \eds L. Guillou, A. Marin \publ Birkhauser \publaddr Boston
\yr 1986
\endref

\ref \key Ch69 \by A. V. Chernavskii
\paper Piecewise linear approximations of embeddings of cells and spheres
in codimensions higher than two
\jour Mat. Sb. \vol 80 \issue (122) \yr 1969 \pages 339--364
\finalinfo Math. USSR Sb. 9 (1969), 321--344 \endref

\ref \key Cl34 \by S.~Claytor
\paper Topological immersions of peanian continua in a spherical surface
\jour Ann. of Math. \vol 35 \yr 1934 \pages 809--835
\endref

\ref \key Cl37 \by S.~Claytor
\paper Peanian continua not embeddable in a spherical surface
\jour Ann. of Math. \vol 38 \yr 1937 \pages 631--646
\endref

\ref \key CF60 \by P. E. Conner and E. E. Floyd
\paper Fixed points free involutions and equivariant maps
\jour Bull. AMS \vol 66 \yr 1960 \pages 416--441
\endref

\ref \key CG83 \by J. H. Conway and C. M. A. Gordon
\paper Knots and links in spatial graphs
\jour J. Graph Theory  \vol 7 \yr 1983 \pages 445--453
\endref

\ref \key Co69 \by G. Cooke
\paper Embedding certain complexes up to homotopy in Euclidean space
\jour Ann. of Math. \vol 90 \yr 1969 \pages 144--156
\endref

\ref \key Co85 \by R. L. Cohen
\paper The immersion conjecture for differentiable manifolds
\jour Ann. of Math. \vol 122 \yr 1985 \pages 237--328
\endref

\ref \key CRS98 \by A.~Cavicchioli, D.~Repov\v s and A.~B.~Skopenkov
\paper Open problems on graphs, arising from geometric topology
\jour Topol. Appl. \yr 1998 \vol 84 \pages 207--226
\endref

\ref  \key CRS04 \by M. Cencelj, D. Repov\v s and A. Skopenkov
\paper On the Browder-Levine-Novikov embedding theorems
\jour Trudy MIRAN \vol 247 \yr 2004  \pages 280--290
\moreref Proc. Steklov Math. Inst. {\bf 247} (2004)
\endref

\ref  \key CRS05 \by M. Cencelj, D. Repov\v s and M. Skopenkov
\paper Classification of framed links in 3-manifolds
\jour submitted \vol \yr   \pages
\endref

\ref  \key CRS \by M. Cencelj, D. Repov\v s and M. Skopenkov
\paper Knotted tori and the $\beta$-invariant
\jour preprint \vol \yr   \pages
\endref

\ref \key CS79 \by S. E. Cappell and J. L. Shaneson \pages 301--303
\paper Imbeddings and immersions of 4-dimensional manifolds in $\R^6$
\inbook  Geometric topology
\eds J. C. Cantrell \publ Academic press \publaddr New York
\yr 1979 \endref


\ref \key CW78 \by F. X. Conolly and B. Williams
\paper Embedding up to homotopy and geometric suspension of manifolds
\jour Quart. J. Math. Oxford (2) \vol 29 \yr 1978 \pages 385--401
\endref

\ref \key Da86 \by R. J. Daverman
\book Decompositions of manifolds
\yr 1986 \publ Academic Press \publaddr New York
\endref

\ref \key Do87 \by S. K. Donaldson
\paper The orientation of Yang-Mills moduli spaces and 4-manifold topology
\jour J. Diff. Geom. \vol 26 \yr 1987 \pages 397--428
\endref

\ref \key DP97 \by T. J. Dodson and P. E. Parker
\book A User's Guide to Algebraic Topology
\publ Kluwer \publaddr Doldrecht-Boston-London \yr 1997
\endref

\ref \key DRS91 \by A.~N.~Drani\v snikov, D.~Repov\v s and E.~V.~\v S\v cepin
\paper On intersection of compacta of complementary dimension in Euclidean space
\jour Topol. Appl \vol 38 \yr 1991 \pages 237--253
\endref

\ref \key DRS93 \by A.~N.~Drani\v snikov, D.~Repov\v s and E.~V.~\v S\v cepin
\paper On intersection of compacta in Euclidean space: the metastable case
\jour Tsukuba J.Math. \vol 17 \yr 1993 \pages 549--564
\endref

\ref \key Ed75 \by R. B. Edwards
\paper The equivalence of close piecewise linear embeddings
\jour General Topology and its Applications \vol 5 \yr 1975 \pages 147--180
\endref

\ref \key FF89 \by A. T. Fomenko and D. B. Fuchs
\book A course in homotopy theory \lang in Russian
\publ Nauka \publaddr Moscow \yr 1989
\endref

\ref \key FKT94 \by M.~H.~Freedman, V.~S.~Krushkal and P.~Teichner
\paper Van Kampen's embedding obstruction is incomplete for 2-complexes
in~$\R^4$ \jour Math. Res. Letters \vol 1 \yr1994 \pages 167--176
\endref

\ref \key FKV87  \by S. M. Finashin, M. Kreck and O. Ya. Viro
\paper Exotic knottings of surfaces in the 4-sphere
\jour Bull. Amer. Math. Soc. \vol 17:2 \yr 1987 \pages 287--290
\endref

\ref \key FKV88  \by S. M. Finashin, M. Kreck and O. Ya. Viro
\paper Non-diffeomorphic but homeomorphic knottings of surfaces in the 4-sphere
\jour Lecture Notes in Math.  \vol 1346 \yr 1988 \pages 157--198
\endref

\ref \key Fl34 \by A. Flores
\paper \"Uber $n$-dimensionale Komplexe die im $E^{2n+1}$ absolute
          Selbstverschlungen sind
\jour Ergeb. Math. Koll. \yr 1934 \vol 6 \pages 4--7
\endref

\ref \key FQ90 \by M.~H.~Freedman and F.~S.~Quinn
\book Topology of 4-Manifolds
\publ Princeton Univ. Press \publaddr Princeton \yr1990
\endref

\ref \key Fr87 \by G. K. Francis
\book A Topological Picturebook
\publ Springer-Verlag \publaddr Berlin \yr 1987
\endref

\ref \key FS59 \by D.~B.~Fuchs and A.~S.~Schwarz
\paper Cyclic powers of polyhedra and the imbedding problem
\jour Dokl. Akad. Nauk SSSR \vol 125 \yr 1959 \pages 285--288
\lang in Russian \endref

\ref \key Fu94 \by F. Fuquan \pages 447--454
\paper Embedding four manifolds in $\R^7$
\yr 1994 \vol 33:3 \jour Topology
\endref

\ref \key Fu02 \by F. Fuquan \pages 927--930
\paper Orientable 4-manifolds topologically embed into $\R^7$
\yr 2002 \vol 41\jour Topology
\endref

\ref \key Ga92 \by M. Galecki
\book On embeddability of CW-complexes in Euclidean space
\publ preprint \publaddr Univ. of Tennessee, Knoxville \yr 1992 \vol \endref

\ref \key Gi71 \by S. Gitler
\paper Embedding and immersion of manifolds
\yr 1971 \vol 22 \jour Proc. Symp. Pura Appl. Math., AMS, Providence
\pages 87--96 \endref

\ref \key Gl63 \by H.~Gluck
\paper Unknotting $S^1$ in $S^4$
\jour Bull. Amer. Math. Soc. \vol 69:1 \yr 1963 \pages 91--94
\endref

\ref \key Gl68 \by H.~Gluck \paper Geometric characterization of
differentiable manifolds in Euclidean space, II \jour Michigan
Math.~J. \vol 15:1 \yr 1968 \pages 33--50
\endref

\ref \key GMR94 \by D.~Gillman, S.~V.~Matveev and D.~Rolfsen
\paper Collapsing and reconstruction of manifolds \jour Contemp.
Math. \vol 164 \yr 1994 \pages 35--39
\endref

\ref \key Go72 \by C.~M.~A.~Gordon
\paper Embeddings of PL-manifolds with boundary
\jour Proc. Camb. Phil. Soc. \vol 72 \yr 1972 \pages 21--25
\endref

\ref \key Gr86 \by M. Gromov
\book Partial differential relations
\publ Ergebnisse der Mathematik und ihrer Grenzgebiete (3), Springer
Verlag \publaddr Berlin-New York \yr 1986
\endref

\ref \key GR92 \by D.~Gillman and D.~Rolfsen
\paper 3-manifolds embed in small 3-complexes
\jour Intern. J. Math. \vol 3 \yr 1992 \pages 179--183
\endref

\ref  \key GS06 \by D. Goncalves and A. Skopenkov
\paper Embeddings of homology equivalent manifolds with boundary
\jour Topol. Appl. \vol   \yr 2006  \pages to appear
\endref

\ref  \key Gu53 \by V. K. A. M. Gugenheim
\paper Piecewise linear isotopy and embedding of elements and spheres, I
\jour Proc. of the London Math. Soc. \vol 3 \yr 1953 \pages 29--53
\endref

\ref \key GW99 \by T. Goodwillie and M. Weiss
\paper Embeddings from the point of view of immersion theory, II
\jour Geometry and Topology \vol 3 \yr 1999 \pages 103--118
\endref

\ref \key Ha61 \by A. Haefliger \pages 47--82
\paper Plongements differentiables de varietes dans varietes
\yr 1961 \vol 36 \jour Comment. Math. Helv. \endref

\ref \key Ha62 \by A.~Haefliger
\paper Knotted $(4k-1)$-spheres in $6k$-space
\jour Ann. of Math. \vol 75 \yr 1962 \pages 452--466
\endref

\ref \key Ha62' \by A. Haefliger \pages 241--244
\paper Differentiable links
\yr 1962 \vol 1 \jour Topology
\endref

\ref \key Ha62'' \by A. Haefliger \pages
\paper Plongements de varietes dans le domain stable
\yr 1962 \vol 245 \jour Seminare Bourbaki
\endref

\ref \key Ha63 \by A.~Haefliger
\paper Plongements differentiables dans le domain stable
\jour Comment. Math. Helv. \vol 36 \yr 1962-63 \pages 155--176
\endref

\ref \key Ha66 \by A. Haefliger
\paper Differentiable embeddings of $S^n$ in $S^{n+q}$ for $q>2$
\pages 402--436 \jour Ann\. Math., Ser.3 \vol 83 \yr 1966
\endref

\ref \key Ha66' \by A.~Haefliger
\paper Enlacements de spheres en codimension superiure a 2
\jour Comment. Math. Helv. \vol 41 \yr 1966-67 \pages 51--72
\endref

\ref \key Ha67 \by A. Haefliger \pages 221--240
\paper Lissage des immersions-I
\yr 1967 \vol 6 \jour Topology
\endref

\ref \key Ha68 \by A. Haefliger \pages 437--445
\paper Knotted spheres and related geometric topic
\yr 1968 \vol \jour In: Proc. of the Int. Cong. Math., Moscow, Mir
\endref

\ref \key Ha68' \by D.~D.~J.~Hacon
\paper Embeddings of $S^p$ in $S^1\times S^q$ in the metastable range
\jour Topology \vol 7 \yr 1968 \pages 1--10
\endref

\ref \key Ha69 \by L.~S.~Harris
\paper Intersections and embeddings of polyhedra
\jour Topology \vol 8 \yr 1969 \pages 1-26
\endref

\ref \key Ha  \by A. Haefliger \pages
\paper Lissage des immersions-II
\yr 1966 \vol  \jour preprint
\endref

\ref \key Ha84 \by N. Habegger \pages 123--130
\paper Obstruction to embedding disks II: a proof of a conjecture by Hudson
\yr 1984 \vol 17 \jour Topol. Appl. \endref

\ref \key Ha86 \by N. Habegger
\paper Classification of links in $S^1\times\R^{n-1}$
\yr 1986 \vol 25:3 \jour Topology \pages 253--260
\endref

\ref \key HH62 \by A. Haefliger and M. W. Hirsch \pages 231--241
\paper Immersions in the stable range
\yr 1962 \vol 75:2 \jour Ann. of Math.
\endref

\ref \key HH63 \by A. Haefliger and M. W. Hirsch \pages 129--135
\paper On existence and classification of differential embeddings
\yr 1963 \vol 27 \jour Topology
\endref

\ref \key HK98 \by N. Habegger and U. Kaiser \pages
\paper Link homotopy in 2--metastable range
\yr 1998 \vol 37:1 \jour Topology \pages 75--94
\endref

\ref \key Hi61 \by M.~W.~Hirsch
\paper The embedding of bounding manifold in euclidean space
\jour Ann. of Math. \vol 74 \yr 1961 \pages 494--497
\endref

\ref \key Hi65 \by M.~W.~Hirsch
\paper On embedding 4-manifolds in $\R^7$
\jour Proc. Camb. Phil. Soc. \vol 61 \yr 1965
\endref

\ref \key Hi66 \by M.~W.~Hirsch
\paper Embeddings and compressions of polyhedra and smooth manifolds
\jour Topology \vol 4:4 \yr 1966 \pages 361--369
\endref

%\ref \key Hi \by S. Hirose
%\paper On diffeomorphisms over surfaces trivially embedded in the 4-sphere
%\yr 2002 \vol \jour preprint  \pages
%\endref

\ref \key HL71 \by J.~F.~P.~Hudson and W. B. R. Lickorish
\paper Extending piecewise linear concordances
\jour Quart. J. Math. (2) \vol 22 \yr 1971 \pages 1--12
\endref

\ref \key HLS65 \by W. C. Hsiang, J. Levine and R. H. Sczarba \pages 173--181
\paper On the normal bundle of a homotopy sphere embedded in Euclidean space
\yr 1965 \vol 3 \jour Topology \endref

\ref \key Ho69 \by K.~Horvati\v c
\paper Embedding manifolds with low-dimensional spine
\jour Glasnik Mat. \vol 4(24):1 \yr 1969 \pages 101--116
\endref

\ref \key Ho71 \by K.~Horvati\v c
\paper On embedding polyhedra and manifolds
\jour Trans. AMS \vol 157 \yr 1971 \pages 417--436
\endref

\ref \key HS64 \by W. C. Hsiang and R. H. Szarba
\paper On the embeddability and non-embeddability of sphere bundles over
spheres
\jour Ann. of Math. \yr 1964 \vol 80:2 \pages 397--402
\endref

\ref \key Hu59 \by S. T. Hu
\book Homotopy Theory \finalinfo Academic Press, New York, 1959
\endref

\ref \key Hu60 \by S. T. Hu
\paper Isotopy invariants of topological spaces
\jour Proc. Royal Soc. \vol A255 \yr 1960 \pages 331--366
\endref

\ref \key Hu63 \by J.~F.~P.~Hudson
\paper Knotted tori \jour Topology \vol 2 \yr 1963 \pages 11--22
\endref

\ref \key Hu63' \by J.~F.~P.~Hudson
\paper Non-embedding theorem \jour Topology \vol 2 \yr 1963 \pages 123--128
\endref

\ref \key Hu66 \by J.~F.~P.~Hudson
\paper Extending piecewise linear isotopies
\jour Proc. London Math. Soc. (3) \vol 16 \yr 1966 \pages 651--668
\endref

\ref \key Hu67 \by  J.~F.~P.~Hudson
\paper PL embeddings \jour Ann. of Math. \vol 85:1 \yr 1967 \pages 1--31
\endref

\ref \key Hu69 \by J. F. P. Hudson
\book Piecewise-Linear Topology
\bookinfo \publ Benjamin \publaddr New York, Amsterdam \yr 1969
\endref

\ref \key Hu70 \by J.~F.~P.~Hudson
\paper Concordance, isotopy and diffeotopy
\jour Ann. of Math. \vol 91:3 \yr 1970 \pages 425--448
\endref

\ref \key Hu70' \by J. F. P. Hudson \pages 407--415 \paper
Obstruction to embedding disks \yr 1970 \vol \jour In: Topology of
manifolds
\endref

\ref \key Hu72 \by J.~F.~P.~Hudson
\paper Embeddings of bounded manifolds \jour Proc. Camb. Phil. Soc.
\vol 72 \yr 1972 \pages 11--20
\endref

\ref \key HZ64 \by J.~F.~P.~Hudson and E.~C.~Zeeman
\paper On regular neighborhoods \jour Proc. Lond. Math. Soc. (3)
\vol 14 \yr 1964 \pages 719--745
\paperinfo Correction to `On regular neighborhoods',
Proc. Lond. Math. Soc. (3), 21 (1970) 513--524.
\endref

\ref \key Ir65 \by M.~C.~Irwin
\paper Embeddings of polyhedral manifolds \jour Ann. of Math. (2)
\vol 82 \yr 1965 \pages 1--14
\endref

\ref \key Ka32 \by E.~R.~van~Kampen
\paper Komplexe in euklidische Raumen \jour Abb. Math. Sem. Hamburg
\vol 9 \yr 1932 \pages72--78 \moreref berichtigung dazu, 152--153
\endref

\ref \key Ke59 \by M. Kervaire
\paper An interpretation of G. Whitehead's generalization of H. Hopf's
invariant
\yr 1959 \vol 62 \jour Ann. of Math. \pages 345--362
\endref

\ref \key Ke69 \by M.~A.~Kervaire
\paper Smooth homology spheres and their fundamental groups
\jour Trans. Amer. Math. Soc. \vol 144 \yr 1969 \pages 67--72
\endref

\ref \key Ke79 \by C.~Kearton
\paper Obstructions to embeddings and isotopy in the metastable range
\jour Math. Ann. \vol 243 \yr 1979 \pages 103--113
\endref

\ref \key Ki84 \by R.~Kirby \paper 4-manifold problems
\jour Contemp. Math. \vol 35 \yr 1984 \pages 513--528
\endref

\ref \key Ki89 \by R.~C.~Kirby
\paper The Topology of 4-Manifolds \jour Lect. Notes Math. \vol 1374
\publ Springer-Verlag \publaddr Berlin \yr1989
\endref

\ref \key KM61 \by A.~Kervaire and J.~W.~Milnor
\paper On 2-spheres in 4-manifolds \yr 1961 \vol 47
\jour Proc. Nat. Acad. Sci. USA \pages 1651--1657
\endref

\ref \key Ko62 \by A. Kosinski
\paper On Alexander's theorem and knotted tori
\jour In: Topology of 3-Manifolds, Prentice-Hall, Englewood Cliffs, Ed.
M.~K.~Fort, N.J. \yr 1962 \pages 55--57
\endref

\ref \key Ko88 \by U. Koschorke
\paper Link maps and the geometry of their invariants
\yr 1988 \vol 61:4 \jour Manuscripta Math. \pages 383--415 \endref

\ref \key Ko90 \by U. Koschorke
\paper On link maps and their homotopy classification
\yr 1990 \vol 286:4 \jour Math. Ann. \pages 753--782 \endref

\ref \key Kr00 \by V. S. Krushkal
\paper Embedding obstructions and 4-dimensional thickenings of 2-complexes
\jour Proc. Amer. Math. Soc. \vol 128:12 \yr 2000 \pages 3683--3691
\endref

\ref  \key KS05 \by M. Kreck and A.  Skopenkov
\paper Inertia groups of smooth embeddings
\jour submitted \vol \yr \pages
\moreref http:// arxiv.org/ abs/ math.GT/ 0512594
\endref

\ref \key Ku30 \by K.~Kuratowski
\paper Sur les probl\`emes des courbes gauche en topologie
\jour Fund. Math. \vol 15 \yr 1930 \pages 271--283
\endref

\ref \key Ku00 \by V.~Kurlin
\paper Basic embeddings into products of graphs
\yr 2000 \jour Topol. Appl. \vol 102 \pages 113--137
\endref

\ref \key KW85 \by J. Keesling and D. C. Wilson
\paper Embedding $T^n$-like continua in Euclidean space
\jour Topol. Appl. \vol 21 \yr 1985 \pages 241--249
\endref

\ref \key La82 \by J.~Lannes
\paper La conjecture des immersions
\jour Ast\'erisque \vol 92/93 \yr 1982 \pages 331--346
\endref

\ref \key La96 \by M. Lackenby
\paper The Whitney Trick
\jour Topol. Appl. \vol 71 \yr 1996 \pages 115--118
\endref

\ref \key Le62 \by J.~Levine
\paper
\jour Notices of the Amer. Math. Soc.  \vol \yr 1962 \pages
\endref

\ref \key Le65 \by J.~Levine
\paper A classification of differentiable knots
\jour Ann. of Math.  \vol 82 \yr 1965 \pages 15--50
\endref

\ref \key Li65 \by W.~B.~R.~Lickorish
\paper The piecewise linear unknotting of cones \yr 1965 \vol 4
\jour Topology \pages 67--91
\endref

\ref \key Li75 \by A. Liulevicius
\paper Immersions up to cobordism
\jour Illinois J. Math \yr 1975 \vol 19 \pages 149--164
\endref

\ref \key LS69 \by W.~B.~R.~Lickorish and L.~C.~Siebenmann
\paper Regular neighborhoods and the stable range
\jour Trans. AMS \vol 139 \yr 1969 \pages 207--230
\endref

\ref \key MA41 \by S.~McLane and V.~W.~Adkisson
\paper Extensions of homeomorphisms on the spheres
\jour Michig. Lect. Topol. \publ Ann~Arbor \yr1941  \pages 223--230
\endref

\ref \key Ma60 \by W. S. Massey
\paper On the Stiefel--Whitney classes of a manifold, 1
\jour Amer. J. Math \vol 82 \yr 1960 \pages 92--102
\endref

\ref \key Ma62 \by W. S. Massey
\paper On the Stiefel--Whitney classes of a manifold, 2
\jour Proc. AMS \vol 13 \yr 1962 \pages 938--942
\endref

\ref \key Ma90 \by W. S. Massey \pages 269--300
\paper Homotopy classification of 3-component links of codimension
greater than 2
\yr 1990 \vol 34 \jour Topol. Appl.
\endref

\ref \key Ma80 \by R. Mandelbaum
\paper Four-Dimensional Topology: An introduction
\jour Bull. Amer. Math. Soc. (N.S.) \vol 2 \yr 1980 \pages 1-159
%The University of Rochester, Department of Mathematics, 1978
\endref

\ref \key Mc67 \by M.~C.~McCord
\paper Embedding P-like compacta in manifolds
\jour Canad. J.~Math. \vol 19 \yr 1967  \pages 321--332
\endref

\ref \key Me02 \by S. Melikhov
\paper On maps with unstable singularities
\yr 2002 \jour Topol. Appl.  \vol 120 \pages 105--156
\endref

\ref \key Me04 \by S. Melikhov
\paper Sphere eversions and realization of mappings
\moreref http:// front.math.ucdavis.edu/ math.GT/ 0305158
\yr 2004 \jour Trudy MIRAN  \vol 247 \pages 1--21, 159--181
\endref

\ref \key Mi54 \by J. Milnor \pages 177--195
\paper Link groups \yr 1954 \vol 59 \jour Ann. of Math.
\endref

\ref \key Mi65 \by J.~Minkus
\paper On embeddings of highly connected manifolds
\jour Trans. of the Amer. Math. Soc. \vol 115 \yr 1965 \pages 525--540
\endref

\ref \key Mi70 \by R. T. Miller
\paper Close isotopies on piecewise-linear manifolds
\jour Trans. of the Amer. Math. Soc. \vol 151:2 \yr 1970 \pages 597--628
\endref

\ref \key Mi72 \by R. T. Miller \pages 406--416
 \paper Approximating codimension 3 embeddings
 \yr 1972 \vol 95 \jour Ann. of Math \endref

\ref \key Mi72' \by R. J. Milgram
\paper On the Haefliger's groups
\yr 1972 \vol 78:5 \jour Bull. of the Amer. Math. Soc. \pages 861--865
\endref

\ref \key Mi75 \by K. C. Millett
\paper Piecewise linear embeddings of manifolds
\jour Illinois J. Math. \vol 19 \yr 1975 \pages 354--369 \endref

\ref \key Mi97 \by P. Minc
\paper Embedding simplicial arcs into the plane
\jour Topol. Proc. \vol 22 \pages 305--340 \yr 1997
\endref

\ref \key MK58 \by J. Milnor and M. Kervaire
\paper Bernoulli numbers, homotopy groups of spheres and a theorem of Rohlin
\yr 1958 \vol  \jour Proc. Int. Cong. Math., Edinburgh  \pages
\endref

%\ref  \key Mo83 \by J. M. Montesinos
%\paper On twins in the four-sphere, I
%\jour Quart. J. Math Oxford (2) \vol 34  \yr 1983  \pages 171-199
%\endref

\ref  \key MR71 \by R. J. Milgram and E. Rees
\paper On the normal bundle to an embedding
\jour Topology \vol 10  \pages 299--308  \yr 1971
\endref

\ref \key MR86 \by W. Massey and D. Rolfsen \paper Homotopy
classification of higher dimensional links \yr 1986 \vol 34 \jour
Indiana Univ. Math. J. \pages 375--391
\endref

\ref  \key MRS03 \by J. Male\v si\v c, D. Repov\v s and A. Skopenkov
\paper On incompleteness of the deleted product obstruction for embeddings
\jour Bol. Soc. Mat. Mexicana (3) \vol 9  \yr 2003  \pages 165--170
\endref

\ref \key MS67 \by S. Marde\v si\' c and J. Segal \pages 171--182
\paper $\varepsilon$-mappings and generalized manifolds
\yr 1967 \vol 14 \jour Michigan Math. J.\endref

\ref \key MS74 \by J. W. Milnor and J. D. Stasheff
\book Characteristic Classes, {\rm Ann. of Math. St. \bf 76} \yr 1974
\publ Princeton Univ. Press \publaddr Princeton, NJ
\endref

\ref  \key MS04 \by J. Mukai and A. Skopenkov
\paper A direct summand in a homotopy group of the mod 2 Moore space
\jour Kyushu J. Math. \vol 58:1  \yr 2004  \pages 203--209
\endref

\ref \key MT95 \by M. Mahowald and R. D. Thompson
\paper  The EHP Sequence and Periodic Homotopy
\jour in: Handbook of Algebraic Topology, ed. I. M. James \yr 1995
\pages Elsevier Science B. V.
\endref

\ref \key Ne68 \by L.~Neuwirth
\paper An algorithm for the construction of 3-manifolds from 2-complexes
\jour Proc. Camb. Phil. Soc. \vol 64 \yr 1968 \pages 603--613
\endref

\ref \key Ne98 \by S. Negami
\paper Ramsay-type theorem for spatial graphs
\jour Jour. Comb. Th. Ser B \vol 72 \yr 1998 \pages 53--62
\endref

\ref \key No61 \by S. P. Novikov
\paper Imbedding of simply-connected manifolds in Euclidean space
\lang in Russian \jour Dokl. Akad. Nauk SSSR \vol 138 \yr 1961 \pages 775--778
\endref

\ref \key Oo86 \by H. Ooshima
\paper Whitehead products in the Stiefel manifolds and Samelson products
in classical groups
\jour Adv. Stud. in Pure Math. \vol 9 \yr 1986 \pages 237--258
\endref

\ref \key Pa56 \by G. F. Paechter
\paper On the groups $\pi_r(V_{mn})$
\jour Quart. J. Math. Oxford, Ser.2
\moreref \paper I    \vol 7:28  \yr 1956 \pages 249--265
\moreref \paper II   \vol 9:33  \yr 1958 \pages 8--27
\moreref \paper III  \vol 10:37 \yr 1959 \pages 17--37
\moreref \paper IV   \vol 10:40 \yr 1959 \pages 241--260
\moreref \paper V    \vol 11:41 \yr 1960 \pages 1--16
\endref

\ref \key Po42 \by L.~S.~Pontryagin
\paper Characteristic cycles of smooth manifolds
\jour Dokl. Akad. Nauk SSSR \vol 35:2 \yr 1942 \pages 35--39 \lang in Russian
\endref

\ref \key Po85 \by M. M. Postnikov
\book Homotopy theory of CW-complexes
\publ Nauka \publaddr Moscow, 1985 \lang in Russian
\endref

\ref \key Pr66 \by T. M. Price
\paper Equivalence of embeddings of $k$-complexes in $E^n$ for $n\ge2k+1$
\jour Michigan Math. J. \vol 13 \yr 1966 \pages 65--69
\endref

\ref \key Pr04 \by V. V. Prasolov
\book Introduction into combinatorial and differential topology
\publ AMS \publaddr Providence, RI \yr 2004
\endref

\ref \key Pr \by V. V. Prasolov
\book Elements of homology theory (in Russian)
\publ  \publaddr  \yr 2005
\endref

\ref  \key PS05 \by V. Prasolov and M. Skopenkov
\paper Ramsay link theory (in Russian)
\jour Mat. Prosveschenie \vol 9 \yr 2005  \pages
\endref

\ref \key PWZ61 \by R.~Penrose, J.~H.~C.~Whitehead and E.~C.~Zeeman
\paper Embeddings of manifolds in a Euclidean space
\jour Ann. of Math. (2) \vol 73 \yr 1961 \pages 613--623
\endref

\ref \key RBS99  \by  D.~Repov\v s, N. Brodsky and A.~B.~Skopenkov
\paper A classification of 3-thickenings of 2-polyhedra
\jour Topol. Appl. \vol 94 \yr 1999 \pages 307--314
\endref

\ref \key Re71 \by E. Rees \pages 152--156
\paper Some embeddings of Lie groups in Euclidean spaces
\yr 1971 \vol 18 \jour Mathematika
\endref

\ref \key Re90 \by E. Rees \pages 72--79
\paper Problems concerning embeddings of manifolds
\yr 1990 \vol 19:1 \jour Advances in Math.
\endref

\ref \key Ro65 \by V. A. Rohlin \paper The embedding of
non-orientable three-dimensional manifolds in the five-dimensional
Euclidean space \lang in Russian \jour Dokl. Akad. Nauk SSSR \vol
160 \yr 1965 \pages 549--551
\endref

\ref \key RS72 \by C. P. Rourke and B. J. Sanderson
\book Introduction to Piecewise-Linear Topology,
{\rm Ergebn. der Math. \bf 69}
\publ Springer-Verlag \publaddr Berlin \yr 1972
\endref

\ref \key RS96 \by D.~Repov\v s and A.~B.~Skopenkov
\paper Embeddability and isotopy of polyhedra in Euclidean spaces
\jour Trudy Math. Inst. Ross. Akad. Nauk \vol 212 \yr 1996 \pages
\moreref Proc. of the Steklov Inst. Math. \vol 212 \yr 1996 \pages 173--188
\endref

\ref \key RS98  \by D. Repov\v s and A.~B.~Skopenkov
\paper A deleted product criterion for approximability of a map by embeddings
\jour Topol. Appl. \vol 87 \yr 1998 \pages 1--19 \endref

\ref \key RS99 \by D. Repovs and A. Skopenkov  \paper New results
on embeddings of polyhedra and manifolds into Euclidean spaces
\lang in Russian \yr 1999 \vol 54:6 \jour Uspekhi Mat. Nauk \pages 61--109
\moreref \paper English transl. \jour Russ. Math. Surv. \pages 1149--1196
\endref

\ref \key RS99'  \by D. Repov\v s and A.~B.~Skopenkov
\paper Borromean rings and embedding obstructions (in Russian)
\jour Trudy Math. Inst. Ross. Akad. Nauk \vol 225 \yr 1999 \pages 331--338
\moreref English transl.:\jour Proc. of the Steklov Inst. Math. \vol 225
\yr 1999 \pages 314--321
\endref

\ref  \key RS00 \by D. Repov\v s and A.  Skopenkov
\paper Obstruction theory for beginners (in Russian)
\jour Mat. Prosveschenie \vol 4 \yr 2000  \pages
\endref

\ref \key RS01 \by C. Rourke and B. Sanderson
\paper The compression theorem. I, II
\pages 399--429, 431--440 \jour Geom. Topol. (electronic) \vol 5 \yr 2001
\endref

\ref  \key RS01' \by D. Repov\v s and A. Skopenkov
\paper On contractible $n$-dimensional compacta, non-embeddable into $\R^{2n}$
\jour Proc. Amer. Math. Soc. \vol 129 \yr 2001 \pages 627--628
\endref

\ref  \key RS02 \by D. Repov\v s and A. Skopenkov \paper On
projected embeddings and desuspension of the $\alpha$-invariant
\jour Topol. Appl. \vol 124  \yr 2002  \pages 69--75
\endref

\ref  \key RS02' \by D. Repov\v s and A.  Skopenkov
\paper Characteristic classes for beginners (in Russian)
\jour Mat. Prosveschenie \vol 6 \yr 2002  \pages 60-77
\endref

\ref \key RSS95 \by D.~Repov\v s, A.~B.~Skopenkov and E.~V.~\v S\v cepin
\paper On embeddability of $X\times I$ into Euclidean space
\jour Houston J.~Math. \vol 21 \yr 1995 \pages 199--204 \endref

\ref \key RSS95' \by D.~Repov\v s, A.~B.~Skopenkov  and E.~V.~\v S\v cepin
\paper On uncountable collections of continua and their span
\jour Colloq. Math. \vol 69:2 \yr 1995 \pages 289--296 \endref

\ref  \key RST95  \by N. Robertson, P. D. Seymor and R. Thomas
\paper Sachs' linkless embedding conjecture
\jour J. Combin. Theory, Ser. B \vol 64:2 \yr 1995  \pages 185--227
\endref

\ref \key Ru73 \by T.~B.~Rushing
\book Topological Embeddings
\yr 1973 \publ Academic Press \publaddr New York
\endref

\ref \key Ru82 \by D. Ruberman
\paper Imbeddings four-manifolds and slicing links
\jour Math. Proc. Camb. Phil. Soc \vol 91 \yr 1982 \pages 107--110
\endref

\ref \key Sa65 \by  R. De Sapio
\paper  Embedding $\pi$-manifolds
\jour Ann. of Math. (2) \vol  82 \yr 1965 \pages  213--224
\endref

\ref \key Sa81 \by H. Sachs
\paper On spatial representation of finite graphs
\jour in: Finite and infinite sets, Colloq. Math. Soc. Janos Bolyai,
North Holland, Amsterdam
\yr 1981 \vol 37 \pages
\endref

\ref \key Sa91 \by K.~S.~Sarkaria
\paper Kuratowski complexes \jour Topology \vol 30 \yr 1991 \pages 67--76
\endref

\ref \key Sa91' \by K.~S.~Sarkaria
\paper A one-dimensional Whitney trick and Kuratowski's graph planarity
criterion \jour Israel J.~Math. \vol 73 \yr 1991 \pages 79--89
\endref

\ref \key Sc77 \by M. Scharlemann
\paper Isotopy and cobordism of homology spheres in spheres
\yr 1977 \vol 16:3 \jour J. London Math. Soc., Ser. 2 \pages 559--567
\endref

\ref \key Sc84 \by E.~V. Schepin
\paper Soft mappings of manifolds
\jour Russian Math. Surveys \vol 39:5 \yr 1984 \pages 209--224 (in Russian)
\endref

\ref \key Sh57 \by A.~Shapiro
\paper Obstructions to the embedding of a complex in a Euclidean space,
I, The first obstruction
\jour Ann. of Math. (2) \vol 66 \yr 1957 \pages 256--269
\endref

\ref \key Si69 \by K.~Sieklucki
\paper Realization of mappings
\jour Fund. Math. \vol 65 \yr 1969 \pages 325--343
\endref

\ref \key Sk94 \by A. Skopenkov
\paper A geometric proof of the Neuwirth theorem on thickenings of
2-polyhedra \jour Mat. Zametki \vol 56:2 \yr 1994 \pages 94--98 \lang in
Russsian \moreref English transl.: Math. Notes, 58:5 (1995), 1244--1247
\endref

\ref \key Sk95 \by A. Skopenkov
\paper A description of continua basically embeddable in $\R^2$
\jour Topol. Appl. \vol 65 \yr 1995 \pages 29--48 \endref

\ref \key Sk97 \by A.~B.~Skopenkov
\paper On the deleted product criterion for embeddability of manifolds in
$\R^m$ \jour Comment. Math. Helv. \vol 72 \yr 1997 \pages 543--555 \endref

\ref \key Sk98  \by A.~B.~Skopenkov
\paper On the deleted product criterion for embeddability in $\R^m$
\jour Proc. Amer. Math. Soc. \vol 126:8 \yr 1998 \pages 2467--2476
\endref

\ref \key Sk00 \by A. Skopenkov
\paper On the generalized Massey--Rolfsen invariant for link maps
\yr 2000 \vol 165 \jour Fund. Math. \pages 1--15 \endref

\ref \key Sk02 \by A. Skopenkov
\paper On the Haefliger-Hirsch-Wu invariants for embeddings and immersions
\yr 2002 \vol 77 \jour Comment. Math. Helv. \pages 78--124
\endref

\ref \key Sk03 \by M. Skopenkov
\paper Embedding products of graphs into Euclidean spaces
\yr 2003 \vol 179 \jour Fund. Math. \pages 191--198
\endref

\ref \key Sk03' \by M. Skopenkov
\paper On approximability by embeddings of cycles in the plane
\yr 2003 \vol 134 \jour Topol. Appl. \pages 1--22
\endref

\ref \key Sk05 \by A. Skopenkov
\paper A new invariant and parametric connected sum of embeddings
\yr \vol \jour submitted to Fund. Math. \pages
\moreref http://arxiv.org/abs/math/0509621
\endref

\ref \key Sk05' \by A. Skopenkov
\paper On the Kuratowski graph planarity criterion
\jour Mat. Prosveschenie  \vol 9 \yr 2005 \pages 116--128
\moreref {\bf 10} (2006), 276--277
\endref

\ref \key Sk05'' \by A. Skopenkov \paper Algebraic topology from
geometrical point of view (in Russian) \jour http:// www.mccme.ru/
ium/ s05/ el\_al\_top.html  \vol  \yr  \pages
\endref

\ref \key Sk06 \by A. Skopenkov
\paper Classification of smooth embeddings of 3-manifolds in the 6-space
\yr \vol \jour  \pages submitted 
\moreref http://arxiv.org/abs/math/0603429
\endref

\ref \key Sk \by A. Skopenkov
\paper Classification of embeddings below the metastable range
\yr 2003 \vol \jour preprint \pages
\endref

\ref \key Sk' \by M. Skopenkov
\paper Explicit formulas for the groups of links and link maps
\yr \vol \jour preprint \pages
\endref

\ref \key Sm42 \by P. A. Smith
\paper Fixed points of periodic transformations
\finalinfo Appendix B in: Algebraic Topology by S.\ Lefschetz
(American Mathematical Society, New York, 1942) \pages 350-373
\endref

\ref \key Sm59 \by S.~Smale
\paper The classification of immersions of spheres in Euclidean spaces
\yr 1959 \vol 69 \jour Ann. of Math. (2) \pages 327--344
\endref

\ref \key Sp90 \by S.~Spie\D z \paper Imbeddings in $\R^{2m}$ of
$m$-dimensional compacta with $\dim(X\times X)< 2m$
\jour Fund. Math. \vol 134 \yr 1990 \pages 105--115
\endref

\ref \key SS83 \by E.~V.~Schepin and M.~A.~Shtanko
\paper A spectral criterion for embeddability of compacta in Euclidean space
\jour Proc. Leningrad Int. Topol. Conf. \publ Nauka \publaddr Leningrad
\yr 1983 \pages 135--142 \lang in Russian
\endref

\ref \key SS92 \by J.~Segal and S.~Spie\D z
\paper Quasi embeddings and embeddings of polyhedra in $\R^m$
\jour Topol. Appl. \vol 45 \yr 1992 \pages 275--282
\endref

\ref \key SSS98 \by J. Segal, A. Skopenkov and S. Spie\D z
\paper Embeddings of polyhedra in $\R^m$ and the deleted product obstruction
\jour Topol. Appl. \vol 85 \yr 1998 \pages 225--234
\endref

\ref \key St63 \by J.~Stallings
\paper On topologically unknotted spheres
\jour Ann. of Math. (2) \vol 77 \yr 1963 \pages 490--503
\endref

\ref \key St65 \by J. Stallings
\paper Homology and central series of groups
\jour J. of Algebra \vol 2 \yr 1965 \pages 170--181
\endref

\ref \key St89 \by Y.~Sternfeld
\paper Hilbert's 13th problem and dimension
\jour Lect. Notes Math \vol 1376 \yr 1989 \pages 1--49
\endref

\ref \key ST91 \by S.~Spie\D z and H.~Toru\'nczyk
\paper Moving compacta in $\R^m$ apart
\jour Topol. Appl. \vol 41 \yr 1991 \pages 193--204
\endref

\ref \key St \by J. Stallings
\paper The embedding of homotopy type into manifolds
\jour mimeographed notes, Princeton Univ.  \vol \yr 1965 \pages
\endref

\ref \key Sz82 \by A.~Sz\"ucs
\paper The Gromov--Eliashberg proof of Haefliger's theorem
\jour St. Sci. Math. Hung. \vol 17 \yr 1982 \pages 303--318
\endref

\ref \key Ta95 \by K.~Taniyama
\paper Homology classification of spatial embeddings of a graph
\jour Topol. Appl. \vol 65 \pages 205--228 \yr 1995
\endref

\ref \key Ti69 \by R. Tindell
\paper Knotting tori in hyperplanes
\jour in: Conf. on Topology of Manifolds, Prindle, Weber and Schmidt \vol
\yr 1969 \pages 147--153
\endref

\ref \key To62 \by H. Toda
\book Composition Methods in the Homotopy Groups of Spheres
\publ Princeton Univ. Press \publaddr Princeton \yr 1962
\endref

\ref \key Va92 \by V.~A.~Vassiliev
\book Complements of discriminants of smooth maps: Topology and applications
\publ Amer. Math. Soc. \publaddr Providence, RI \yr 1992
\endref

\ref \key Vr77  \by J. Vrabec
\paper Knotting a $k$-connected closed PL $m$-manifolds in $\R^{2m-k}$
\jour Trans. Amer. Math. Soc. \vol 233 \yr 1977  \pages 137--165
\endref

\ref \key Vr89  \by J. Vrabec
\paper Deforming of a PL submanifold of a Euclidean space into a hyperplane
\jour Trans. Amer. Math. Soc. \vol 312:1 \yr 1989 \pages 155--178
\endref

\ref \key Wa64 \by C.~T.~C.~Wall
\paper Differential topology, IV (theory of handle decompositions)
\jour Cambridge \vol  \yr 1964 \pages mimeographed notes
\endref

\ref \key Wa65 \by C.~T.~C.~Wall
\paper All 3-manifolds imbed in 5-space
\jour Bull. Amer. Math. Soc. \vol 71 \yr 1965 \pages 490--503
\endref

\ref \key Wa65' \by C.~T.~C.~Wall
\paper Unknotting spheres in codimension two and tori in codimension one
\jour Proc. Camb. Phil. Soc. \vol 61 \yr 1965 \pages 659--664
\endref

\ref \key Wa66 \by C.~T.~C.~Wall
\paper Classification problems in differential topology, IV, Thickenings
\jour Topology \vol 5 \yr 1966 \pages 73--94
\endref

\ref \key Wa66' \by C.~T.~C.~Wall
\paper Classification problems in differential topology, V,
On certain 6-manifolds
\jour Invent. Math. \vol 1 \yr 1966 \pages 355--374
\endref

\ref \key Wa70  \by C. T. C. Wall
\book Surgery on compact manifolds
 \yr 1970 \publ Academic Press \publaddr London
\endref

\ref \key We67 \by C.~Weber
\paper Plongements de poly\`edres dans le domain metastable
\jour Comment. Math. Helv. \vol 42 \yr 1967 \pages 1--27
\endref

\ref \key We68 \by C.~Weber
\paper Deux remarques sur les plongements d'un AR dans un \'espace euclidien
\jour Bull. Acad. Polon. Sci. Ser. Sci. Math. Astronom. Phys
\vol 16 \yr 1968 \pages 851--855
\endref

\ref \key We \by M. Weiss
\paper Second and third layers in the calculus of embeddings
\jour  \vol  \yr \pages preprint
\endref

\ref \key Wh35 \by H.~Whitney
\paper Differentiable manifolds in Euclidean space
\jour Proc. Nat. Acad. Sci. USA \vol 21:7 \yr 1935 \pages 462--464
\endref

\ref \key Wh44 \by H.~Whitney
\paper The self-intersections of a smooth $n$-manifolds in $2n$-space
\jour Ann. of Math (2) \vol 45 \yr 1944 \pages 220--246
\endref

\ref \key Wu58 \by W.~T.~Wu
\paper On the realization of complexes in a Euclidean space, I
\jour Sci Sinica \vol 7\yr 1958 \pages 251--297
\moreref \paper II
\jour Sci Sinica \vol 7\yr 1958 \pages 365--387
\moreref \paper III
\jour Sci Sinica \vol 8\yr 1959 \pages 133--150
\endref

\ref \key Wu59 \by W.~T.~Wu
\paper \jour Science Record, N.S. \vol 3 \yr 1959 \pages 342--351
\endref

\ref \key Wu65 \by W.~T.~Wu \book A Theory of Embedding, Immersion and
Isotopy of Polytopes in an Euclidean Space
\yr 1965 \publ Science Press \publaddr Peking
\endref

\ref \key Ya54 \by C. T. Yang \pages 262--282
\paper On theorems of Borsuk-Ulam, Kakutani-Yamabe-Yujobo and Dyson, I
\yr 1954 \vol 60 \jour Ann. of Math.
\endref

\ref \key Ze60 \by E.~C.~Zeeman
\paper Unknotting spheres
\jour Ann. of Math. (2) \vol 72 \yr 1960 \pages 350--360
\endref

\ref \key Ze62  \by E.~C.~Zeeman
\paper Isotopies and knots in manifolds
\jour In: Topology of 3-Manifolds, Prentice-Hall, Englewood Cliffs, Ed.
M.~K.~Fort, N.J. \yr 1962
\endref

\ref \key Ze63 \by E.~C.~Zeeman
\paper Unknotting combinatorial balls
\jour Ann. of Math.  \vol 78 \yr 1963 \pages 501--526
\endref

\ref \key Zh75 \by A. V. Zhubr \pages 839--856
\paper A classification of simply-connected spin 6-manifolds (in Russian)
\yr 1975 \vol 39:4 \jour Izvestiya AN SSSR \endref

%\ref \key Zh80 \by A. V. Zhubr \pages 1312--1315
%\paper A classification of simply-connected 6-manifolds (in Russian)
%\yr 1980 \vol 255:6 \jour Doklady AN SSSR \endref

\ref \key Zh89 \by A. V. Zhubr \pages 325--339
\paper Classification of simply-connected topological 6-manifolds
\yr 1989 \vol 1346 \jour Lecture Notes in Math. \endref

\ref \key Zh94 \by L.~Zhongmou
\paper Every 3-manifold with boundary embeds in
$\text{Triod}\times\text{Triod}\times I$
\jour Proc. AMS \vol 122:2 \yr 1994 \pages 575--579
\endref

\endRefs
\enddocument
\end